\RequirePackage{amsmath} 
\documentclass[12pt]{iopart}

\usepackage{iopams}  

\usepackage{subfigure}
\usepackage{graphicx}
\usepackage{color}      
\usepackage{tikz}
\usetikzlibrary{plotmarks}

\newcommand{\s}{2} 
\newcommand{\f}{0.3} 
\newcommand{\n}{1} 
\newcommand{\p}{(1-sqrt(3/5))/2}  

\begin{document}
\title[Hot spots conjecture]{The hot spots conjecture can be false: Some numerical examples}

\author{Andreas Kleefeld$^1$\footnote{Author to whom any correspondence should be addressed.}}

\address{$^1$ Forschungszentrum J\"ulich GmbH, 
	J\"ulich Supercomputing Centre, 52425 J\"ulich, Germany}
\eads{\mailto{a.kleefeld@fz-juelich.de}}
\vspace{10pt}
\begin{indented}
\item[]4 January 2021
\end{indented}

\begin{abstract}
	The hot spots conjecture is only known to be true for special geometries. It can be shown numerically that the hot spots conjecture can fail to be true for easy to construct bounded domains with one hole.
	The underlying eigenvalue problem for the Laplace equation with Neumann boundary condition is solved with boundary integral equations yielding a non-linear eigenvalue problem. 
	Its discretization via the boundary element collocation method in combination with the algorithm by Beyn yields highly accurate results both for the first non-zero eigenvalue and 
	its corresponding eigenfunction which is due to superconvergence. Additionally, it can be shown numerically that the ratio between the maximal/minimal value inside the domain and 
	its maximal/minimal value on the boundary can be larger than $1+10^{-3}$. Finally, numerical examples for easy to construct domains with up to five holes are provided which fail the hot spots conjecture as well.
\end{abstract}

\ams{35J25, 35P20, 65F15, 65M38, 78A46}
%
\vspace{2pc}
\noindent{\it Keywords}: interior Neumann eigenvalues, Helmholtz equation, potential theory, boundary integral equations, numerics
%
\vspace{2pc}\\
\noindent \submitto{arXiv}
%
%
%
\section{Introduction}\label{intro}
The hot spots conjecture has been given in 1974 by Jeffrey Rauch \cite{rauch1974} and explicitly stated a decade later in Kawohl \cite{kawohl1985}. Refer also to the paper by Ba{\~n}uelos \& Burdzy \cite{Bauelos1999OnT} from 1999. 
Since then a lot of researchers have worked on this challenging problem (see Judge \& Mondal \cite{judge} for a recent overview from 2020). 

Before stating it in mathematical terms, we explain it in a simple fashion.
Imagine that we have a flat piece of metal $D$ where 
$D$ is a bounded subset of the two-dimensional space (Euclidean domain) which can have holes with a sufficiently smooth boundary. 
Next, an (almost) arbitrary initial temperature distribution is provided on $D$ (refer to \cite[p. 2]{Bauelos1999OnT}). 
Assume that the domain is insulated, then the hottest and coldest 
spot of $D$ will appear on the boundary when waiting for a long time.

\begin{figure}[!ht]  
\subfigure[Solution at time $t_1=1/200$]{
\includegraphics[width=0.22\textwidth]{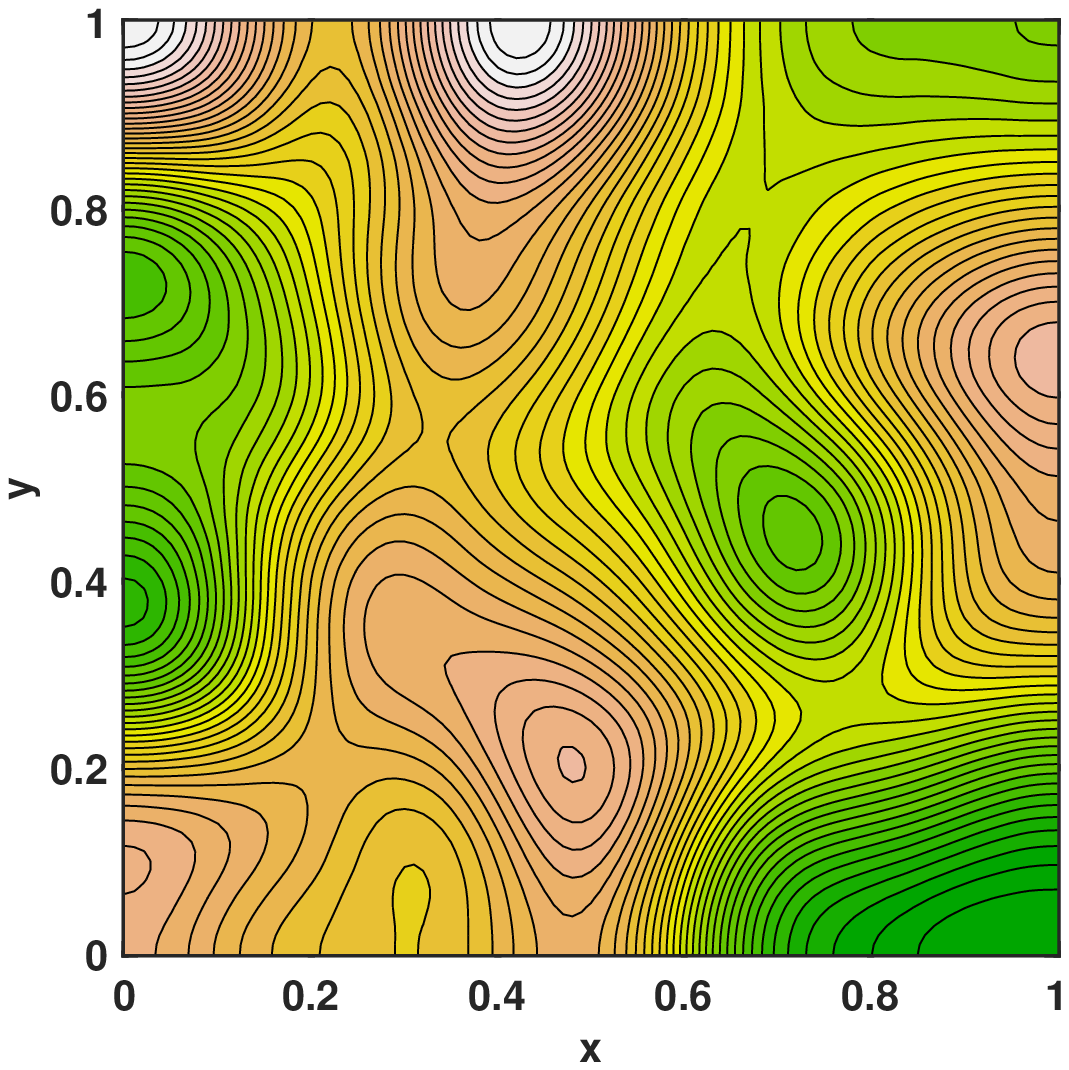}
}
\subfigure[Solution at time $t_2=1/10$]{
\includegraphics[width=0.22\textwidth]{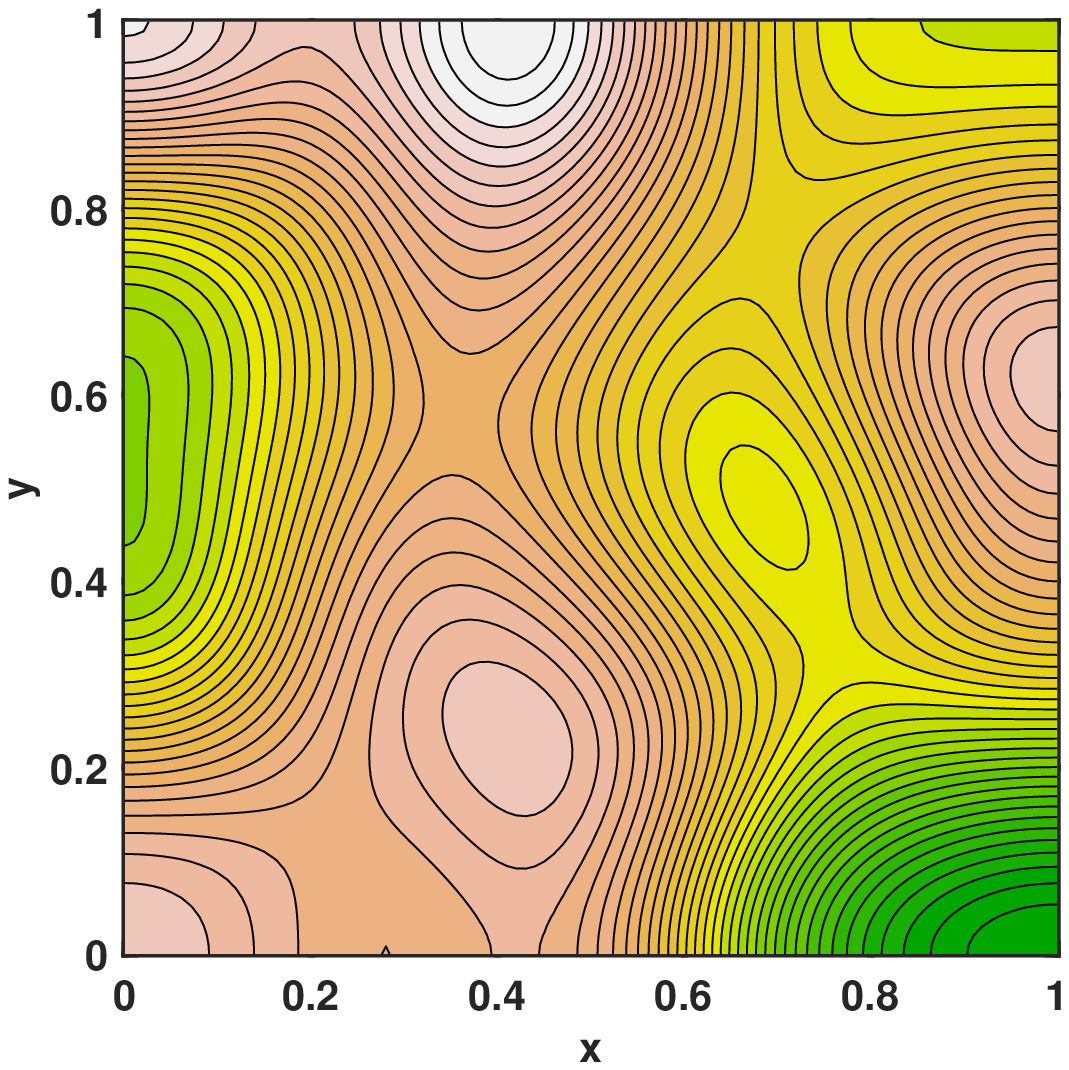}
}  
\subfigure[Solution at time $t_3=1/2$]{
\includegraphics[width=0.22\textwidth]{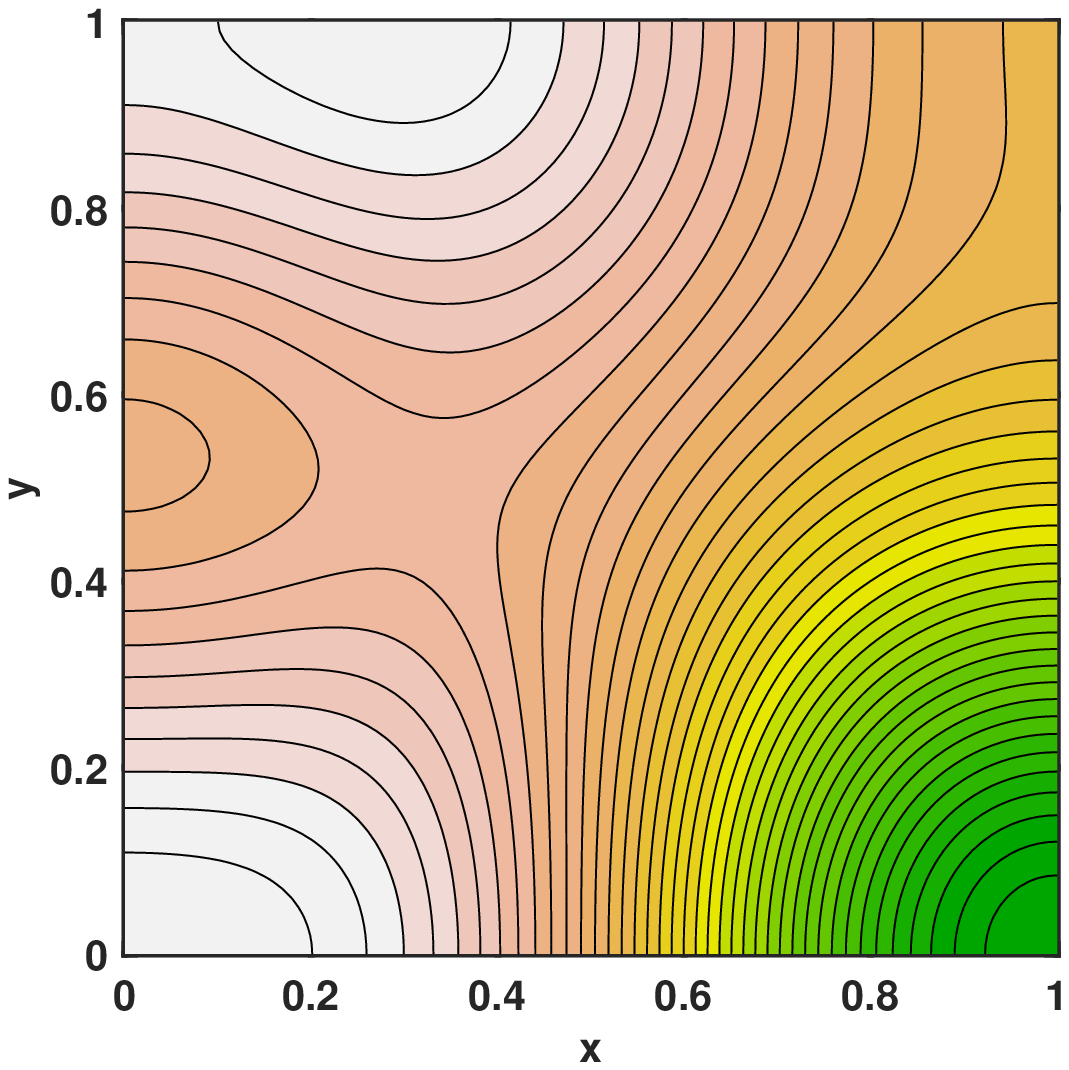}
}
\subfigure[Solution at time $t_4=2$]{
\includegraphics[width=0.22\textwidth]{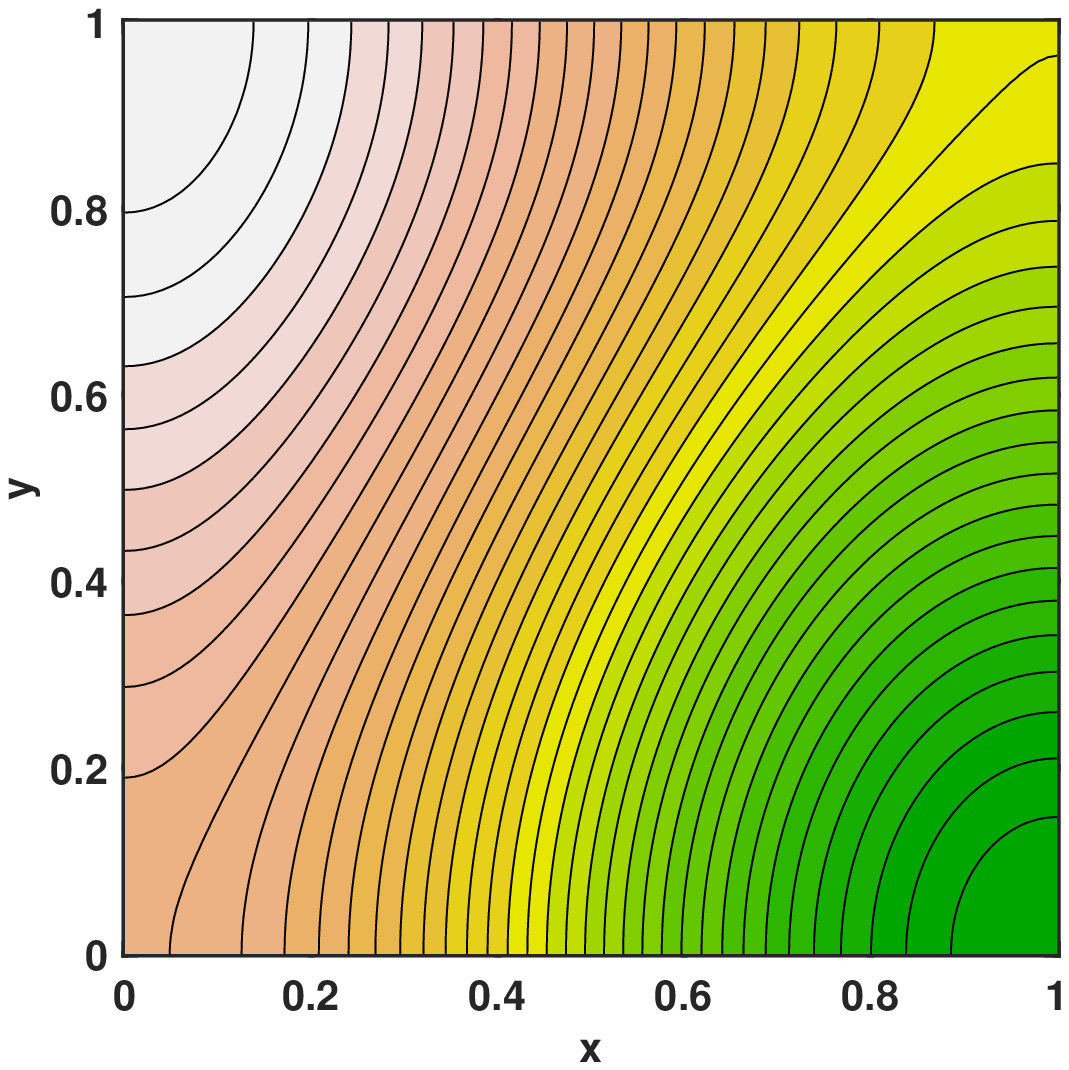}
}
\caption{\label{introex} The numerical solution of the heat equation $\partial_t u=\frac{1}{10}\Delta u$ and 
initial condition with standard normal random numbers and homogeneous Neumann boundary condition for time
$t_1=1/200$, $t_2=1/10$, $t_3=1/2$, and $t_4=2$ using $h=1/100$ and $k=1/100$ 
(see also \cite[p. 11]{asante} for more details on the implementation of the exponential time differencing method). The maximal and minimal value for $T_4$ appear on the boundary. 
Note that the solution $T_4$ is approximately representing the first non-zero Neumann eigenfunction (see also Figure \ref{sqtr}).
}
\end{figure}

Now, we go into the mathematical detail:
That means we have to solve the heat equation $\partial_t u=\Delta u$ for $t\rightarrow \infty$ with homogeneous 
Neumann boundary condition $\partial_{\nu}u=0$ and `almost' arbitrary initial condition in an open connected bounded $D$ with Lipschitz boundary for its equilibrium (see \cite[p. 2]{Bauelos1999OnT} and \cite{sauter} for the definition of a Lipschitz domain). 
Refer to Figure \ref{introex} for an example.

Precisely, we have to find the smallest non-trivial eigenvalue of the Laplacian
with homogeneous Neumann boundary condition and the corresponding eigenfunction. 
Note that the smallest eigenvalue of the Laplacian with homogeneous Neumann boundary condition 
is zero with corresponding eigenfunction $u_0=\mathrm{const}$.
Mathematically, we have to find a solution $u\neq 0$ and the smallest $k\in \mathbb{R}_{>0}$ such that the Helmholtz 
equation $\Delta u+k^2u=0$ in $D$ with $\partial_{\nu}u=0$ 
on the boundary $\Gamma$. All solutions $k$ are called non-trivial interior Neumann eigenvalues 
and $\lambda_i=k^2_i$ will be the $i$-th non-trivial Neumann eigenvalue of the Laplacian. Its corresponding eigenfunctions 
are denoted by $u_i$. Further, it is known that the eigenvalues satisfy 
$0=\lambda_0<\lambda_1\leq \lambda_2\leq\lambda_3\leq \ldots$ when $D$ is a bounded planar domain with Lipschitz boundary (see for example \cite[p. 449]{hempel} and the references therein, specifically \cite{reed}). If $\lambda_1<\lambda_2$ and $\left\langle u(0,x),u_1\right\rangle\neq 0$, then 
$u(t,x)=\mathrm{e}^{\lambda_1 t}\left\langle u(0,x),u_1\right\rangle+$ lower terms.
Note that the first non-trivial eigenvalue can have multiplicity more than one which means that there can be more than one eigenfunction.

Now, the conjecture can be stated as (refer also to \cite[p. 2]{Bauelos1999OnT}): Let $D\subset \mathbb{R}^2$ be an open connected bounded domain with Lipschitz boundary $\Gamma$. Then:
\begin{itemize}
 \item[C1:] For each eigenfunction $u_2(x)$ corresponding to $\lambda_2$ which is not identically zero, we have
 $$\inf_{x\in \Gamma}u_2(x)<u_2(y)<\sup_{x\in \Gamma} u_2(x)\quad \forall y\in D\,.$$
 \item[C2:] For each eigenfunction $u_2(x)$ corresponding to $\lambda_2$ which is not identically zero, we have
 $$\inf_{x\in \Gamma}u_2(x)\leq u_2(y)\leq \sup_{x\in \Gamma} u_2(x)\quad \forall y\in D\,.$$
 \item[C3:] There exist an eigenfunction $u_2(x)$ corresponding to $\lambda_2$ which is not identically zero, such that
 $$\inf_{x\in \Gamma}u_2(x)\leq u_2(y)\leq \sup_{x\in \Gamma} u_2(x)\quad \forall y\in D\,.$$
\end{itemize}

Here, C1 is the original conjecture of Rauch.
The hypothesis has been shown to be true for some special geometries such as parallelepipeds, balls, rectangles, cylinders \cite{kawohl1985}, obtuse triangles \cite{Bauelos1999OnT}, some convex and 
non-convex domains with symmetry \cite{Bauelos1999OnT}, wedges \cite{atar}, lip domains \cite{atarburdzy}, convex domains with two axes of symmetry \cite{jerison}, convex $C^{1,\alpha}$ domains ($0<\alpha<1$) 
with one axis of symmetry \cite{pascu}, a certain class of planar convex domains \cite{miya2}, subequilateral isosceles triangles \cite{miyamoto}, a certain class of acute triangles \cite{siudeja}, Euclidean triangles \cite{judge}, and strips on two-dimensional 
Riemannian manifolds \cite{david}.

It is assumed that the hot spots conjecture is true for arbitrary convex domains, but a proof is still open.
The hot spots conjecture is assumed to be true also for simply-connected bounded non-convex domains, but no successful attempts 
(neither theoretically nor numerically) have been made to prove this conjecture or to find a counterexample.

It has been shown that for some domains with one hole that the hot spots conjecture is true (for example an annulus \cite{kawohl1985}), but that there are also domains 
with one or more holes where the hot spots conjecture is false (see Burdzy \cite{burdzy}, Burdzy \& Werner \cite{burdzywerner}, and Bass \& Burdzy \cite{bass2000}, respectively). 
For domains on manifolds, we refer the reader to \cite{freitas}.

However, the proofs in \cite{burdzy,burdzywerner,bass2000} are very technical and are based on stochastic arguments. 
No numerical results support their counterexamples, 
since their domains are too complicated for being constructed. To be precise, they are very thin and have a polygonal structure. Further, the first non-trivial Neumann eigenvalue is assumed to be simple. If this is not the case, the proof collapses.

The only non-published numerical results given so far are for triangles in the PolyMath project 7 `Hot spots conjecture' from 2012 to 2013 using the finite element method.

Further work in the direction of better understanding the conjecture is given for example by Steinerberger \cite{steinerberger}. Related results on graphs are given by Lederman \& Steinerberger \cite{lederman2019}.

\section*{Contribution}
It is the goal of this paper to construct `simple' domains with one hole and show numerically with high precision (due to superconvergence) that those domains do not 
satisfy the hot spots conjecture. The method based on boundary integral equations is very efficient and its convergence is faster than expected.
Additionally, we show the influence on the location of the hot spots by changing the boundary of the domain
in order to understand this connection. It is believed that this might help researchers to provide 
assumptions when the hot spots conjecture will be true or false for arbitrary 
bounded simply-connected domains which are not necessarily convex. We show that it is possible to construct domains with one hole such that the ratio between the maximum/minimum in the interior 
and its maximum/minimum on the boundary is larger than $1+10^{-3}$. Finally, numerical results are given that show that there exist domains with up to five holes which do not satisfy the hot spots conjecture as well.
The Matlab programs including the produced data are available at github \texttt{https://github.com/kleefeld80/hotspots} and can be used by any researcher trying their own geometries and to reproduce the numerical results within this article. 

\section*{Outline of the paper}
In Section \ref{algo}, we explain the algorithm in order to compute the first non-zero Neumann Laplace eigenvalue and its corresponding eigenfunction for an arbitrary domain with or without a hole using 
boundary integral equations resulting in a non-linear eigenvalue problem. Further, it is shown in detail how to discretize the boundary integral equations via the boundary element collocation method and how to 
numerically solve the non-linear eigenvalue problem. Extensive numerical results are provided in Section \ref{numa} showing the superconvergence and highly accurate results for 
domains with one or no hole. Domains are provided that show the failure of the hot spots conjecture and further interesting results. The extension to domains with up to five holes is straightforward and given at the end 
of this section as well. A short summary and outlook is given in Section \ref{suma}.

\section{The algorithm}\label{algo}
In this section, we explain the algorithm to compute numerically non-trivial interior Neumann eigenvalues and its corresponding eigenfunction to high accuracy for bounded domains with one hole (the extension to more than one hole is straightforward) very efficiently. The ingredients are boundary 
integral equations and its approximation via boundary element collocation method; that is, a two-dimensional problem is reduced to a one-dimensional problem. The resulting non-linear eigenvalue problem is solved using complex-valued contour integrals integrating over the resolvent reducing
the non-linear eigenvalue problem to a linear eigenvalue problem which is possible due to Keldysh's theorem (see Beyn \cite{beyn}).

\subsection{Notations}
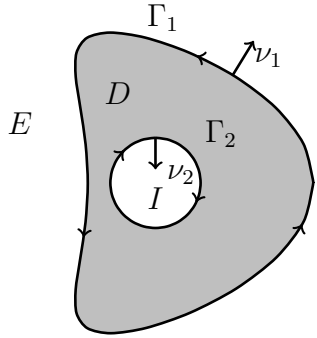
\begin{figure}[!ht]
\centering
\begin{tikzpicture}
    \draw[line width=1pt, scale=\s, domain=0:6.29, smooth, variable=\x, black, fill=lightgray] plot ({0.75*cos(\x r)+0.3*cos(2*\x r)}, {sin(\x r)});
    \node[scale=\n] at (-0.25*\s,0.6*\s) {$D$};
    \node[scale=\n] at (-0.9*\s,0.4*\s) {$E$};
    \draw[line width=1pt, scale=\s, ->] ({0.75*cos(1 r)+0.3*cos(2*1 r)},{sin(1 r)}) -- + (-0.003,0.001);
    \draw[line width=1pt, scale=\s, ->] ({0.75*cos(3.5 r)+0.3*cos(2*3.5 r)},{sin(3.5 r)}) -- +(-0.001,-0.008);
    \draw[line width=1pt, scale=\s, ->] ({0.75*cos(6 r)+0.3*cos(2*6 r)},{sin(6 r)}) -- +(0.001,0.002);
    \draw[line width=1pt, scale=\s, draw=none] ({0.75*cos(1.2 r)+0.3*cos(2*1.2 r)},{sin(1.2 r)}) -- node[scale=\n, above] {$\Gamma_1$} + (0,0);
    \draw[line width=1pt, scale=\s, ->] ({0.75*cos(0.8 r)+0.3*cos(2*0.8 r)},{sin(0.8 r)}) -- node[scale=\n, right] {$\nu_1$} + ({0.2*(cos(0.8 r))},{-0.2*(-0.75*sin(0.8 r)-0.6*sin(2*0.8 r))});
    \draw[line width=1pt, scale=\f*\s, domain=0:6.284, smooth, variable=\t, black, fill=white] plot ({sin(\t r)}, {cos(\t r)});
    \node[line width=1pt, scale=\n] at (0,-0.2) {$I$};
    \draw[line width=1pt, scale=\f*\s, ->] ({sin(2 r)},{cos(2 r)}) -- +(-0.001,-0.003);
    \draw[line width=1pt, scale=\f*\s, ->] ({sin(5.5 r)},{cos(5.5 r)}) -- +(0.001,0.001);
    \draw[line width=1pt, scale=\f*\s, draw=none] ({sin(1 r)},{cos(1 r)}) -- node[scale=\n, above right] {$\Gamma_2$} + (0,0);
    \draw[line width=1pt, scale=\f*\s, ->] ({sin(0 r)},{cos(0 r)}) -- node[scale=\n, below right] {$\nu_2$} + ({0.7*sin(0 r)},{0.7*-cos(0 r)});
\end{tikzpicture}
\caption{\label{setup} Used notations for a bounded domain with one hole.}
\end{figure}

We consider a bounded Lipschitz domain $D\subset \mathbb{R}^2$ with one hole. The outer boundary $\Gamma_1$ is assumed to be sufficiently smooth that is oriented counter-clockwise
and a sufficiently smooth inner boundary $\Gamma_2$ that is oriented clockwise. The normal $\nu_1$ on the boundary $\Gamma_1$ 
is pointing into the unbounded exterior $E$. The normal $\nu_2$ on the boundary $\Gamma_2$ is pointing into the bounded exterior $I$. We refer the reader to Figure \ref{setup}. 
The boundary of $D$ is given by $\Gamma=\Gamma_1\cup\Gamma_2$ 
($\Gamma_1 \cap \Gamma_2\neq \emptyset$).
Note that we also consider bounded domains without a hole. In this case, we have $\Gamma_2=\emptyset$ and hence $\Gamma=\Gamma_1$ and $I=\emptyset$.

\subsection{Boundary integral equation}
The solution to the Helmholtz equation $\Delta u+k^2 u=0$ (reduced wave equation) in the domain $D$ for a given wave number $k$ with $\mathrm{Im}(k)\geq 0$ is given by (see \cite[Theorem 2.1]{coltonkress})
\begin{eqnarray*}
u(x)=\int_{\Gamma}\partial_{\nu(y)} u(y)\cdotp \Phi_k(x,y)-u(y)\cdotp \partial_{\nu(y)}\Phi_k(x,y)\,\mathrm{d}s(y)\,,\quad x\in D
\end{eqnarray*}
which can be written in our notation as
\begin{eqnarray}
u(x)&=&\int_{\Gamma_1}\partial_{\nu_1(y)} u(y)\cdotp \Phi_k(x,y)-u(y)\cdotp \partial_{\nu_1(y)}\Phi_k(x,y)\,\mathrm{d}s(y)\nonumber\\
&+&\int_{\Gamma_2}\partial_{\nu_2(y)} u(y)\cdotp \Phi_k(x,y)-u(y)\cdotp \partial_{\nu_2(y)}\Phi_k(x,y)\,\mathrm{d}s(y)\,,\quad x\in D
\label{greensrep}
\end{eqnarray}
where $\Phi_k(x,y)=\mathrm{i}H_0^{(1)}(k\|x-y\|)/4$, $x\neq y$ denotes the fundamental solution of the Helmholtz equation in two dimensions (see \cite[p. 66]{coltonkress}). 
Here, $H_0^{(1)}$ denotes the first-kind Hankel function of order zero.
We denote $u(y)$ for $y\in\Gamma_1$ as $u_1(y)$ and similarly $u(y)$ for $y\in\Gamma_2$ as $u_2(y)$. Hence, we can write (\ref{greensrep}) as 
\begin{eqnarray}
u(x)&=&-\int_{\Gamma_1}u_1(y)\cdotp \partial_{\nu_1(y)}\Phi_k(x,y)\,\mathrm{d}s(y)\nonumber\\
    &&-\int_{\Gamma_2}u_2(y)\cdotp \partial_{\nu_2(y)}\Phi_k(x,y)\,\mathrm{d}s(y)\,,\quad x\in D
\label{greensrep2}
\end{eqnarray}
where we also used the homogeneous Neumann boundary conditions $\partial_{\nu_1} u=0$ and $\partial_{\nu_2} u=0$.
We rewrite (\ref{greensrep2}) as  
\begin{eqnarray}
u(x)=-\mathrm{DL}_k^{\Gamma_1}u_1(x)-\mathrm{DL}_k^{\Gamma_2}u_2(x)\,,\quad x\in D
\label{start}
\end{eqnarray}
where we used the notation
\[\mathrm{DL}_k^{\Gamma_i} \psi_i(x)=\int_{\Gamma_i}\psi_i(y)\cdotp \partial_{\nu_i(y)}\Phi_k(x,y)\,\mathrm{d}s(y)\,,\quad x\in D\,,\quad i=1,2\]
for the \textit{acoustic double layer potential} with density $\psi_i$ (see \cite[p. 39]{coltonkress}). 
Assume for a moment that $k$ is given. The functions $u_1$ and $u_2$ are still unknown. Once we know them, we can compute the solution $u$ inside 
the domain of $D$ at any point we want using (\ref{start}).

Now, we explain how to obtain those functions $u_1$ and $u_2$ on the boundary.
Letting $x\in D$ approach the boundary $\Gamma_1$ and using the jump relation of the acoustic double layer operator (see \cite[p. 39]{coltonkress} for the smooth boundary case, otherwise \cite{sauter}), yields the 
boundary integral equation
\begin{eqnarray}
u_1(x)=-\left(\mathrm{D}_k^{\Gamma_1\rightarrow \Gamma_1}u_1(x)-\left(1-\Omega_1(x)\right)u_1(x)\right)-\mathrm{D}_k^{\Gamma_2\rightarrow \Gamma_1}u_2(x)\,,\quad x\in \Gamma_1\nonumber\\
\label{biestart}
\end{eqnarray}
where we used the notation 
\[\mathrm{D}_k^{\Gamma_i\rightarrow \Gamma_j} \psi_i(x)=\int_{\Gamma_i}\psi_i(y)\cdotp \partial_{\nu_i(y)}\Phi_k(x,y)\,\mathrm{d}s(y)\,,\quad x\in \Gamma_j\,,\quad i,j=1,2\]
for the \textit{double layer operator} (see \cite[p. 41]{coltonkress}). Here, $\Omega_1(x)$ denotes the interior solid angle at a point $x$ on $\Gamma_1$. When the boundary is smooth at this point, then $\Omega_1(x)=1/2$. In fact, 
$\Omega_1(x)$ is $1/2$ almost everywhere for Lipschitz domains.

Similarly, we obtain for $x\in D$ approaching the boundary $\Gamma_2$ and using the jump relation for the double layer operator
\begin{eqnarray}
u_2(x)=-\mathrm{D}_k^{\Gamma_1\rightarrow \Gamma_2}u_1(x)-\left(\mathrm{D}_k^{\Gamma_2\rightarrow \Gamma_2}u_2(x)-\left(1-\Omega_2(x)\right)u_2(x)\right)\,,\quad x\in \Gamma_2\,.\nonumber\\
\label{biestart2}
\end{eqnarray}
We can rewrite (\ref{biestart}) and (\ref{biestart2}) as a $2\times 2$ system of boundary integral equations in the form
\begin{eqnarray}
\underbrace{\left(
\overbrace{
\left(\begin{matrix}
       \Omega_1 & 0\\
       0 & \Omega_2
      \end{matrix}
\right)}^{\mathrm{C}}
\overbrace{\left(\begin{matrix}
       \mathrm{I} & 0\\
       0 & \mathrm{I}
      \end{matrix}
\right)}^{\mathrm{I}_\mathrm{B}}+
\overbrace{\left(\begin{matrix}
        \mathrm{D}_k^{\Gamma_1\rightarrow \Gamma_1}& \mathrm{D}_k^{\Gamma_2\rightarrow \Gamma_1}\\
        \mathrm{D}_k^{\Gamma_1\rightarrow \Gamma_2}&\mathrm{D}_k^{\Gamma_2\rightarrow \Gamma_2} 
      \end{matrix}
\right)}^{\mathrm{K}(k)}\right)}_{\mathrm{M}(k)}\underbrace{\left(\begin{matrix}
       u_1\\
       u_2
      \end{matrix}\right)}_{u}
=\left(\begin{matrix}
       0\\
       0
      \end{matrix}
\right)\quad\text{on }\Gamma
 \label{nonlinear}
\end{eqnarray}
where $\mathrm{I}$ and $\mathrm{I}_\mathrm{B}$ denotes the identity and the $2\times 2$ block identity operator, respectively.
Hence, we have to numerically solve the non-linear eigenvalue problem (\ref{nonlinear}) written as $\mathrm{M}(k)u=0$
to find the smallest non-trivial (real) eigenvalue $k$ and the corresponding eigenfunction $u$. Then, we can numerically evaluate (\ref{start}) to compute the eigenfunction at any point in the interior we want.
As in \cite[p. 188]{kleefeld2019shape}, we can argue that the compact operator $\mathrm{K}(k)$ maps from $\mathcal{H}^{-1/2}(\Gamma_1)\times \mathcal{H}^{-1/2}(\Gamma_2)$ to $\mathcal{H}^{1/2}(\Gamma_1)\times \mathcal{H}^{1/2}(\Gamma_2)$. Here, $\mathcal{H}^{s}(\Gamma)$ denotes a Sobolev space of order $s\in\mathbb{R}$ on the domain $\Gamma$ which are defined via Bessel potentials (see \cite[pp. 75--76]{McLea2000} for more details). The operator 
$\mathrm{M}(k)=\mathrm{C}\cdotp\mathrm{I}_{\mathrm{B}}+\mathrm{K}(k)$ is Fredholm of index zero for $k\in \mathbb{C}\backslash \mathbb{R}_{\leq 0}$ and therefore the theory of eigenvalue problems
for holomorphic Fredholm operator-valued functions applies to $\mathrm{M}(k)$.
\subsection{Discretization}
\begin{figure}[!ht]
	\begin{center}
		\begin{tikzpicture}
		\draw[line width=1pt,scale=\s, domain=0:6.284, smooth, variable=\t, black, fill=white] plot ({sin(\t r)}, {cos(\t r)});
		\draw[line width=1pt] (\s*1-0.2,0) -- (\s*1+0.2,0);
		\draw[line width=1pt] (0,\s*1-0.2) -- (0,\s*1+0.2);
		\draw[line width=1pt] (\s*-1-0.2,0) -- (\s*-1+0.2,0);
		\draw[line width=1pt] (0,\s*-1-0.2) -- (0,\s*-1+0.2);
		\node[scale=\n] at ( 0.9*\s,0.9*\s) {$\Delta_1$};
		\node[scale=\n] at (-0.9*\s,0.9*\s) {$\Delta_2$};
		\node[scale=\n] at (-0.9*\s,-0.9*\s) {$\Delta_3$};
		\node[scale=\n] at ( 0.9*\s,-0.9*\s) {$\Delta_4$};
		\end{tikzpicture}
	\end{center}
	\caption{\label{Rand}Subdivision of the given boundary $\Gamma_1$ into four pieces $\Delta_1$, $\Delta_2$, $\Delta_3$, and $\Delta_4$.} 
\end{figure}
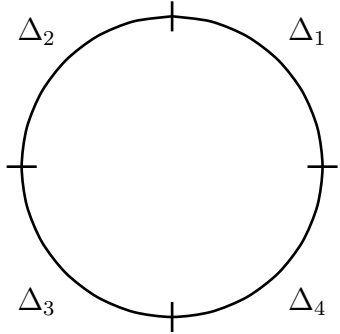

In this section, we explain how to discretize (\ref{nonlinear}) using quadratic interpolation of the boundary, 
but using piecewise quadratic interpolation with $\alpha=(1-\sqrt{3/5})/2$ (see \cite{kleefeldlin2} for the 3D case) instead of quadratic interpolation for the unknown $u$ on each of the $n_f$ boundary elements which ultimately leads to the non-linear eigenvalue 
\begin{eqnarray}
\mathbf{M}(k)\vec{u}=\vec{0}
\label{nonlin}
\end{eqnarray}
where the matrix is of size $3\cdotp 2\cdotp n_f\times 3\cdotp 2\cdotp n_f$. The size of the matrix is slightly larger than the one given in Kleefeld \cite{kleefeld2019shape}, but it has the advantage that no singular 
integral has to be evaluated numerically, since we can use a similar singularity subtraction technique as explained in Kleefeld and Lin \cite[pp. A1720--A1721]{kleefeldlin} and the convergence rate is slightly higher.
The details are about to follow for a domain without a hole for simplicity. In this case, we have to solve a boundary integral equation of the second kind of the form
\begin{eqnarray}
\Omega_1(x)u_1(x)+\int_{\Gamma_1}u_1(y)\partial_{\nu_1(y)}\Phi_k(x,y)\,\mathrm{d}s(y)=0\,,\quad x\in \Gamma_1\,.
\label{jetztgehtslos}
\end{eqnarray}
First, we subdivide the boundary $\Gamma_1$ into $n_f$ pieces denoted by $\Delta_j$ with $j=1,\ldots,n_f$. A subdivision into four pieces is shown in Figure \ref{Rand} for the unit circle.

\begin{figure}[!ht]
	\begin{center}
		\begin{tikzpicture}
	\draw[line width=1pt,scale=\s, domain=0:2*pi, smooth, variable=\t, black, fill=white] plot ({sin(\t r)}, {cos(\t r)});
	\draw[line width=1pt] (\s*1-0.1,-0.1) -- (\s*1+0.1,0.1);
	\draw[line width=1pt] (\s*1-0.1,0.1) -- (\s*1+0.1,-0.1);
	\draw[line width=1pt] (-0.1,\s*1-0.1) -- (0.1,\s*1+0.1);
	\draw[line width=1pt] (0.1,\s*1-0.1) -- (-0.1,\s*1+0.1);
	\draw[line width=1pt] (\s*-1-0.1,-0.1) -- (\s*-1+0.1,0.1);
	\draw[line width=1pt] (\s*-1-0.1,0.1) -- (\s*-1+0.1,-0.1);
	\draw[line width=1pt] (-0.1,\s*-1-0.1) -- (0.1,\s*-1+0.1);
	\draw[line width=1pt] (0.1,\s*-1-0.1) -- (-0.1,\s*-1+0.1);
	\draw[line width=1pt] ({\s*sin(pi/4 r)},{\s*cos(pi/4 r)}) circle (3pt);
	\draw[line width=1pt] ({\s*sin(3*pi/4 r)},{\s*cos(3*pi/4 r)}) circle (3pt);
	\draw[line width=1pt] ({\s*sin(5*pi/4 r)},{\s*cos(5*pi/4 r)}) circle (3pt);
	\draw[line width=1pt] ({\s*sin(7*pi/4 r)},{\s*cos(7*pi/4 r)}) circle (3pt);
	\node[scale=\n] at ({\s*1.2*sin(1*pi/4 r)},{\s*1.2*cos(1*pi/4 r)}) {$v_2$};
	\node[scale=\n] at ({\s*1.2*sin(3*pi/4 r)},{\s*1.2*cos(3*pi/4 r)}) {$v_8$};
	\node[scale=\n] at ({\s*1.2*sin(5*pi/4 r)},{\s*1.2*cos(5*pi/4 r)}) {$v_6$};
	\node[scale=\n] at ({\s*1.2*sin(7*pi/4 r)},{\s*1.2*cos(7*pi/4 r)}) {$v_4$};
	\node[scale=\n] at ({\s*1.2*sin(2*pi/4 r)},{\s*1.2*cos(2*pi/4 r)}) {$v_1$};
	\node[scale=\n] at ({\s*1.2*sin(4*pi/4 r)},{\s*1.2*cos(4*pi/4 r)}) {$v_7$};
	\node[scale=\n] at ({\s*1.2*sin(6*pi/4 r)},{\s*1.2*cos(6*pi/4 r)}) {$v_5$};
	\node[scale=\n] at ({\s*1.2*sin(8*pi/4 r)},{\s*1.2*cos(8*pi/4 r)}) {$v_3$};
		\end{tikzpicture}
	\end{center}
	\caption{\label{Rand2}The eight nodes $v_1,\ldots,v_8$, the four vertices $v_1$, $v_3$, $v_5$, and $v_7$ (marked with $\times$), and the four midpoints $v_2$, $v_4$, $v_6$, and $v_8$ (marked with $\circ$) for the four faces $\Delta_1$, $\Delta_2$, $\Delta_3$, and $\Delta_4$.} 
\end{figure}
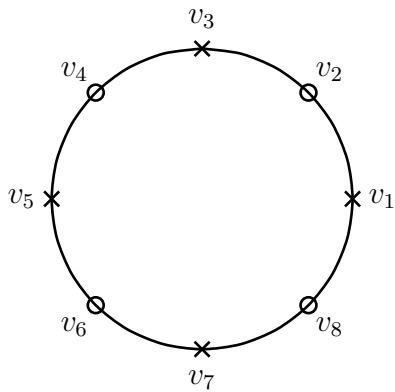

\begin{figure}[!ht]
	\begin{center}
		\begin{tikzpicture}
		\draw[line width=1pt,scale=\s, domain=0:pi/2, smooth, variable=\t, black, fill=white] plot ({sin(\t r)}, {cos(\t r)});
		\draw[line width=1pt,scale=\s, domain=0:1, smooth, dashed, variable=\t, black] plot ({(1-\t)*(2*(1-\t)-1)+4*sqrt(2)/2*\t*(1-\t)},{4*sqrt(2)/2*\t*(1-\t)+\t*(2*\t-1)});
		\draw[line width=1pt] (\s*1-0.1,-0.1) -- (\s*1+0.1,0.1);
		\draw[line width=1pt] (\s*1-0.1,0.1) -- (\s*1+0.1,-0.1);
		\draw[line width=1pt] (-0.1,\s*1-0.1) -- (0.1,\s*1+0.1);
		\draw[line width=1pt] (0.1,\s*1-0.1) -- (-0.1,\s*1+0.1);
		\draw[line width=1pt] ({\s*sin(pi/4 r)},{\s*cos(pi/4 r)}) circle (3pt);
		\node[scale=\n] at ({\s*1.2*sin(1*pi/4 r)},{\s*1.2*cos(1*pi/4 r)}) {$\frac{\pi}{4}$};
		\node[scale=\n] at ({\s*1.2*sin(2*pi/4 r)},{\s*1.2*cos(2*pi/4 r)}) {$0$};
		\node[scale=\n] at ({\s*1.2*sin(8*pi/4 r)},{\s*1.2*cos(8*pi/4 r)}) {$\frac{\pi}{2}$};
		\end{tikzpicture}
	\end{center}
	\caption{\label{Rand3}The approximation of the first part of the boundary $\Delta_1$ (solid line) via a quadratic interpolation polynomial using two vertices and the midpoint (dashed line).} 
\end{figure}
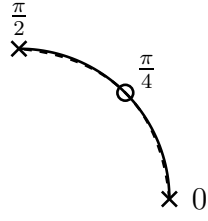

\begin{figure}[!ht]
	\begin{center}
		\begin{tikzpicture}
		\draw[scale=\s, domain=0:pi/2, smooth, variable=\t, black, fill=white] plot ({sin(\t r)}, {cos(\t r)});
		\draw[scale=\s, domain=0:1, smooth, dashed, variable=\t, black] plot ({(1-\t)*(2*(1-\t)-1)+4*sqrt(2)/2*\t*(1-\t)},{4*sqrt(2)/2*\t*(1-\t)+\t*(2*\t-1)});
		\draw[line width=1] (\s*1-0.1,-0.1) -- (\s*1+0.1,0.1);
		\draw[line width=1] (\s*1-0.1,0.1) -- (\s*1+0.1,-0.1);
		\draw[line width=1] (-0.1,\s*1-0.1) -- (0.1,\s*1+0.1);
		\draw[line width=1] (0.1,\s*1-0.1) -- (-0.1,\s*1+0.1);
		\draw[mark=triangle*, mark size=3 pt] plot coordinates {({\s*sin(pi/4 r)},{\s*cos(pi/4 r)})};
		\draw[mark=triangle*, mark size=3 pt] plot coordinates {({\s*((1-\p)*(2*(1-\p)-1)+4*sqrt(2)/2*\p*(1-\p))},{\s*(4*sqrt(2)/2*\p*(1-\p)+\p*(2*\p-1))})};
		\draw[mark=triangle*, mark size=3 pt] plot coordinates {({\s*((\p)*(2*(\p)-1)+4*sqrt(2)/2*(1-\p)*(\p))},{\s*(4*sqrt(2)/2*(1-\p)*(\p)+(1-\p)*(2*(1-\p)-1))})};
		\end{tikzpicture}
	\end{center}
	\caption{\label{Rand4}The first three collocation nodes (solid triangles) on the first part of the approximated boundary (dashed line) for the unit circle using $\alpha=(1-\sqrt{3/5})/2$. The exact boundary is shown with a solid line including the two vertices marked by a cross.} 
\end{figure}
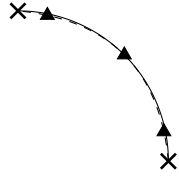

Then equation (\ref{jetztgehtslos}) can be equivalently written as
\begin{eqnarray*}
\Omega_1(x)u_1(x)+\sum_{j=1}^{n_f}\int_{\Delta_j}u_1(y)\partial_{\nu_1(y)}\Phi_k(x,y)\,\mathrm{d}s(y)=0\,,\quad x\in \Gamma_1\,.
\end{eqnarray*}
For each $j$ there exists a unique map $m_j$ which maps from the standard interval $\sigma=[0,1]$ to $\Delta_j$. Then, we can apply a simple change of variables to each integral over $\Delta_j$ giving
\begin{eqnarray*}
&&\Omega_1(x)u_1(x)\nonumber\\
&+&\sum_{j=1}^{n_f}\int_{\sigma}u_1(m_j(s))\partial_{\nu_1(m_j(s))}\Phi_k(x,m_j(s))J(s)\,\mathrm{d}s(s)=0\,,\; x\in \Gamma_1
\end{eqnarray*}
with the Jacobian given by $J(s)=\|\partial_s m_j(s)\|$.
In most cases, we can explicitly write down this map. However, we approximate each $m_j(s)$ by a quadratic interpolation polynomial 
$m_j(s)\approx\widetilde{m}_j(s)=\sum_{i=1}^3 v_{(2i-j) \,\mathrm{mod}\, (2n_f)} L_i(s)$ where the Lagrange basis functions are 
$$L_1(s)=u\cdotp (1-2s)\,,\quad L_2(s)=4s\cdotp u\,,\quad \text{ and }\quad L_3(s)=s\cdotp (2s-1)$$
with $u=1-s$.
Here, $\gamma\, \mathrm{mod}\,\delta=\gamma-\lfloor\frac{\gamma}{\delta}\rfloor\cdotp \delta$ with $\lfloor\cdotp\rfloor$ the floor function. The nodes $v_\ell$ ($\ell=1,\ldots,2n_f$) are the given vertices and midpoints of the $n_f$ faces. We refer the reader to Figure \ref{Rand2} for an example with four faces and eight nodes (four vertices and midpoints, respectively). 
An example how the approximation of $\Delta_1$ via a quadratic interpolation polynomial using two vertices and the midpoint looks like is shown in Figure \ref{Rand3}.
Next, we define the `collocation nodes' $\widetilde{v}_{j,k}$ by $\widetilde{v}_{j,k}=\widetilde{m}_j(q_k)$ for $j=1,\ldots,n_f$ and for $k=1,2,3$ where $q_1=\alpha$, $q_2=1/2$, and $q_3=1-\alpha$ with $0<\alpha<1/2$ 
a given and fixed constant. This ensures that the collocation nodes are always lying within a piece of the boundary and at those points the interior solid angle is $1/2$. For a specific choice of $\alpha$ the overall convergence 
rate can be improved.
The first three collocation nodes on the approximated boundary for the unit circle using $\alpha=(1-\sqrt{3/5})/2$ are shown in Figure \ref{Rand4}. 

The unknown function $u_1(\widetilde{m}_j(s))$ is now approximated on each of the $j$ pieces by a quadratic interpolation polynomial of the form $\sum_{k=1}^{3}u_1(\widetilde{m}_j(q_k))\widetilde{L}_k(s)$ which can be written as $\sum_{k=1}^{3}u_1(\widetilde{v}_{j,k})\widetilde{L}_k(s)$ where
the Lagrange basis functions are given by
\begin{eqnarray*}
	\widetilde{L}_1(s)= \frac{u-\alpha}{1-2\alpha}\frac{1-2s}{1-2\alpha}\,,\;
	\widetilde{L}_2(s)=4\frac{s-\alpha}{1-2\alpha}\frac{u-\alpha}{1-2\alpha}\,,\;
	\widetilde{L}_3(s)= \frac{s-\alpha}{1-2\alpha}\frac{2s-1}{1-2\alpha}
\end{eqnarray*}
with $u=1-s$.
We obtain
\begin{eqnarray*}
	&&\Omega_1(x)u_1(x)\\
	&+&\sum_{j=1}^{n_f}\sum_{k=1}^3\int_{\sigma}\partial_{\nu_1(\widetilde{m}_j(s))}\Phi_k(x,\widetilde{m}_j(s))\|\partial_s \widetilde{m}_j(s)\|\widetilde{L}_k(s)\,\mathrm{d}s(s) u_1(\widetilde{v}_{j,k})=r(x)
	\end{eqnarray*}
with $r(x)$ the residue which is due to the different approximations. We set $r(\widetilde{v}_{i,\ell})=0$ and since $\Omega_1(\widetilde{v}_{i,\ell})=1/2$ always by the choice of the collocation nodes, we obtain the linear system of size $3n_f\times 3n_f$
\begin{eqnarray*}
	\frac{1}{2}u_1(\widetilde{v}_{i,\ell})
	+\sum_{j=1}^{n_f}\sum_{k=1}^3a_{i,\ell,j,k} u_1(\widetilde{v}_{j,k})=0
\end{eqnarray*}
with the resulting integrals
\begin{eqnarray}
a_{i,\ell,j,k}=\int_{\sigma}\partial_{\nu_1(\widetilde{m}_j(s))}\Phi_k(\widetilde{v}_{i,\ell},\widetilde{m}_j(s))\|\partial_s \widetilde{m}_j(s)\|\widetilde{L}_k(s)\,\mathrm{d}s(s)
\label{integrals}
\end{eqnarray}
which will be approximated by the adaptive Gauss-Kronrod quadrature (see \cite{shampine}). This can be written abstractly as $\mathbf{M}(k)\vec{u}=\vec{0}$. Note that the integrand of the integral of the right-hand side (\ref{integrals}) can easily by written down as
$$\frac{\mathrm{i}k H_1^{(1)}(kr)}{4r}\left(a\cdotp n_{1}+b\cdotp n_2\right)\widetilde{L}_k(s)$$
with $ H_1^{(1)}$ the first-kind Hankel function of order one and with
\begin{eqnarray*}
a&=&\left[\widetilde{v}_{i,\ell}-\widetilde{m}_j(s)\right]_1\,,\quad b=\left[\widetilde{v}_{i,\ell}-\widetilde{m}_j(s)\right]_2\,,\quad
r=\sqrt{a^2+b^2}\,,\\
n_{1}&=&\left[\partial_s \widetilde{m}_j(s)\right]_2\,,\quad 
n_{2}=-\left[\partial_s \widetilde{m}_j(s)\right]_1\,.
\end{eqnarray*}
Note that the Jacobian cancels out. A word has to be spent on the following issue: When $P\neq Q$ with $P=\widetilde{v}_{i,\ell}$ and $Q=\widetilde{m}_j(s)$, then the integrand of the integral of the right-hand side (\ref{integrals}) is smooth. 
However, for the case $P=Q$ 
a singularity within the integral of the right-hand side (\ref{integrals}) is present. 
In this
case, we can use the singularity subtraction method to rewrite the singular integral in the following form
\begin{eqnarray*}
& &\int_{\sigma}\partial_{\nu_1(\widetilde{m}_j(s))}\Phi_k(\widetilde{v}_{i,\ell},\widetilde{m}_j(s))\|\partial_s \widetilde{m}_j(s)\|\widetilde{L}_k(s)\,\mathrm{d}s(s)\\
&=&\int_{\sigma}\partial_{\nu_1(\widetilde{m}_j(s))}\left(\Phi_k(\widetilde{v}_{i,\ell},\widetilde{m}_j(s))-\Phi_0(\widetilde{v}_{i,\ell},\widetilde{m}_j(s))\right)\|\partial_s \widetilde{m}_j(s)\|\widetilde{L}_k(s)\,\mathrm{d}s(s)\\
&+&\int_{\sigma}\partial_{\nu_1(\widetilde{m}_j(s))}\Phi_0(\widetilde{v}_{i,\ell},\widetilde{m}_j(s))\|\partial_s \widetilde{m}_j(s)\|\widetilde{L}_k(s)\,\mathrm{d}s(s)=I_{i,\ell}^\text{smooth}+I_{i,\ell}^\text{sing}
\end{eqnarray*}
where $\Phi_0(P,Q)=-\log(|P-Q|)/(2\pi)$ is the fundamental solution of the Laplace equation.
The integral $I_{i,\ell}^\text{smooth}$ has a smooth kernel (no singularity present) and is converging rapidly to zero when increasing the number of faces (independent of the wave number $k$). Hence, we directly set $I_{i,\ell}^\text{smooth}=0$. The integral $I_{i,\ell}^\text{sing}$ (a singularity is present) can be rewritten as a sum of integrals without any singularity.
This is due to the fact that we have $\mathrm{D}_0^{\Gamma_1\rightarrow \Gamma_1}\, 1\,(x)=-\Omega_1(x)$ mit $\psi=1$, $\forall x\in \Gamma_1$ (see \cite[p. 363]{seybert}) and hence, we approximately use for all $i,\ell$:
$$\sum_{j=1}^{n_f}\sum_{k=1}^3\int_{\sigma}\partial_{\nu_1(\widetilde{m}_j(s))}\Phi_0(\widetilde{v}_{i,\ell},\widetilde{m}_j(s))\|\partial_s \widetilde{m}_j(s)\|\widetilde{L}_k(s)\,\mathrm{d}s(s)\approx-\Omega_1(\widetilde{v}_{i,\ell})=-\frac{1}{2}$$
(that is, the row sum of the matrix $\mathbf{M}(0)$ obtained from the discretization of the double layer for the Laplace equation shall be $-1/2$)
and hence, we can compute $I_{i,\ell}^\text{sing}$ as
\begin{eqnarray}
I_{i,\ell}^\text{sing}&\approx&-\underbrace{\Omega_1(\widetilde{v}_{i,\ell})}_{\frac{1}{2}}\nonumber\\
&&-\mathop{\sum_{j=1}^{n_f}}_{j\neq i}\mathop{\sum_{k=1}^3}_{k\neq \ell}\int_{\sigma}\partial_{\nu_1(\widetilde{m}_j(s))}\Phi_0(\widetilde{v}_{i,\ell},\widetilde{m}_j(s))\|\partial_s \widetilde{m}_j(s)\|\widetilde{L}_k(s)\,\mathrm{d}s(s)\,.
\label{integrals2}
\end{eqnarray}
Each of the integrands within the integrals of the right-hand side (\ref{integrals2}) are smooth and can be computed with the previously mentioned Gauss-Kronrod quadrature. 
Note that we never have to compute the interior solid angle since it is always $1/2$ by the choice of the collocation points. In fact, we never have to use the value $1/2$ since it cancels out with the $1/2$ 
within the definition of 
$\mathbf{M}(k)$.

Finally, if we use a domain with a hole, we obtain the system of size $m\times m=3\cdotp 2\cdotp n_f\times 3\cdotp 2\cdotp n_f$ when including the boundary $\Gamma_2$ and the unknown function $u_2$. Written abstractly we obtain the non-linear eigenvalue problem
\begin{eqnarray}
\mathbf{M}(k)\vec{u}=\vec{0}
\label{nonlin2}
\end{eqnarray}

Of course, the extension to more than one hole is obvious. In this case, the matrix $\mathbf{M}(k)$ within (\ref{nonlin2}) will be of size $(q-1)\cdotp 2\cdotp n_f\times (q-1)\cdotp 2\cdotp n_f$ where $q$ denotes the number of 
holes within the domain.

\subsection{Non-linear eigenvalue problem}
The non-linear eigenvalue problem (\ref{nonlin2}) is solved with the Beyn algorithm \cite{beyn}. It is based on complex-valued contour integrals integrating over the resolvent reducing the non-linear eigenvalue to a linear one (of very small size) which can be achieved by Keldysh's theorem. Precisely, a user-specified $2\pi$-periodic contour $\gamma$ of class $C^1$ within the complex plane has to be given. We need a contour that is enclosing a part of the real line where the smallest non-zero eigenvalue is expected. We usually use a circle with radius $R$ and center $(\mu,0)$ (in order to exclude the eigenvalue zero, we choose $\mu>R$). In this case, we have $\varphi(t)=\mu+R\cos(t)+\mathrm{i}R\sin(t)$ which satisfies $\varphi\in C^\infty$. The number of eigenvalues including their multiplicity within the contour $\gamma$ is denoted by $n(\gamma)$. With the randomly chosen matrix $\hat{\mathbf{V}}\in \mathbb{C}^{m\times \ell}$ with $m\gg\ell\geq n(\gamma)$ the two contour integrals of the form
\begin{eqnarray*}
	\mathbf{A}_0&=&\frac{1}{2\pi\mathrm{i}}\int_\gamma \mathbf{M}^{-1}(k)\hat{\mathbf{V}}\,\mathrm{d}s(k)\,,\\
	\mathbf{A}_1&=&\frac{1}{2\pi\mathrm{i}}\int_\gamma k\mathbf{M}^{-1}(k)\hat{\mathbf{V}}\,\mathrm{d}s(k)
\end{eqnarray*}
over the given contour $\gamma$ are now rewritten as
\begin{eqnarray*}
	\mathbf{A}_0&=&\frac{1}{2\pi\mathrm{i}}\int_{0}^{2\pi} \mathbf{M}^{-1}(\varphi(t))\hat{\mathbf{V}}\varphi'(t)\,\mathrm{d}s(t)\,,\\
		\mathbf{A}_1&=&\frac{1}{2\pi\mathrm{i}}\int_{0}^{2\pi} \varphi(t)\mathbf{M}^{-1}(\varphi(t))\hat{\mathbf{V}}\varphi'(t)\,\mathrm{d}s(t)
\end{eqnarray*}
and approximated by the trapezoidal rule yielding
\begin{eqnarray*}
	\mathbf{A}_{0,N}&=&\frac{1}{\mathrm{i}N}\sum_{j=0}^{N-1} \mathbf{M}^{-1}(\varphi(t_j))\hat{\mathbf{V}}\varphi'(t_j),\\
	\mathbf{A}_{1,N}&=&\frac{1}{\mathrm{i}N}\sum_{j=0}^{N-1} \varphi(t_j)\mathbf{M}^{-1}(\varphi(t_j))\hat{\mathbf{V}}\varphi'(t_j)\,,
\end{eqnarray*}
where the parameter $N$ is given and the equidistant nodes are $t_j=2\pi j/N$, $j=0,\ldots,N$. Note that the choice $N=24$ is usually sufficient, which is due to the exponential convergence rate (\cite[Theorem 4.7]{beyn}). The next step is the computation of a singular value decomposition of $\mathbf{A}_{0,N}=\mathbf{V}\mathbf{\Sigma}\mathbf{W}^\mathrm{H}$ with $\mathbf{V}\in \mathbb{C}^{m\times \ell}$, $\mathbf{\Sigma}\in \mathbb{C}^{\ell\times \ell}$, and $\mathbf{W}\in \mathbb{C}^{\ell\times \ell}$. Then, we perform a rank test for the matrix $\mathbf{\Sigma}=\mathrm{diag}(\sigma_1,\sigma_2,\ldots,\sigma_\ell)$ for a given tolerance $\epsilon=\text{tol}_\text{rank}$ (usually $\epsilon=10^{-4}$). That is, find $n(\gamma)$ such that $\sigma_1\geq\ldots\geq \sigma_{n(\gamma)}>\epsilon>\sigma_{n(\gamma)+1}\geq \ldots \geq \sigma_\ell$. Define $\mathbf{V}_0=(\mathbf{V}_{ij})_{1\leq i\leq m,1\leq j\leq n(\gamma)}$, $\mathbf{\Sigma}_0=(\mathbf{\Sigma}_{ij})_{1\leq i\leq n(\gamma),1\leq j\leq n(\gamma)}$, and $\mathbf{W}_0=(\mathbf{W}_{ij})_{1\leq i\leq \ell,1\leq j\leq n(\gamma))}$ and compute the $n(\gamma)$ eigenvalues $k_i$ and eigenvectors $\vec{s}_i$ of the matrix
$\mathbf{B}=\mathbf{V}_0^\mathrm{H}\mathbf{A}_{1,N}\mathbf{W}_0\mathbf{\Sigma}_0^{-1}\in\mathbb{C}^{n(\gamma)\times n(\gamma)}$. The $i$-th non-linear eigenvector $\vec{u}_i$ is given by $\mathbf{V}_0\vec{s}_i$. 

We refer the reader to \cite[p. 3849]{beyn} for more details on the implementation of this algorithm and the detailed analysis behind it including the proof of exponential convergence.

\subsection{Eigenfunction}
After we obtain the smallest non-zero eigenvalue $k$ and the corresponding function $u$ on the boundary from (\ref{nonlin2}), we insert this into (\ref{start}) to compute the eigenfunction inside the domain at any point we want. The discretization of the integrals is done as explained previously. Precisely, we have
\begin{eqnarray*}
	u(x)=-\mathrm{DL}_k^{\Gamma_1}u_1(x)-\mathrm{DL}_k^{\Gamma_2}u_2(x)
	\approx  \sum_{j=1}^{n_f}\sum_{k=1}^3\left(\hat{a}_{j,k} u_1(\widetilde{v}_{j,k})+\hat{b}_{j,k} u_2(\widetilde{v}_{j,k})\right)
\end{eqnarray*}
with 
\begin{eqnarray*}
	\hat{a}_{j,k}&=&\int_{\sigma}\partial_{\nu_1(\widetilde{m}_j(s))}\Phi_k(x,\widetilde{m}_j(s))\|\partial_s \widetilde{m}_j(s)\|\widetilde{L}_k(s)\,\mathrm{d}s(s)\\
	\hat{b}_{j,k}&=&\int_{\sigma}\partial_{\nu_2(\widetilde{m}_j(s))}\Phi_k(x,\widetilde{m}_j(s))\|\partial_s \widetilde{m}_j(s)\|\widetilde{L}_k(s)\,\mathrm{d}s(s)
	\end{eqnarray*}
for an arbitrary point $x\in D$.
In fact, we can find maximal or minimal values by maximizing or minimizing this function. This is done by the Nelder-Mead algorithm (in Matlab by the fminsearch function), refer also to \cite{lagarias}.
  
\subsection{Superconvergence}
The convergence of the method is out of the scope of this paper. Standard convergence results are available for boundary integral equations of the second kind using boundary element collocation method under 
suitable assumptions on the boundary (for example the boundary is at least of class $C^2$) 
and the boundary condition for the Laplace equation (see \cite{Atkinson1997}). Quadratic approximations of the boundary and the boundary function yield cubic convergence (refer to \cite{Atkinson1997} 
for the Laplace equation and \cite{kleefeldlin2} for the Helmholtz equation). In fact, the convergence results can be improved as shown in \cite{kleefeldlin2} for the three-dimensional case. 
However, the exact theoretical convergence rate for the eigenvalue is not known. It is expected that it is at least of order three, but we see later in the numerical results that it is better than three (for a sufficiently smooth 
boundary). 
Future work in this direction could be done using the ideas of Steinbach \& Unger \cite{unger}. 

Finally, note that a specific choice of $0<\alpha<1/2$ can improve the overall convergence rate for smooth boundaries (see \cite{kleefeldlin2}), but since we are happy with the cubic convergence rate with the pick 
$\alpha=(1-\sqrt{3/5})/2$ (a Gauss-quadrature point within the interval $[0,1]$), 
we have not investigated this any further.

\section{Numerical results}\label{numa}
\subsection{Simply-connected convex domains}
First, we check the correctness and the convergence of the underlying method for the unit circle. It is known that the first non-trivial interior Neumann eigenvalue is the smallest positive root of $J_1'$ 
(the first derivative of the first kind Bessel function of order one). The root can be computed to arbitrary precision with Maple with the command 
\begin{verbatim}
 restart; Digits:=16: fsolve(diff(BesselJ(1,x),x),x=1..2);
\end{verbatim}
It is approximately given by 
\begin{equation}
 1.841\,183\,781\,340\,659\;.
 \label{firstresult}
\end{equation}
For the Beyn algorithm we use the parameters $N=24$, $R=1/2$, $\mu=2$, and $\ell=10$ for various number of faces $n_f$ and number of collocation points $n_c$. 
With the definition of the absolute error $E_{n_f}^{(i)}$ of the $i$-th eigenvalue approximation, we compute the error of the first non-trivial eigenvalue $E_{n_f}^{(1)}$ of our method compared with (\ref{firstresult}). 
Additionally, we define the 
estimated error of convergence $\mathrm{EOC}^{(i)}=\log(E_{n_f}^{(i)}/E_{2\cdotp n_f}^{(i)})/\log(2)$ of the $i$-th eigenvalue approximation and compute $\mathrm{EOC}^{(1)}$.
As we can see in Table \ref{circle1}, the first non-trivial Neumann eigenvalue 
can be made accurate up to eleven digits.
\begin{table}[!ht]  
\caption{\label{circle1} Absolute error and estimated order of convergence of the first non-trivial interior Neumann eigenvalue for a unit circle using different number of faces and collocation points.}
\begin{indented}
 \item[]\begin{tabular}{@{}rrrl}
 \br
  $n_f$  & $n_c$ & abs. error $E_{n_f}^{(1)}$ & $\mathrm{EOC}^{(1)}$ \\
  \mr
     5 &   15 & $5.8503_{-3} $ &\\
    10 &   30 & $4.7818_{-4} $ &3.6129\\
    20 &   60 & $4.5775_{-5} $ &3.3849\\
    40 &  120 & $5.0168_{-6} $ &3.1897\\
    80 &  240 & $5.9173_{-7} $ &3.0838\\
   160 &  480 & $7.2096_{-8} $ &3.0369\\
   320 &  960 & $8.9069_{-9} $ &3.0169\\
   640 & 1920 & $1.1072_{-9} $ &3.0080\\
  1280 & 3840 & $1.3803_{-10}$ &3.0039\\
  \br
 \end{tabular}
 \end{indented}
\end{table}
The estimated order of convergence is at least of order three. Note that we can go beyond eleven digits accuracy by further increasing $n_f$, but it is not necessary here.

Next, we show the influence of the algebraic multiplicity of the eigenvalue. 
The first non-trivial interior Neumann eigenvalue for the unit circle has algebraic multiplicity two. The same is true for the second non-trivial interior Neumann eigenvalue. 
It is obtained by computing the first positive root of $J'_2$ approximately given by $3.054\,236\,928\,227\,140$. 
The third non-trivial interior Neumann eigenvalue is simple and obtained by computing the second root of $J_0'$ given by $3.831\,705\,970\,207\,512$. Note that
the first root of $J_0'$ is zero which corresponds to the interior Neumann eigenvalue zero with a corresponding eigenfunction which is a constant. In Table \ref{circle23}, we show the absolute error and the estimated order of convergence for 
the second and third non-zero interior Neumann eigenvalue for a unit circle where we used different number of faces and different number of collocation points using the parameters $N=24$, $R=1/2$, $\mu=3$, and $\ell=10$.

\begin{table}[!ht]  
\caption{\label{circle23} Absolute error and estimated order of convergence of the second and third non-trivial interior Neumann eigenvalue for a unit circle using different number of faces and collocation points.}
\begin{indented}
 \item[]\begin{tabular}{@{}rrrlrl}
 \br
  $n_f$  & $n_c$ & abs. error $E_{n_f}^{(2)}$ & $\mathrm{EOC}^{(2)}$ & abs. error $E_{n_f}^{(3)}$ & $\mathrm{EOC}^{(3)}$\\
  \mr
     5 &   15 & $1.2817_{-2} $ &       & $1.6335_{-2} $ &       \\
    10 &   30 & $1.1543_{-3} $ & 3.4729& $1.3081_{-3} $ & 3.6424\\
    20 &   60 & $1.2187_{-4} $ & 3.2437& $1.3133_{-4} $ & 3.3162\\
    40 &  120 & $1.4173_{-5} $ & 3.1041& $1.5147_{-5} $ & 3.1161\\
    80 &  240 & $1.7199_{-6} $ & 3.0428& $1.8416_{-6} $ & 3.0401\\
   160 &  480 & $2.1229_{-7} $ & 3.0182& $2.2788_{-7} $ & 3.0146\\
   320 &  960 & $2.6386_{-8} $ & 3.0082& $2.8373_{-8} $ & 3.0057\\
   640 & 1920 & $3.2895_{-9} $ & 3.0038& $3.5406_{-9} $ & 3.0024\\
  1280 & 3840 & $4.1066_{-10}$ & 3.0018& $4.4225_{-10}$ & 3.0011\\
  \br
 \end{tabular}
 \end{indented}
\end{table}

Again, we notice the estimated order of convergence of at least order three and that we can achieve the order and the high accuracy regardless of the algebraic multiplicity of the eigenvalue.

\begin{figure}[!ht]
\subfigure[First eigenfunction of a unit circle]{
\includegraphics[width=0.31\textwidth]{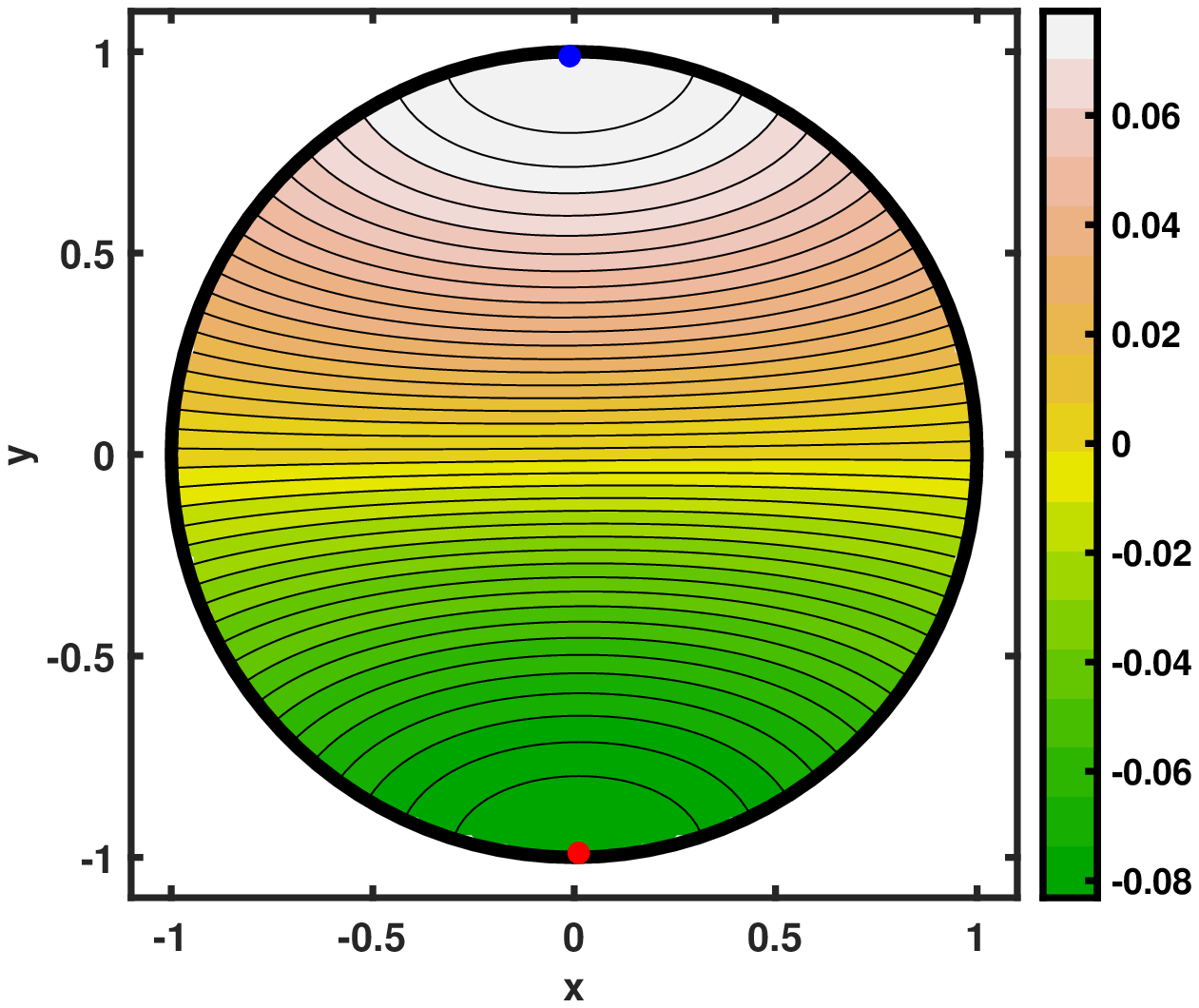}
}
\subfigure[Second eigenfunction of a unit circle]{
\includegraphics[width=0.31\textwidth]{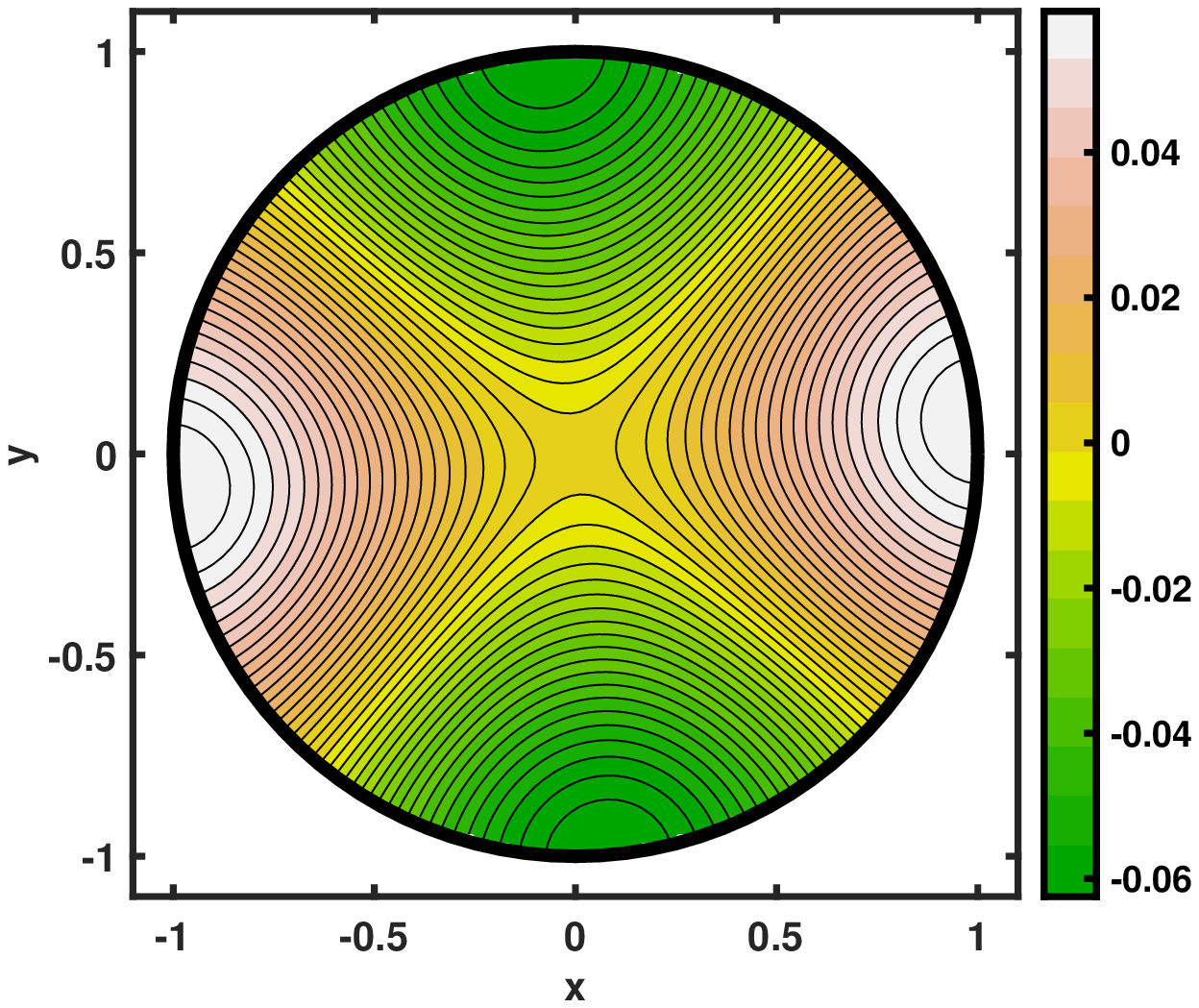}
}
\subfigure[Third eigenfunction of a unit circle]{
\includegraphics[width=0.31\textwidth]{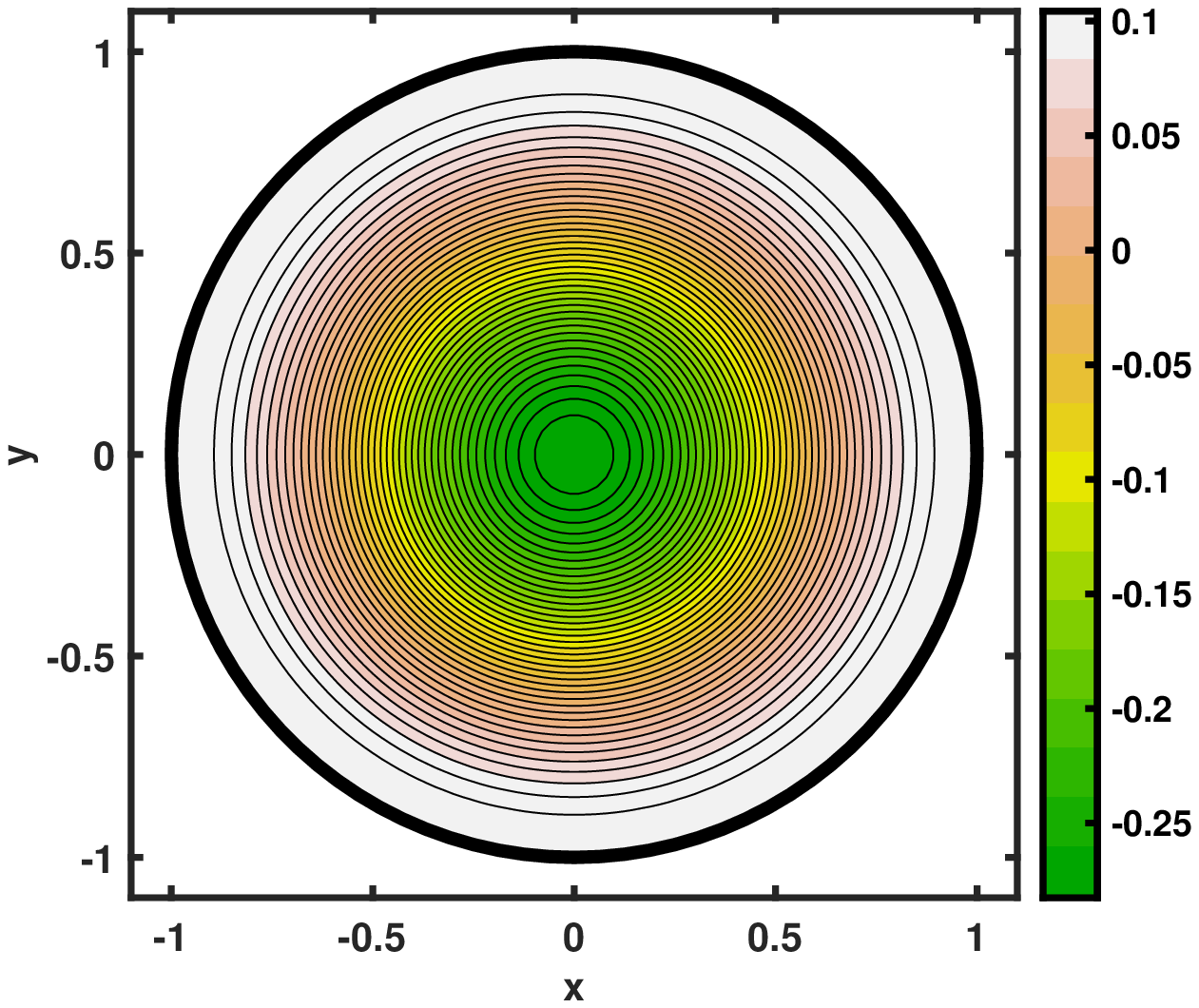}
}
\caption{\label{circle} The first three eigenfunctions corresponding to the first three non-trivial interior Neumann eigenvalues $1.841\,184$, $3.054\,237$, and $3.831\,706$ for the unit circle.}
\end{figure}

We define $\square_\kappa=[-\kappa,\kappa]\times [-\kappa,\kappa]$ with $\kappa\in\mathbb{R}_{>0}$.
In general, we use a resolution of $100\times 100$ equidistantly distributed points, here within $\square_{1.1}$ and compute for each point that is located inside the unit circle the value of the eigenfunction.
In Figure \ref{circle}, we show the eigenfunctions corresponding to the first three non-trivial interior Neumann eigenvalues as a contour plot with 40 contour lines. 
We also include the location of the maximum and minimum of the eigenfunction that corresponds to the
first non-trivial interior Neumann eigenvalue as a red and blue dot, respectively. 
We can see that the extreme values for the first non-trivial interior Neumann eigenfunction of the unit circle are obtained on the boundary as it is conjectured for simply-connected convex domains. Note that the second eigenfunction 
corresponding to the first non-trivial interior Neumann eigenvalue is a rotated version of the first eigenfunction.

\begin{figure}[!ht]
\subfigure[First eigenfunction of an ellipse]{
\includegraphics[width=0.31\textwidth]{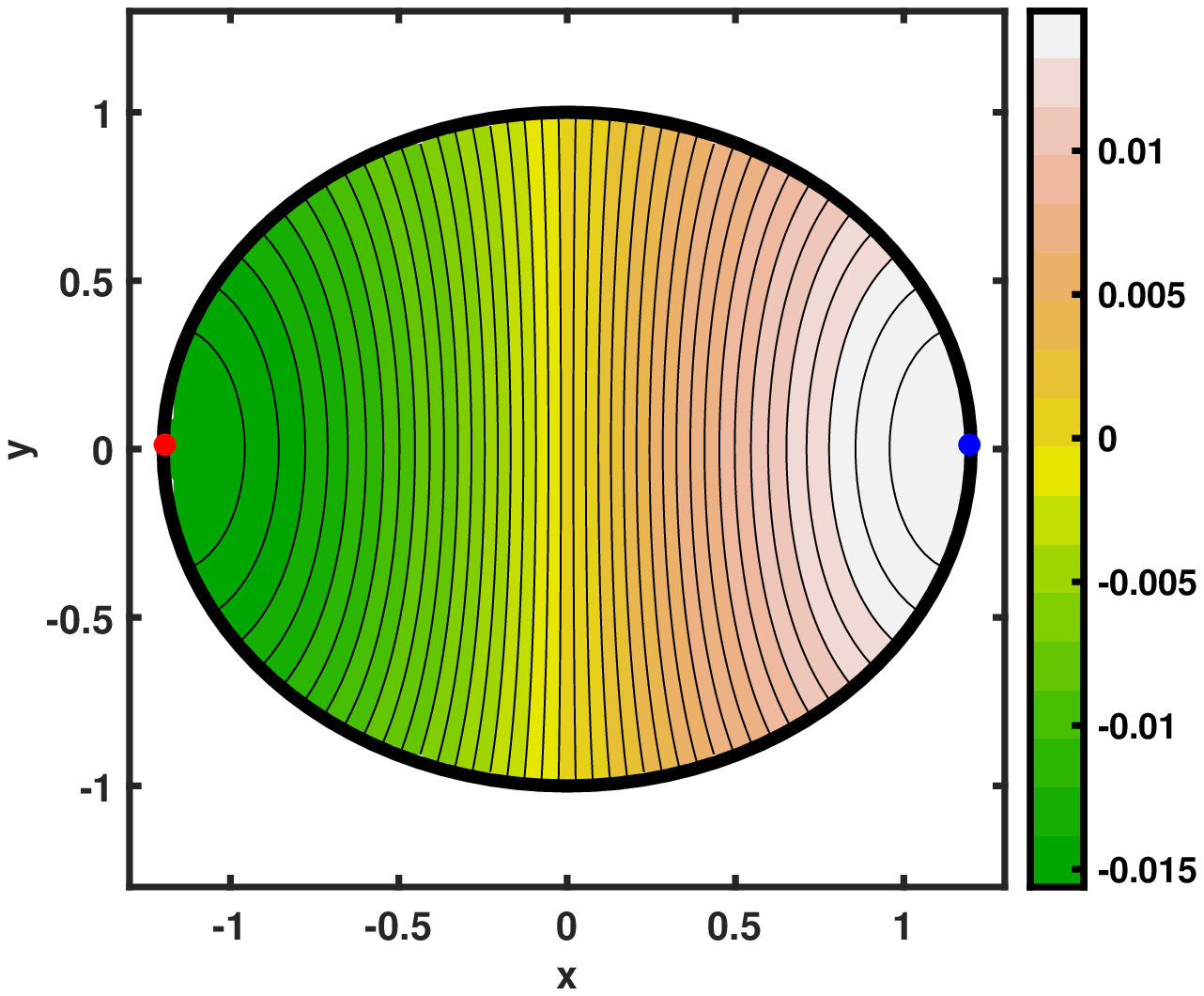}
}
\subfigure[First eigenfunction of a deformed ellipse with $\varepsilon=0.1$]{
\includegraphics[width=0.31\textwidth]{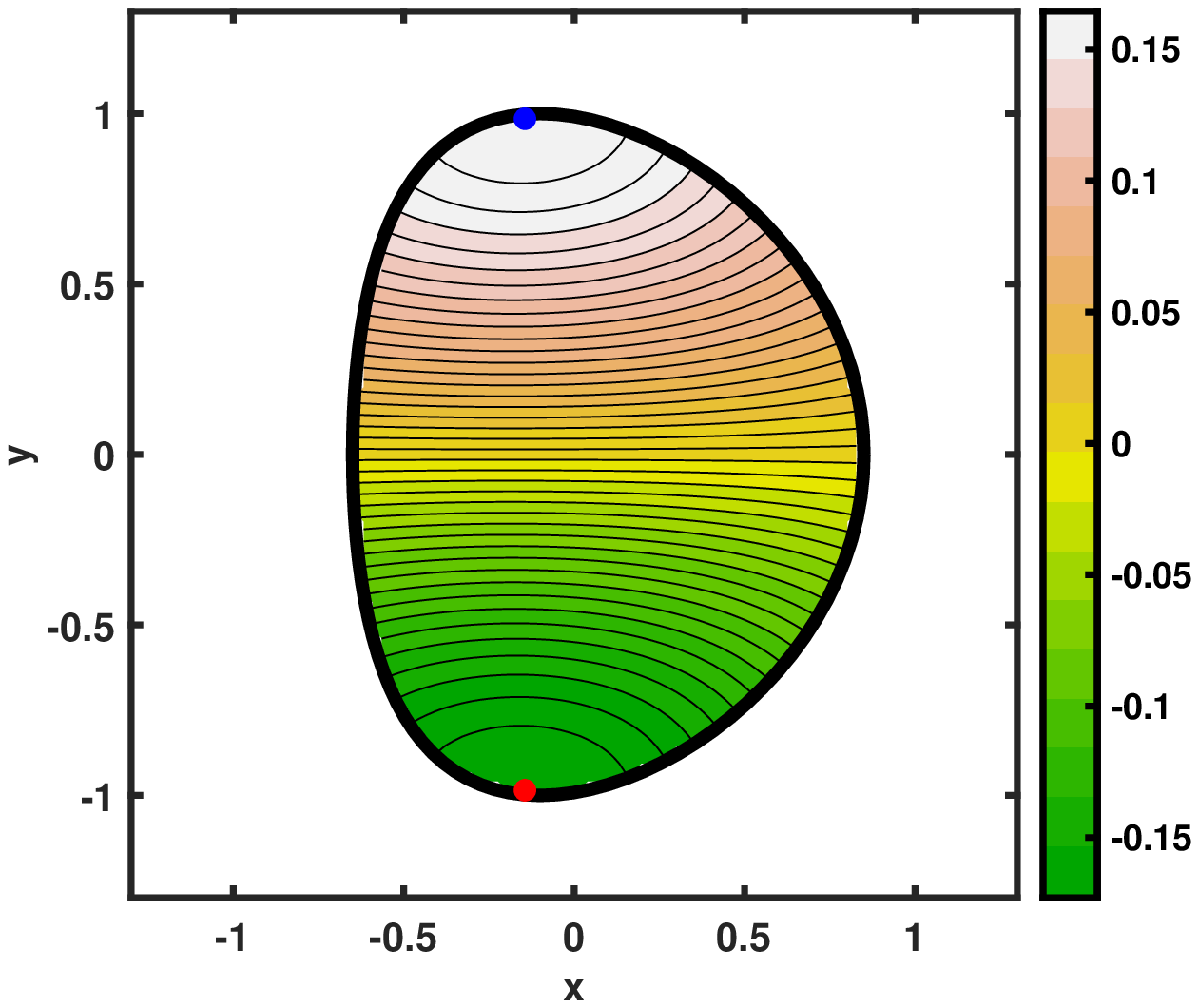}
}
\subfigure[First eigenfunction of a deformed ellipse with $\varepsilon=0.2$]{
\includegraphics[width=0.31\textwidth]{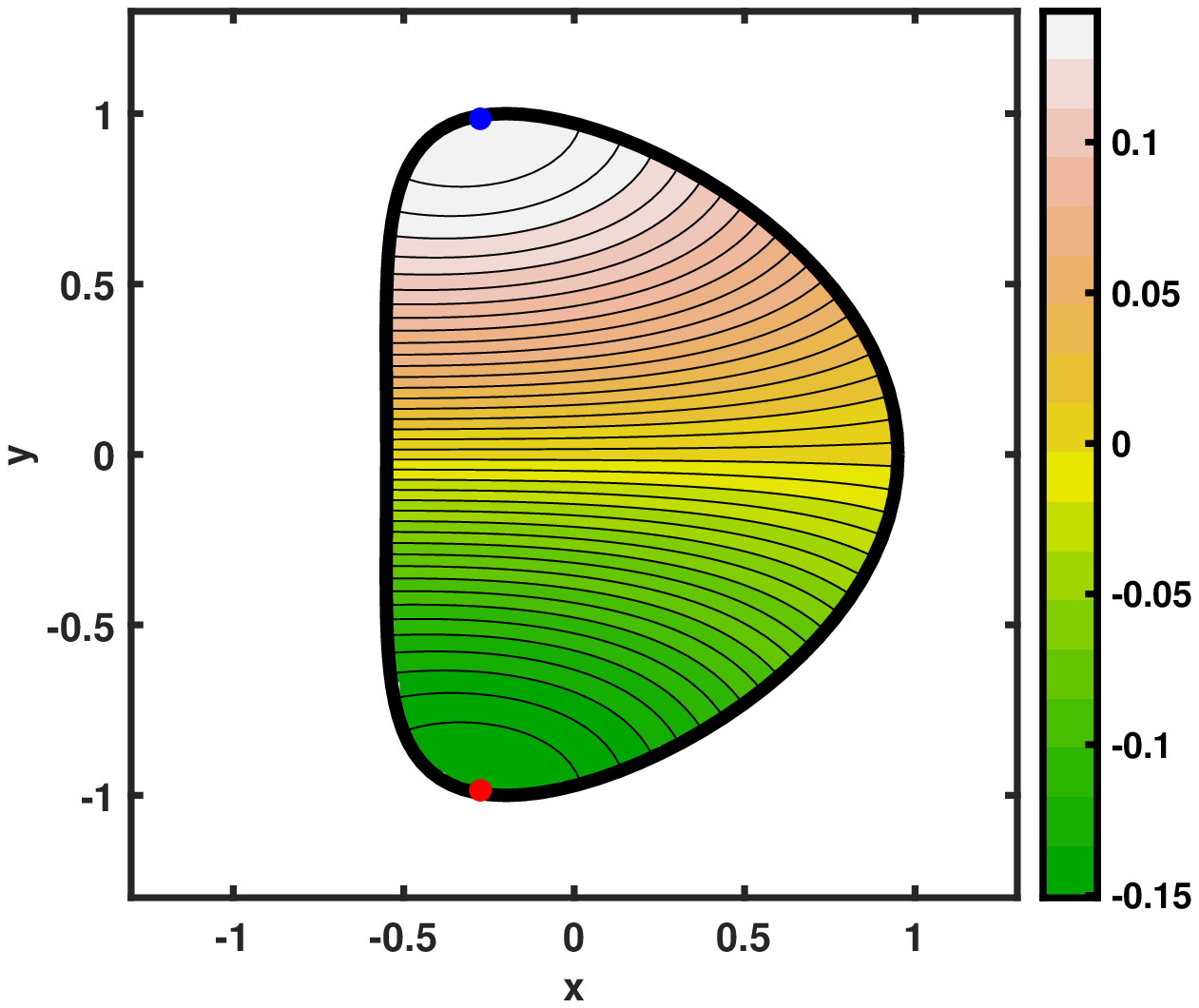}
}
\caption{\label{otherconvex} The eigenfunctions corresponding to the first non-trivial interior Neumann eigenvalues for the ellipse and two deformed ellipses with $\varepsilon=0.1$ and $\varepsilon=0.2$ with corresponding non-trivial 
interior Neumann eigenvalue $1.544\,422$, $1.849\,064$, and $1.819\,478$.}
\end{figure}

We also show the eigenfunctions including the maximal and minimal value for a variety of other simply-connected convex domains in Figure \ref{otherconvex} such as an ellipse and two deformed ellipses. 
The boundary of the 
ellipse is given in parametric form as $(6\cos(t)/5,\sin(t))^\top$ with $t\in [0,2\pi)$. We use the parameters as before for Beyn's algorithm except $\mu=1.5$ and consider $\square_{1.3}$. The first 
non-trivial interior Neumann eigenvalue is given by $1.544\,422$ which has algebraic multiplicity one. The parametrization of the deformed ellipse's boundary is given by $(0.75\cos(t)+\varepsilon\cos(2t),\sin(t))^\top$ with $t\in [0,2\pi)$, 
where the parameter $\varepsilon$ is chosen to be $0.1$ and $0.2$ (see \cite{cakonikress,kleefeldpieronek} for its first and second use). Using $n_f=320$, the first non-trivial interior Neumann eigenvalue of 
the deformed ellipses with $\varepsilon=0.1$ and $\varepsilon=0.2$ are $1.849\,064$ and $1.819\,478$, respectively. Again, they both have algebraic multiplicity one.
It is generally believed that the hot spots conjecture for general simply-connected convex domains is true,
but a general proof is still open.
In all our numerical results for simply-connected convex domains, we obtain the extrema on the boundary as one can
see in Figure \ref{otherconvex}.

The same is true when we consider piecewise smooth convex domains such as the unit square and the equilateral triangle with side length one. The first non-trivial interior 
Neumann eigenvalues are known to be $\pi$ and $4\pi/3$ (see \cite{grebenkov}). Their multiplicity is two. 

In Table \ref{squaretriangle} we see that our algorithm works fine with these two piecewise smooth domains. The estimated order of 
convergence is better than four and hence better than for the previously discussed smooth domains. Since we use quadratic interpolation of the smooth boundary, there is an approximation error limiting the convergence rate. 
For the considered piecewise smooth domains, the boundary is approximated exactly since it is a linear function thus explaining the better convergence. 
\begin{table}[!ht]  
\caption{\label{squaretriangle} Absolute error and estimated order of convergence of the first non-trivial interior Neumann eigenvalue for a unit square and an equilateral triangle with side length 
one using different number of faces and collocation points.}
\begin{indented}
 \item[]\begin{tabular}{@{}rrrlrrrl}
 \br
  $n_f$  & $n_c$ & abs. error $E_{n_f}^{\square}$ & $\mathrm{EOC}^{\square}$ & $n_f$  & $n_c$ & abs. error $E_{n_f}^{\Delta}$ & $\mathrm{EOC}^{\Delta}$\\
  \mr
     4 &   12 & $7.3586_{-3}  $ &       &    3 &    9 & $2.7096_{-2}  $ &       \\
     8 &   24 & $4.4226_{-4}  $ & 4.0565&    6 &   18 & $2.0723_{-3}  $ & 3.7088\\
    16 &   48 & $1.9060_{-5}  $ & 4.5363&   12 &   36 & $9.7843_{-5}  $ & 4.4046 \\
    32 &   96 & $7.6974_{-7}  $ & 4.6301&   24 &   72 & $4.1556_{-6}  $ & 4.5573\\
    64 &  192 & $3.0527_{-8}  $ & 4.6562&   48 &  144 & $1.7064_{-7}  $ & 4.6060\\
   128 &  384 & $1.2042_{-9}  $ & 4.6639&   96 &  288 & $6.9001_{-8}  $ & 4.6282\\
   256 &  768 & $4.7429_{-11} $ & 4.6662&  192 &  576 & $2.7581_{-10} $ & 4.6449\\
   512 & 1536 & $1.8665_{-12} $ & 4.6674&  384 & 1152 & $1.0880_{-11} $ & 4.6640\\
  \br
 \end{tabular}
 \end{indented}
\end{table}

\begin{figure}[!ht]
\subfigure[First eigenfunction of a unit square]{
\includegraphics[width=0.31\textwidth]{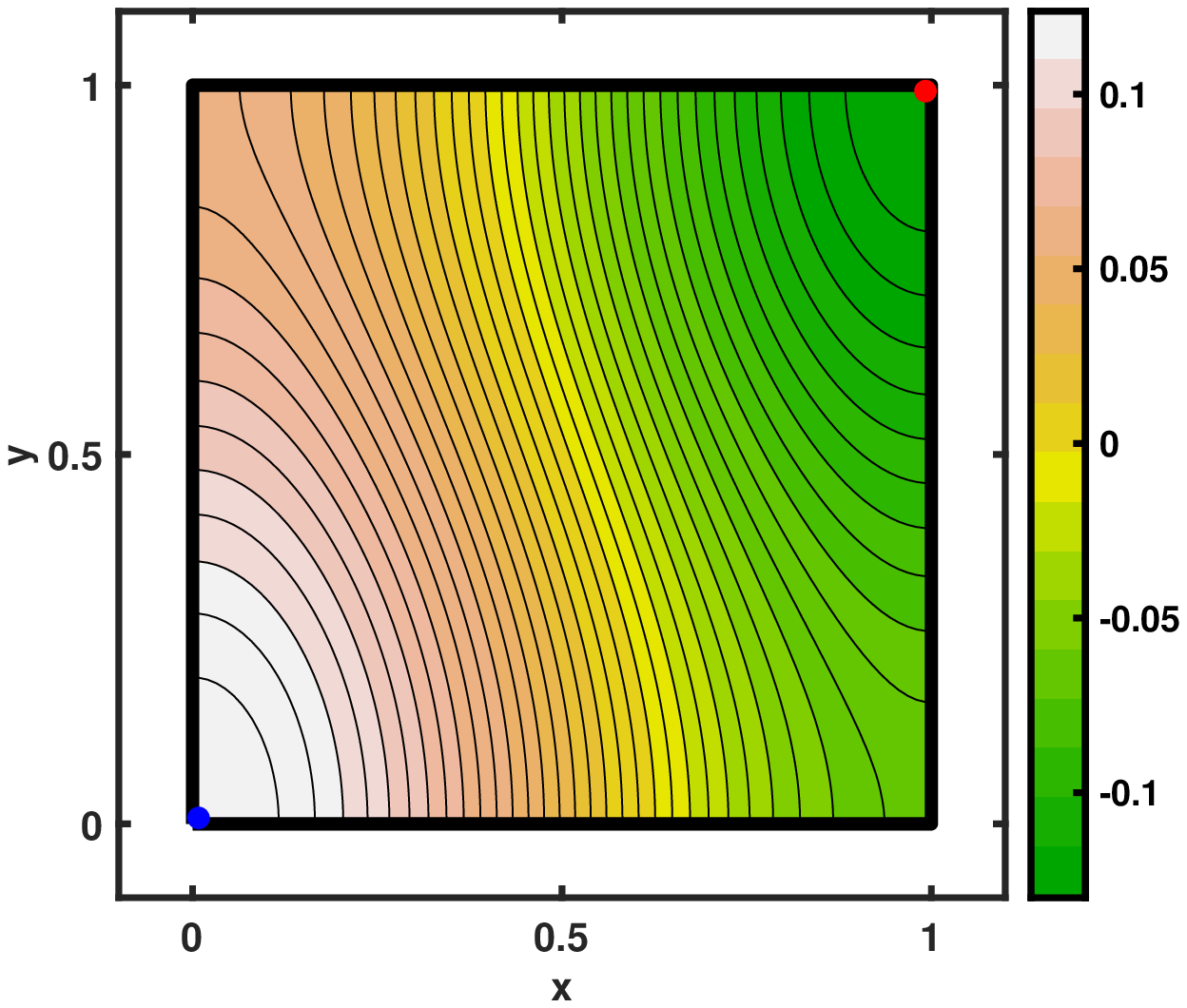}
}
\subfigure[First eigenfunction of an equilateral triangle with side length one]{
\includegraphics[width=0.31\textwidth]{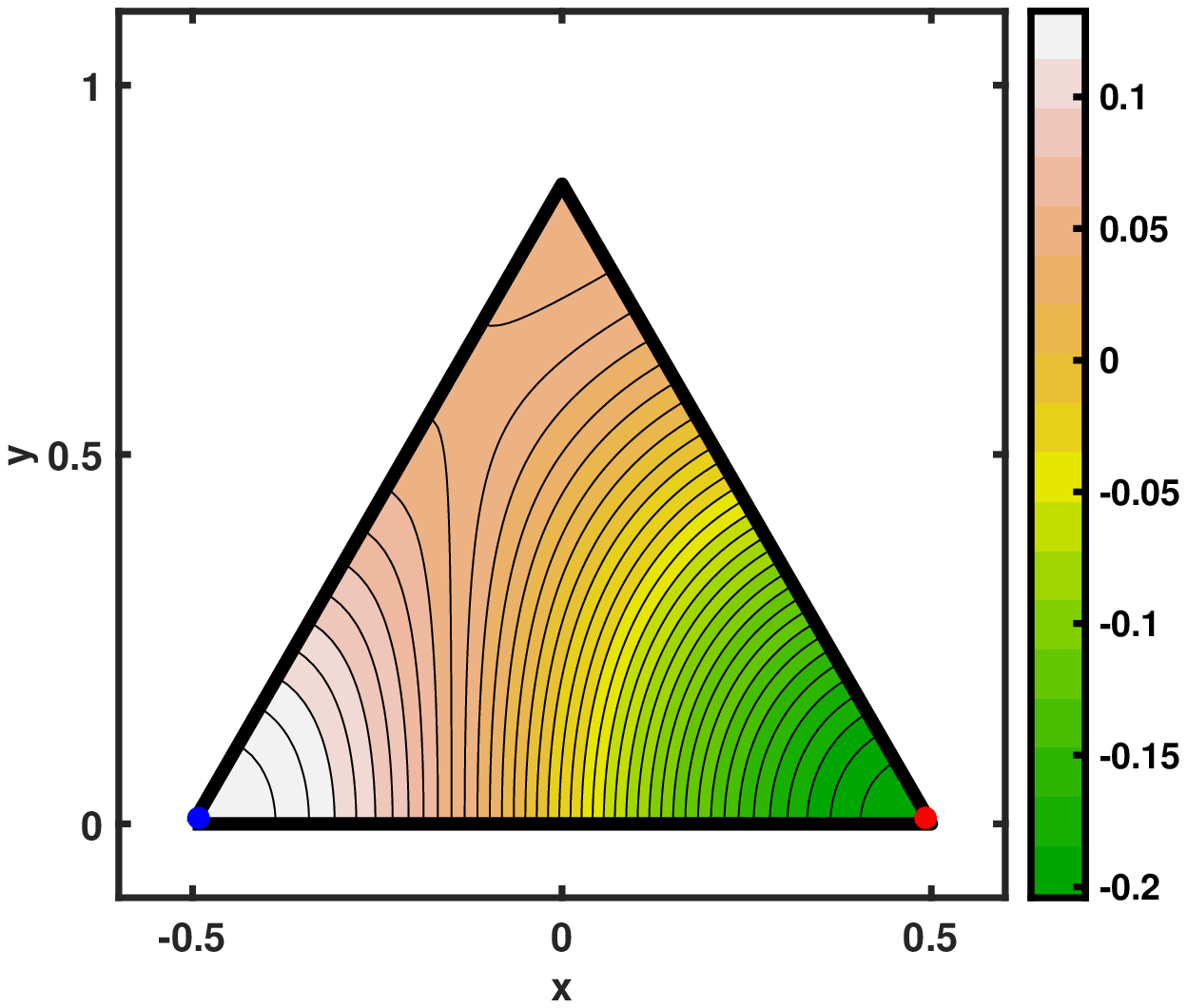}
}
\caption{\label{sqtr} One of the eigenfunctions corresponding to the first non-trivial interior Neumann eigenvalues for the unit square and the equilateral triangle with side length one 
with corresponding non-trivial 
interior Neumann eigenvalue $\pi$ and $4\pi/3$.}
\end{figure}

Later, we see that these nice convergence rates depend on the regularity of the solution at a corner and we obtain worse approximation results. 
In Figure \ref{sqtr} we show one of the corresponding eigenfunctions for the unit square (refer also to Figure \ref{introex}) and the equilateral triangle with side length one including the location of the maximum and minimum. 
As we can see, they are located on the boundary.

\subsection{Simply-connected non-convex domains}
\begin{figure}[!ht]
\subfigure[First eigenfunction of a deformed ellipse with $\varepsilon=0.3$]{
\includegraphics[width=0.31\textwidth]{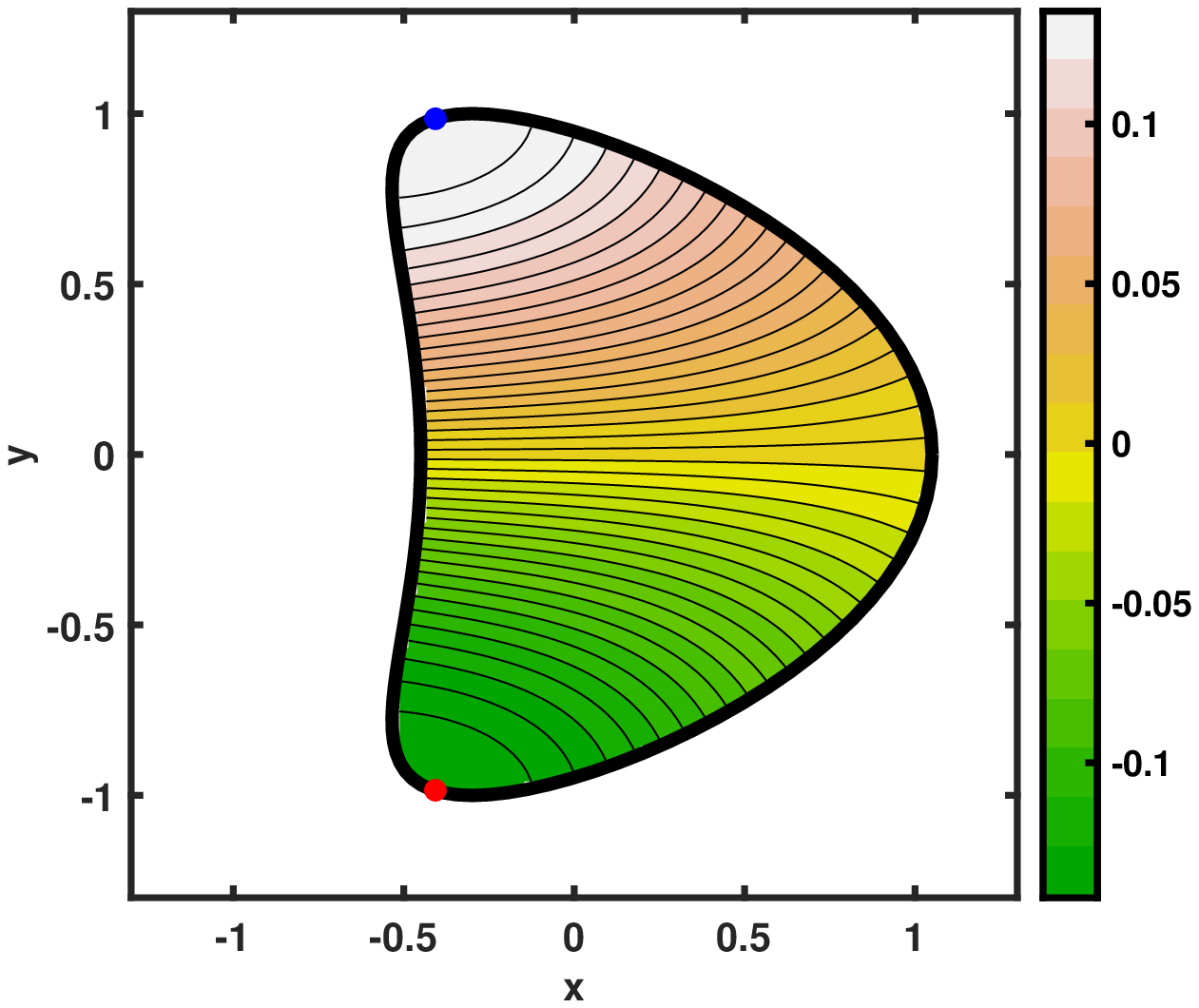}
}
\subfigure[First eigenfunction of a peanut]{
\includegraphics[width=0.31\textwidth]{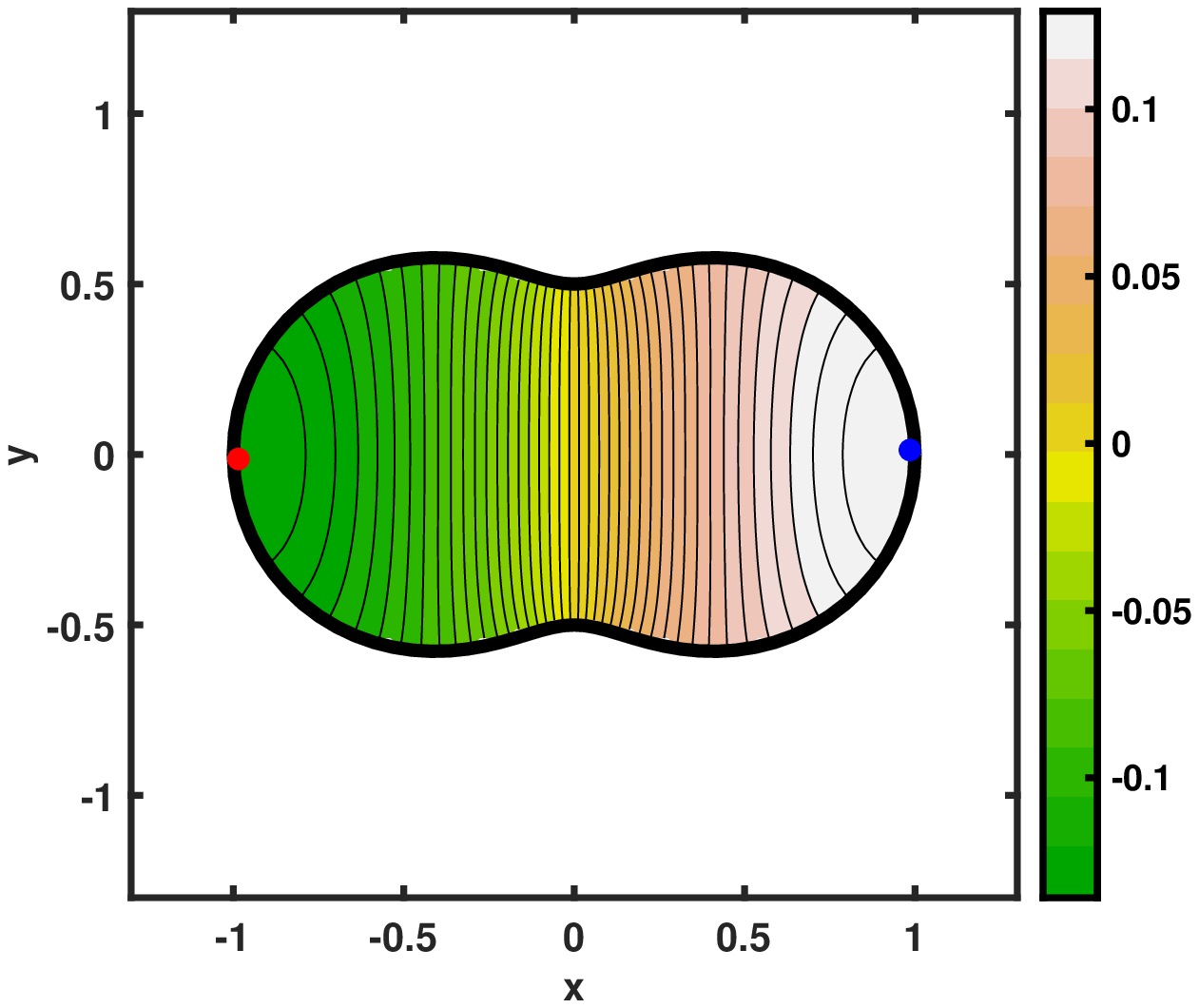}
}
\subfigure[First eigenfunction of an apple]{
\includegraphics[width=0.31\textwidth]{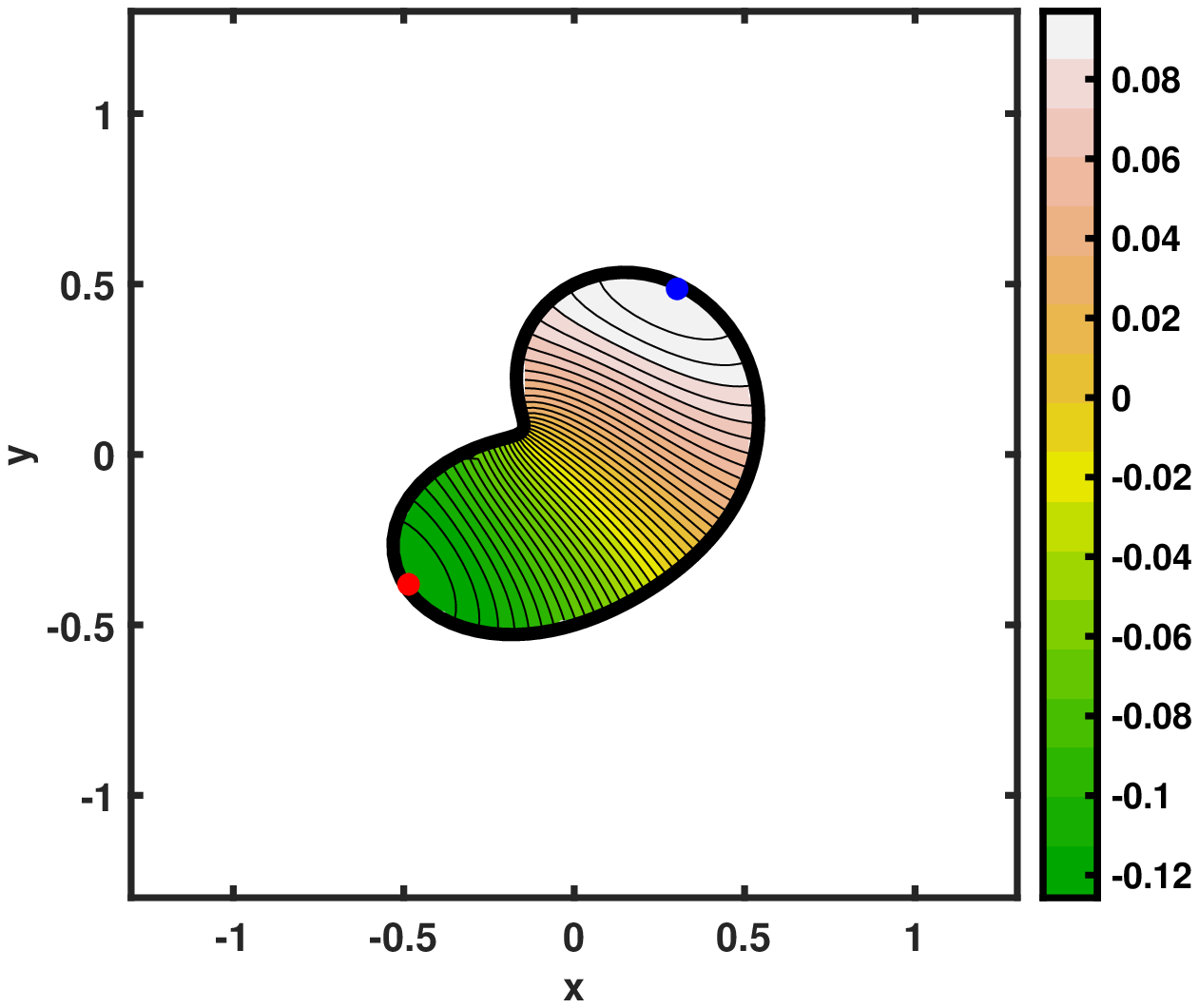}
}
\caption{\label{nonconvex} The eigenfunctions corresponding to the first non-trivial interior Neumann eigenvalues for the deformed ellipse with $\varepsilon=0.3$, the peanut-shaped, and the apple-shaped domain 
with corresponding non-trivial 
interior Neumann eigenvalue $1.770\,906$, $1.721\,292$, and $2.761\,274$.}
\end{figure}

No simply-connected non-convex domain has yet been found (neither theoretically nor numerically) 
that fails the hot spots conjecture. 
Now, we concentrate on this 
case. We consider the deformed ellipse from the previous section with $\varepsilon=0.3$, the peanut-shaped domain 
and the apple-shaped domain. The boundaries of the last two domains are given parametrically as 
$\sqrt{\cos^2(t)+\sin^2(t)/4}\left(\cos(t),\sin(t)\right)^\top$ and 
$\frac{0.5+0.4\cos(t)+0.1\sin(2t)}{1+0.7\cos(t)}\left(\cos(t),\sin(t)\right)^\top$ with $t\in [0,2\pi)$, respectively 
(see \cite{haiwen} for its use). We use the parameters as before with $\mu=1.5$ for the first two domains and $\mu=3$ for the apple-shaped domain.
Using $n_f=320$, the first non-trivial interior Neumann eigenvalue for the three domains are $1.770\,906$, $1.721\,292$, and $2.761\,274$, respectively. They are all simple.
As we can see in Figure \ref{nonconvex}, the maximum and minimum are obtained on the boundary of the domains using $\square_{1.3}$.

Interesting domains have been constructed by Kleefeld (refer to \cite{kleefeld2019shape} for more details) and extended by Abele and Kleefeld \cite{abele} for the purpose of finding new shape optimizers for certain 
non-trivial Neumann eigenvalues. The boundaries of the domains considered in those articles are given by `generalized' equipotentials which are implicit curves. The simplest equipotential is of the form 
$\sum_{i=1}^m \|x-P_i\|^{-1}=c$ where the points $P_i\in \mathbb{R}^2$, 
$i=1,\ldots,m$, the number of points $m$, and the 
parameter $c$ are given. All $x\in \mathbb{R}^2$ that satisfy the equation describe the boundary of the domain. We use this idea to construct three non-symmetric and non-convex 
simply-connected domains, say $D_1$, $D_2$, and $D_3$.
For the boundary of the domain $D_1$, we use the parameter $m=3$ and $c=14/5$. The three points are $(-1,1/2)^\top$, $(1,1/3)^\top$, and 
$(0,4/5)^\top$. Using $\mu=1$, we obtain the first non-trivial interior Neumann eigenvalue $1.051\,055$. For the plot of the eigenfunction, we use $\square_{1.6}$.
The boundary of the second domain is constructed through the use of $m=4$ and $c=7/2$ and the points are $(-5/4,1/10)^\top$, $(5/4,0)^\top$, $(1/10,-1)^\top$, and $(0,1)^\top$. 
We obtain the first non-trivial interior Neumann eigenvalue $1.086\,037$ when using $\mu=1$ and the eigenfunction using $\square_{1.8}$. The third domain's boundary is constructed using $m=5$ and $c=18/5$ and the points are
$(-1,1/2)^\top$, $(1,1/2)^\top$, $(0,-1)^\top$, $(3/2,-1)^\top$, 
and $(-3/2,-6/5)^\top$. Using $\mu=1$, we obtain the first non-trivial interior Neumann eigenvalue $0.861\,858$. A plot of the corresponding eigenfunction is shown within $\square_{2.1}$.

\begin{figure}[!ht]
\subfigure[First eigenfunction of $D_1$]{
\includegraphics[width=0.31\textwidth]{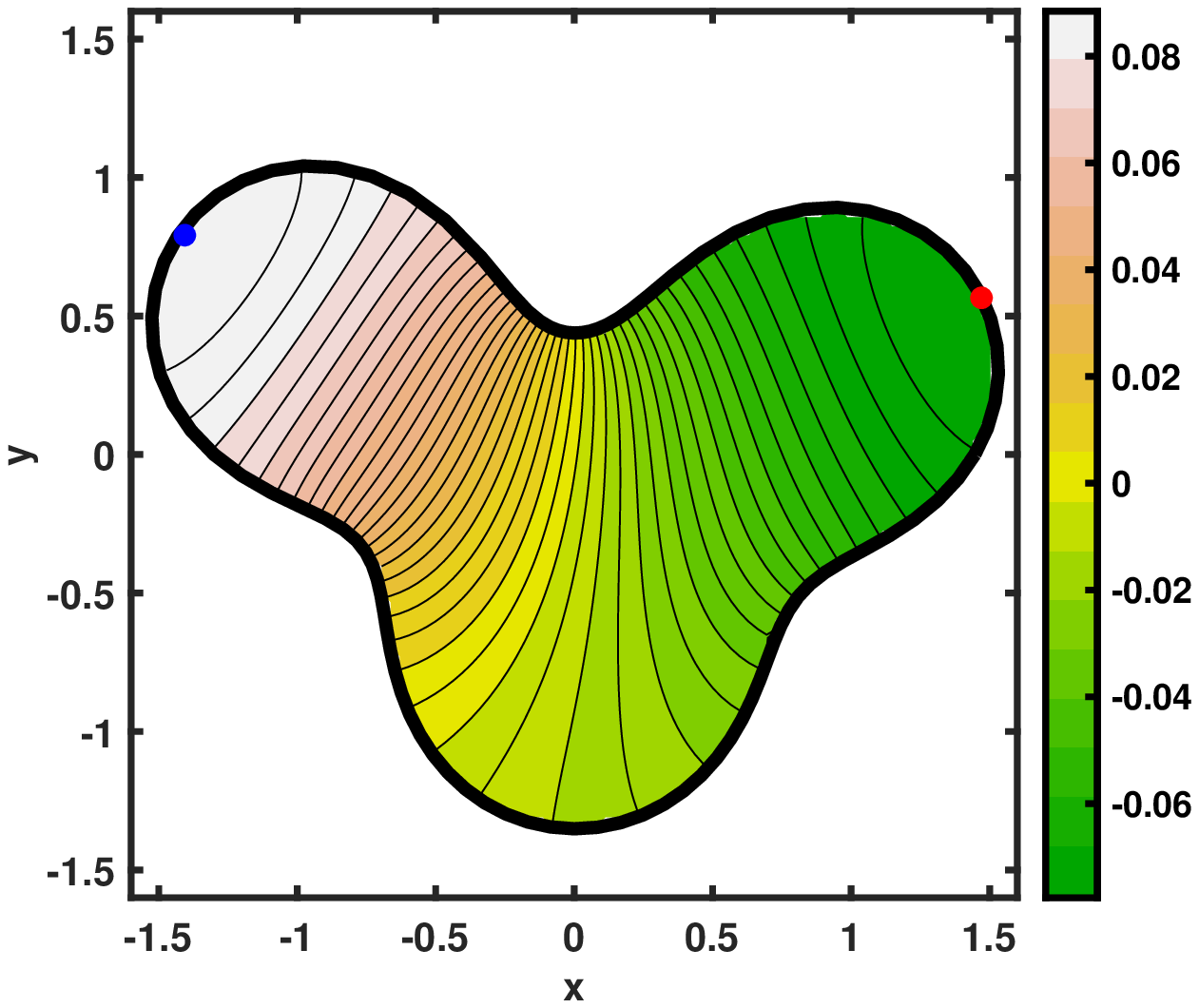}
}
\subfigure[First eigenfunction of $D_2$]{
\includegraphics[width=0.31\textwidth]{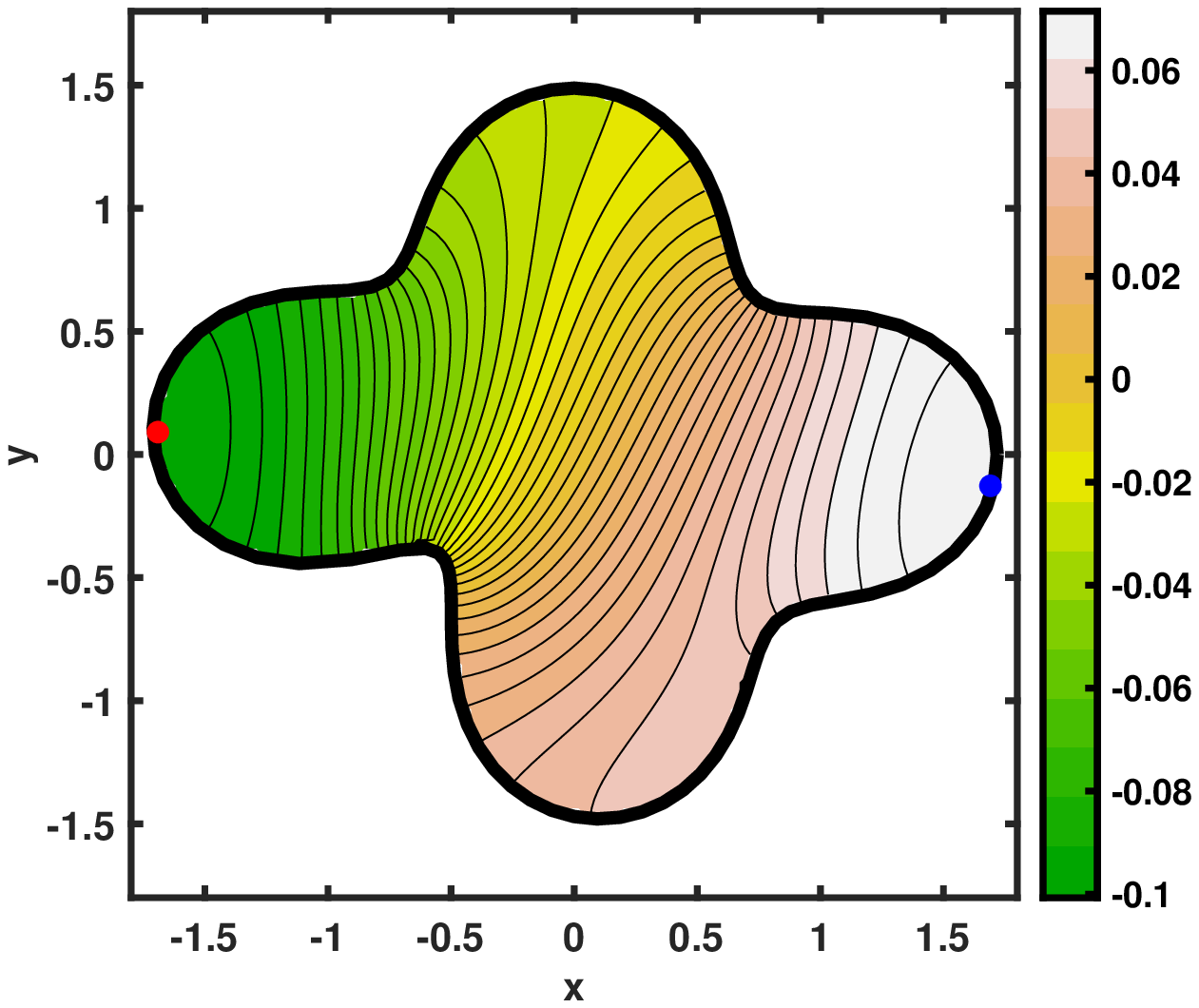}
}
\subfigure[First eigenfunction of $D_3$]{
\includegraphics[width=0.31\textwidth]{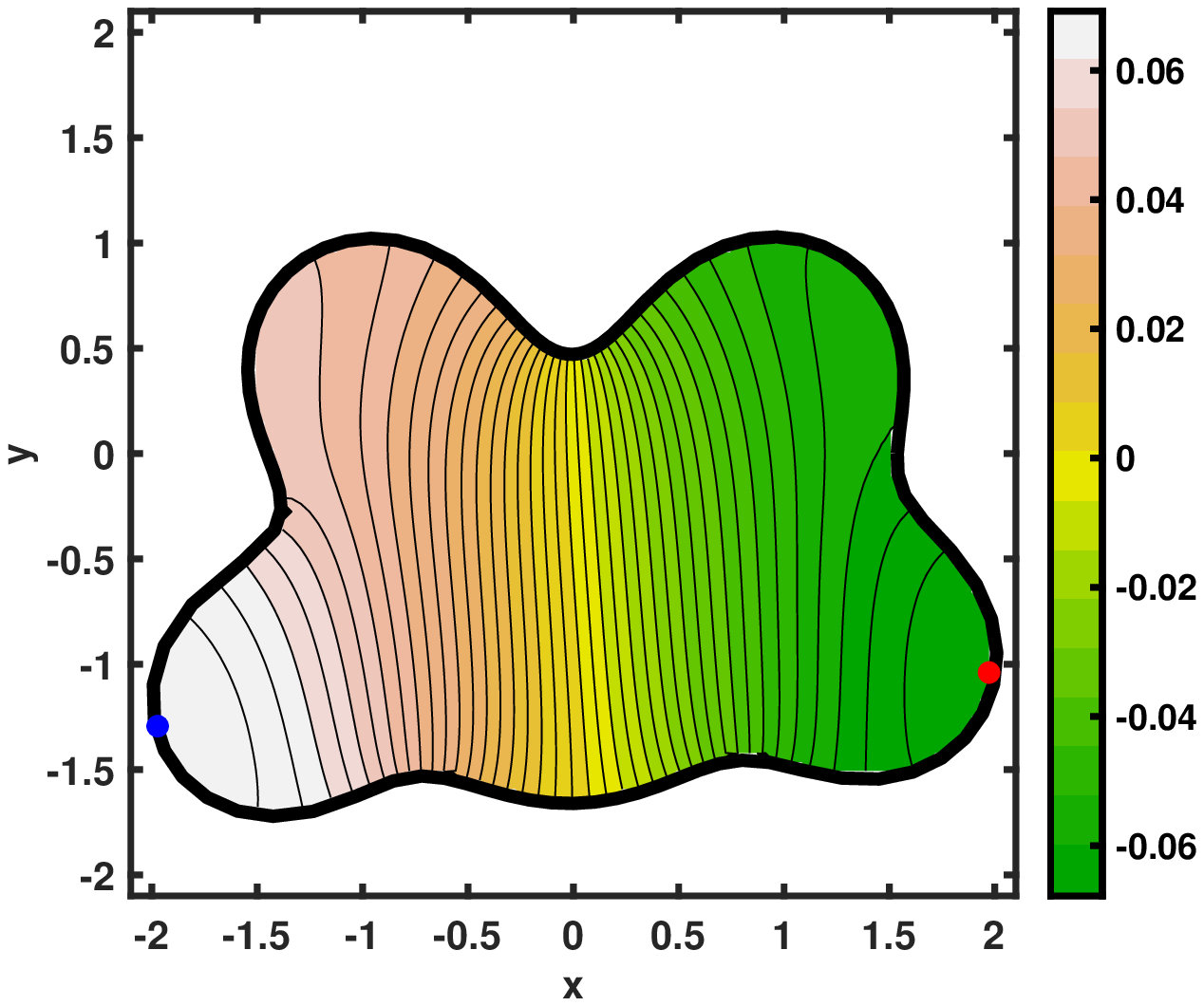}
}
\caption{\label{morecomplex} The eigenfunctions corresponding to the first non-trivial interior Neumann eigenvalues for domains $D_1$, $D_2$, and $D_3$
with corresponding non-trivial 
interior Neumann eigenvalue $1.051\,055$, $1.086\,037$, and $0.861\,858$.}
\end{figure}

\begin{figure}[!ht]
\subfigure[First eigenfunction of $\llcorner$]{
\includegraphics[width=0.31\textwidth]{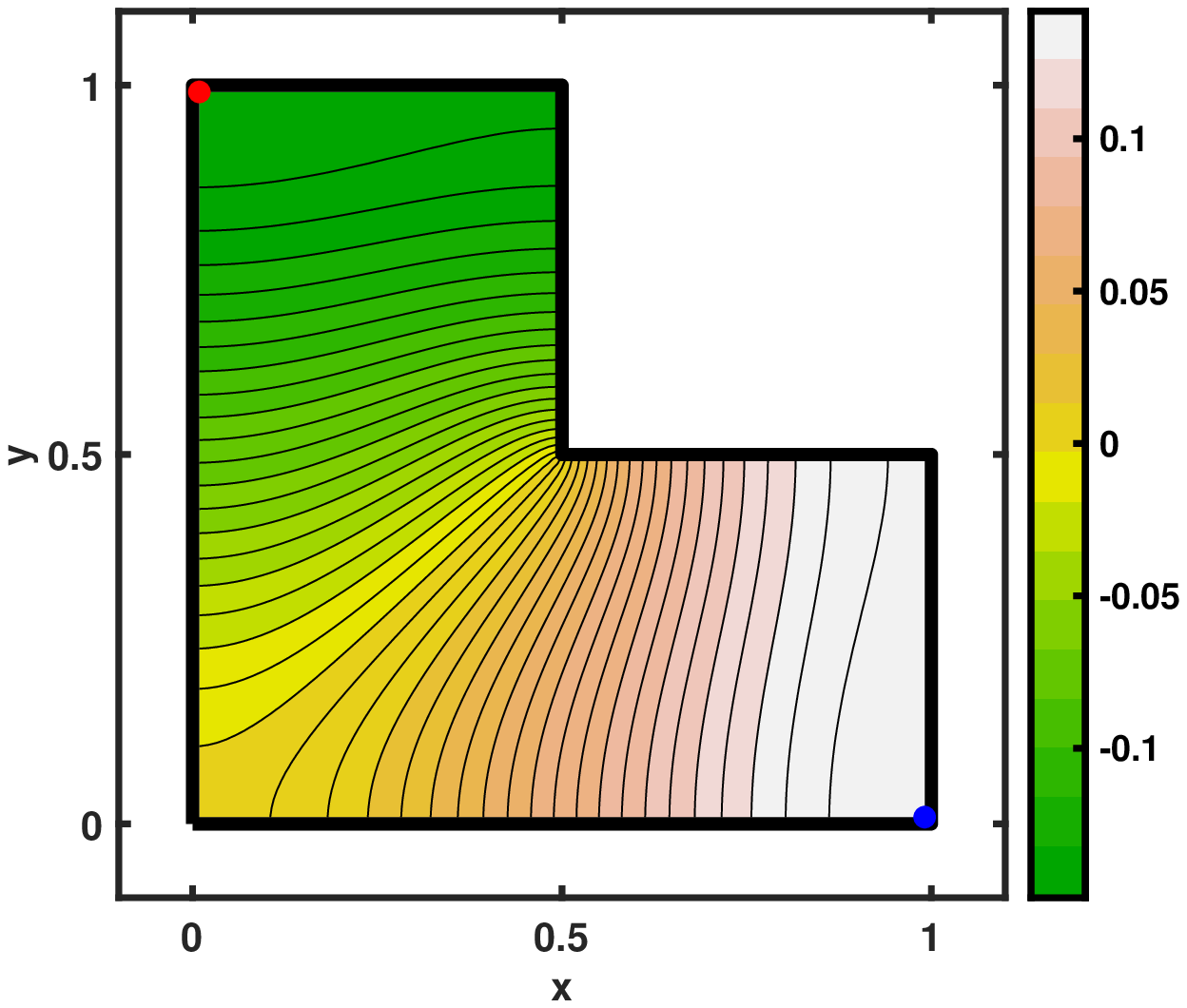}
}
\subfigure[First eigenfunction of $\llcorner_2$]{
\includegraphics[width=0.31\textwidth]{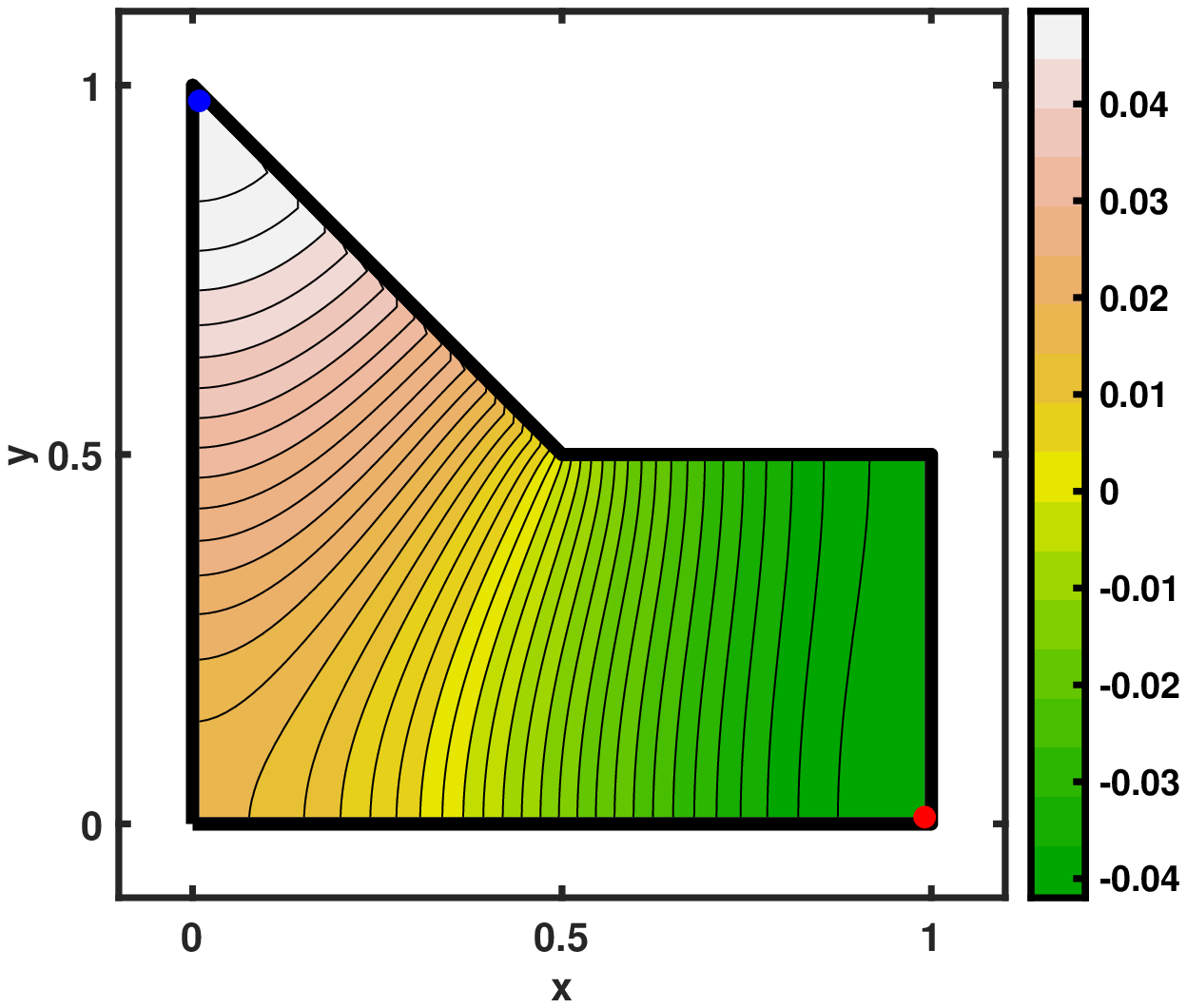}
}
\subfigure[First eigenfunction of $\llcorner_3$]{
\includegraphics[width=0.31\textwidth]{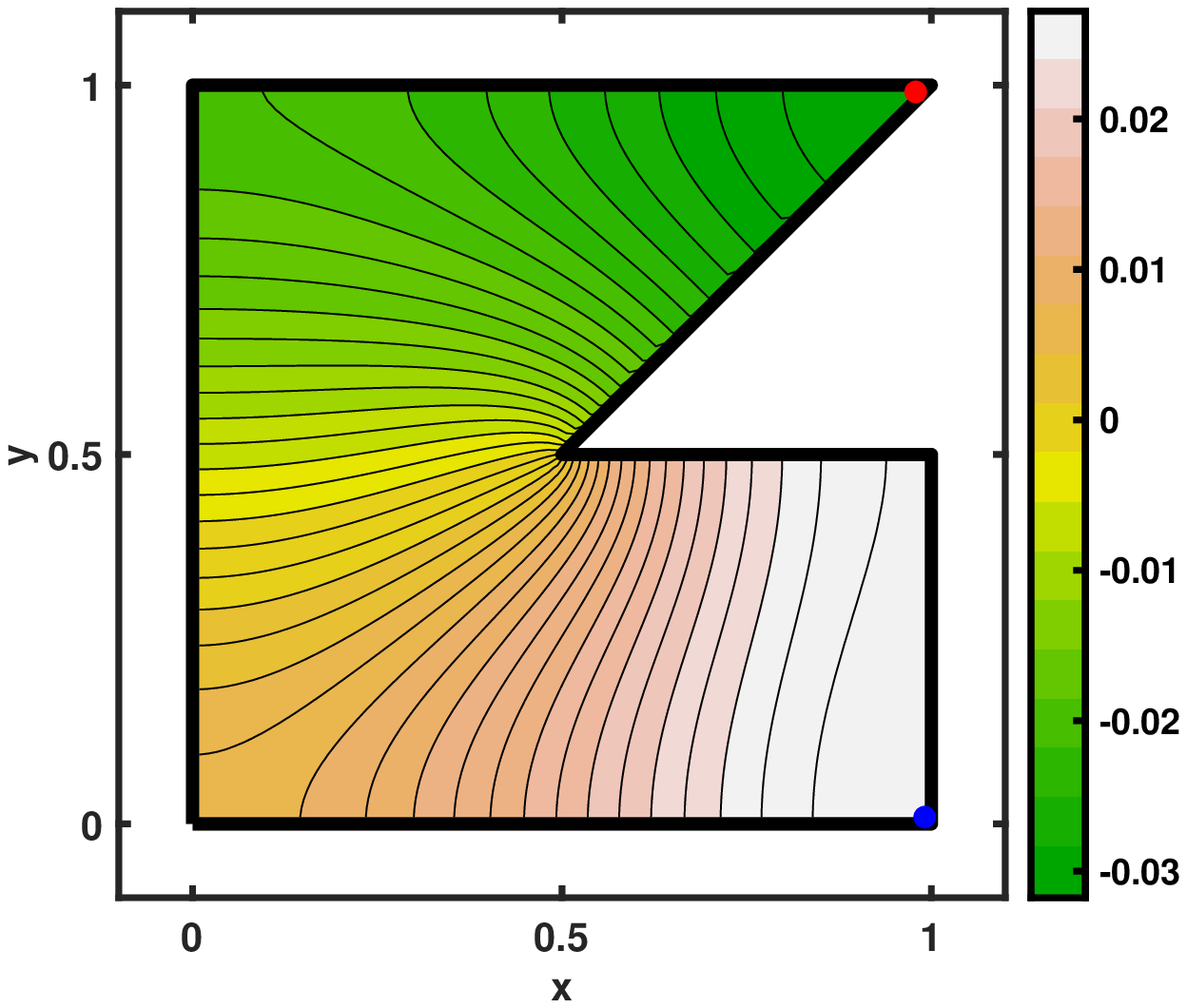}
}
\caption{\label{leigen} The eigenfunction corresponding to the first non-trivial interior Neumann eigenvalues for the L-shaped domain and its two variants
with corresponding non-trivial 
interior Neumann eigenvalue $2.429\,474$, $2.725\,559$, and $2.207\,276$, respectively.}
\end{figure}

All three eigenfunctions including the location of the maximal and minimal value are shown in Figure \ref{morecomplex}.
As we can see again, the extreme values are obtained on the boundary. The same is true for the L-shaped domain that is given by $\llcorner=[0,1]^2-[0.5,1]^2$. 
We obtain the first non-trivial interior Neumann eigenvalue $2.429\,474$ 
and compare it with the well-known value (see \cite{gilette} for the approximation $2\sqrt{1.475\,621\,845}\approx 2.429\,504$). As we can see in Table \ref{lshape} 
the absolute error decreases dramatically since we have less regularity of the solution at the corner. We are only able to achieve four digits accuracy with a convergence rate that seems to be $4/3$. 
\begin{table}[!ht]  
\caption{\label{lshape} Absolute error and estimated order of convergence of the first non-trivial interior Neumann eigenvalue for an L-shaped domain using different number of faces and collocation points.}
\begin{indented}
 \item[]\begin{tabular}{@{}rrrl}
 \br
  $n_f$  & $n_c$ & abs. error $E_{n_f}^{\llcorner}$ & $\mathrm{EOC}^{\llcorner}$ \\
  \mr
     6 &   18 & $7.0712_{-3}  $ &       \\
    12 &   36 & $3.1014_{-3}  $ & 1.1891\\
    24 &   72 & $1.1976_{-3}  $ & 1.3728\\
    48 &  144 & $4.7500_{-4}  $ & 1.3341\\
    96 &  288 & $1.8860_{-4}  $ & 1.3326\\
   192 &  576 & $7.4872_{-5}  $ & 1.3329\\
   384 & 1152 & $2.9724_{-5} $ & 1.3327\\
  \br
 \end{tabular}
 \end{indented}
\end{table}
In Figure \ref{leigen} we show the corresponding eigenfunction.
We also tried two different domains $\llcorner_2$ and $\llcorner_3$ as shown in Figure \ref{lshape} b) and c). 
Using the same parameters as for the L-shaped domain $\llcorner$, we
obtain the first non-trivial Neumann eigenvalues $2.725\,559$ and $2.207\,276$, respectively. Interestingly, the 
approximate convergence rates are $2.8432$, $2.0744$, $1,6656$, $1.6411$, $1.7140$, and $2.0115$ for $\llcorner_2$
and $1.0716$, $1.1645$, $1.1745$, $1.2202$, $1.3329$, and $1.6816$ for $\llcorner_3$.
That is why we concentrate on smooth boundaries.

In sum, we are not able to construct a simply-connected non-convex domain that fails the hot spots conjecture. 
Next, we concentrate on non-simply-connected domains.
\subsection{Non-simply-connected domains}
Now, we consider an annulus with inner radius $R_1=1/2$ and outer radius $R_2=2$. For this domain, we can again compute a reference solution to arbitrary precision.
The non-trivial interior Neumann eigenvalues are obtained through the roots of
\[J_n'(R_1x)Y_n'(R_2x)-J_n'(R_2x)Y_n'(R_1x)=0\,,\quad n=0,1,2,\ldots\]
where $Y_n'$ denotes the derivative of the second kind Bessel function of order $n$.
(see for example \cite[Equation 4.16]{tsai}). The first two roots are obtained with the Maple command
\begin{verbatim}
 restart; Digits:=16:
 Jp:=unapply(diff(BesselJ(1,x),x),x):
 Yp:=unapply(diff(BesselY(1,x),x),x):
 fsolve(Jp(x/2)*Yp(2*x)-Jp(2*x)*Yp(x/2),x=1);
 Jp:=unapply(diff(BesselJ(2,x),x),x):
 Yp:=unapply(diff(BesselY(2,x),x),x):
 fsolve(Jp(x/2)*Yp(2*x)-Jp(2*x)*Yp(x/2),x=1);
\end{verbatim}
The two smallest roots are approximately given by 
\[ 0.822\,252\,688\,623\,884\quad \text{and}\quad 1.504\,647\,782\,189\,479\,, \]
respectively. They both have multiplicity two. Again, we show in Table \ref{annulus} that our method is 
able to achieve ten digits accuracy with a cubic convergence order for the first two non-trivial interior Neumann eigenvalues for an annulus using the parameter $\mu=6/5$ for various number of faces. Note that twice 
the number of faces is needed due to the need of two boundary curves.
\begin{table}[!ht]  
\caption{\label{annulus} Absolute error and estimated order of convergence of the first and second non-trivial interior Neumann eigenvalue for an annulus with $R_1=1/2$ and $R_2=2$ 
using different number of faces and collocation points.}
\begin{indented}
 \item[]\begin{tabular}{@{}rrrlrl}
 \br
  $2n_f$  & $2n_c$ & absolute error $E_{n_f}^{(1)}$ & $\mathrm{EOC}^{(1)}$ & absolute error $E_{n_f}^{(2)}$ & $\mathrm{EOC}^{(2)}$\\
  \mr
    10 &   30 & $2.4700_{-3}$ &       & $6.1268_{-3}$ &   \\
    20 &   60 & $1.9490_{-4}$ & 3.6637& $5.4708_{-4}$ & 3.4853\\
    40 &  120 & $1.7741_{-5}$ & 3.4575& $5.7347_{-5}$ & 3.2540\\
    80 &  240 & $1.8725_{-6}$ & 3.2441& $6.6443_{-6}$ & 3.1095\\
   160 &  480 & $2.1655_{-7}$ & 3.1122& $8.0539_{-7}$ & 3.0443\\
   320 &  960 & $2.6159_{-8}$ & 3.0493& $9.9900_{-8}$ & 3.0111\\
   640 & 1920 & $3.2351_{-9}$ & 3.0154& $1.2988_{-8}$ & 2.9433\\
  \br
 \end{tabular}
 \end{indented}
\end{table}

\begin{figure}[!ht]
\subfigure[First eigenfunction of an annulus]{
\includegraphics[width=0.31\textwidth]{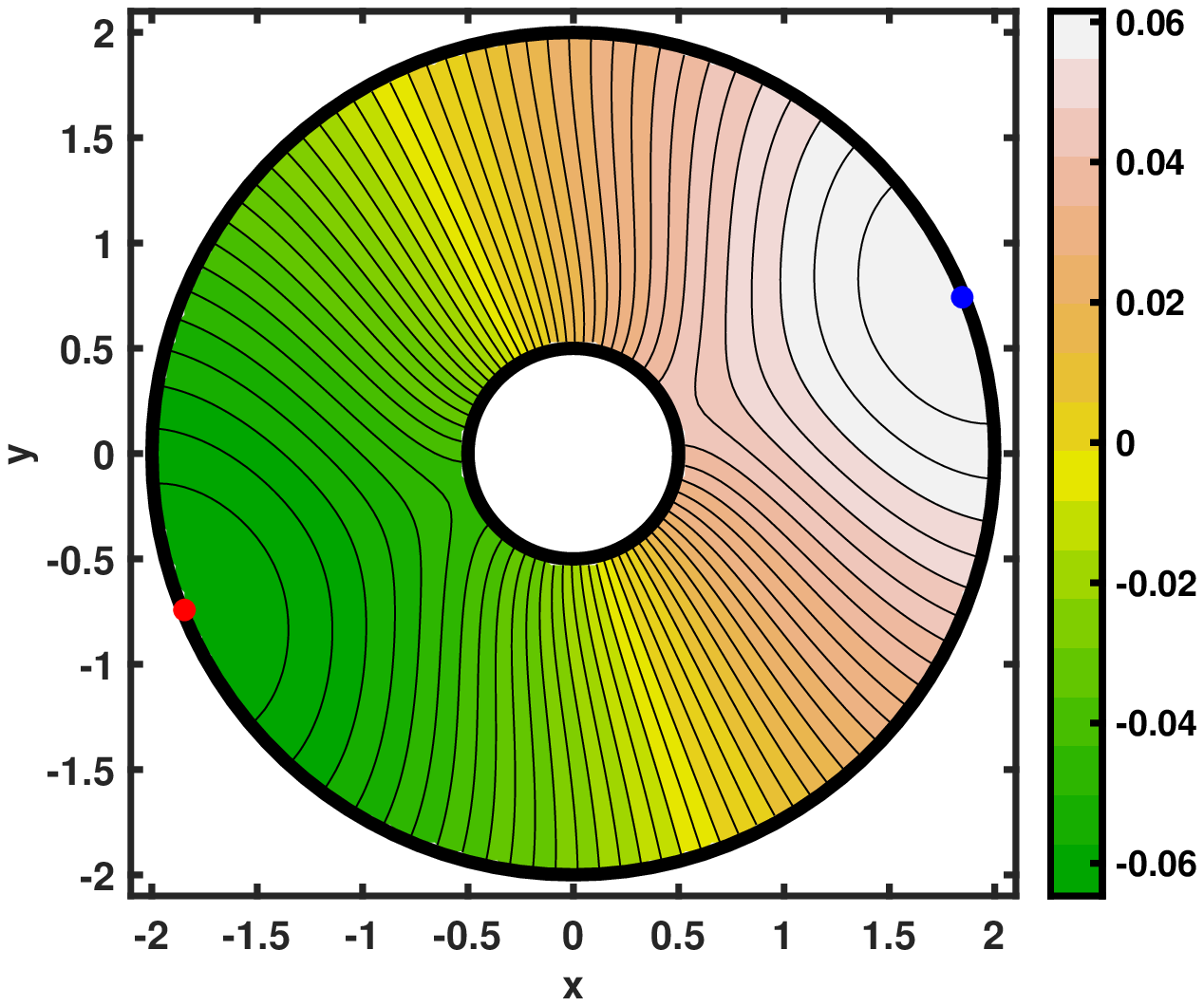}
}
\subfigure[Second eigenfunction of an annulus]{
\includegraphics[width=0.31\textwidth]{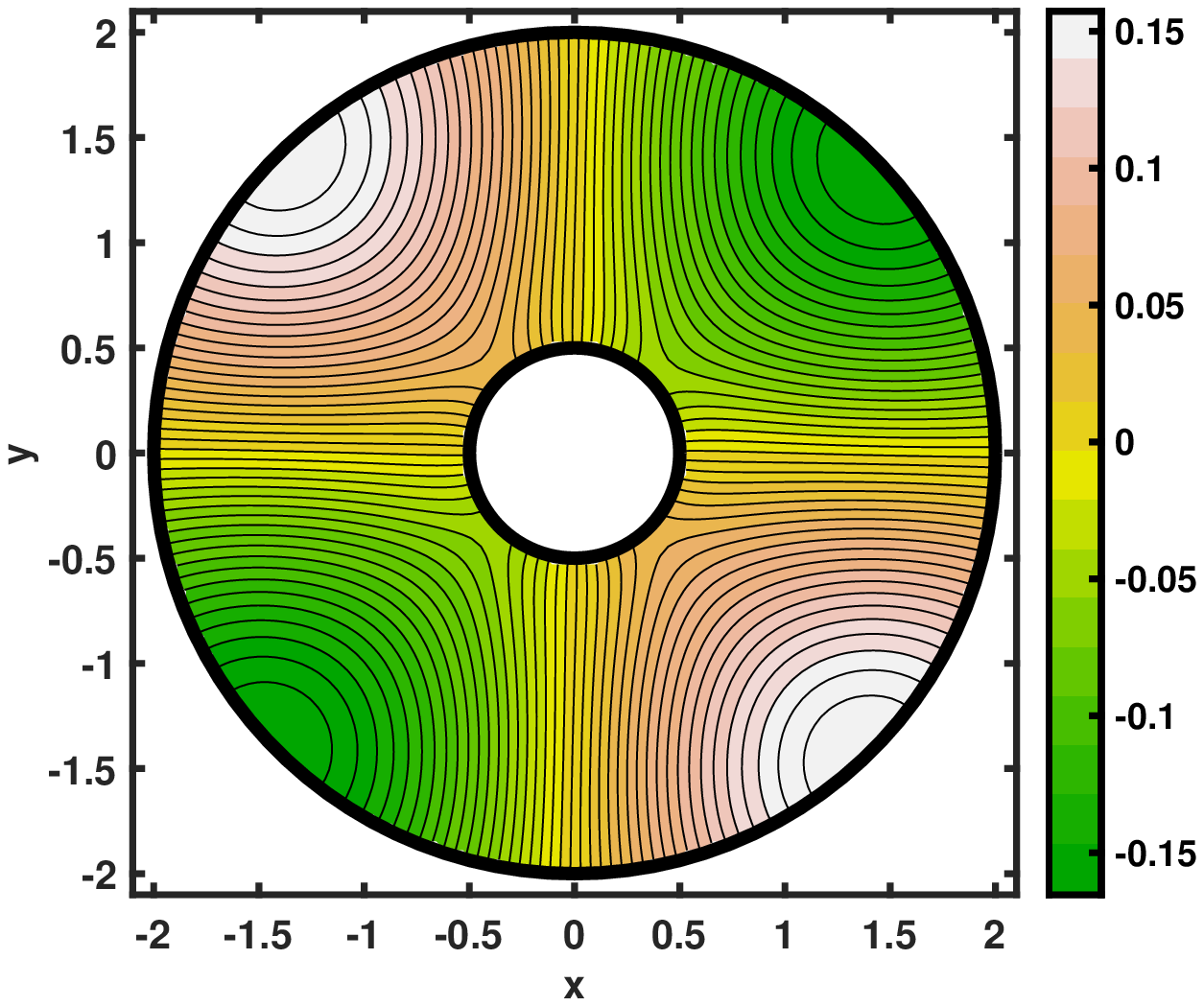}
}
\subfigure[Third eigenfunction of an annulus]{
\includegraphics[width=0.31\textwidth]{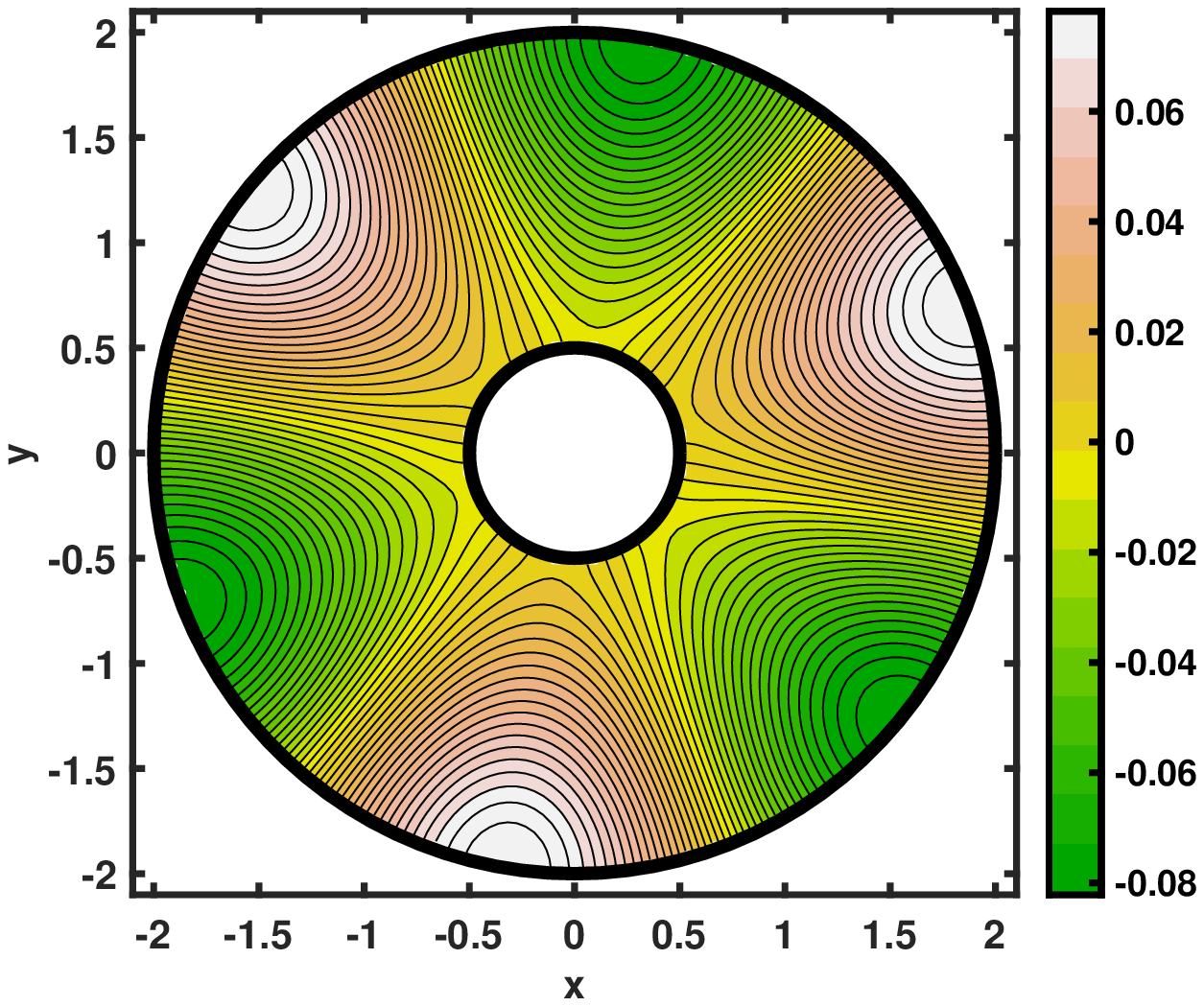}
}
\caption{\label{annulus2} The first three eigenfunctions corresponding to the first three non-trivial interior Neumann eigenvalues $0.822\,253$, $1.504\,648$, and $2.096\,773$ for an annulus.}
\end{figure}

In Figure \ref{annulus2} we show the first three eigenfunctions for $\square_{2.1}$ corresponding to the non-trivial interior Neumann 
eigenvalues $0.822\,253$, $1.504\,648$, and $2.096\,773$. To compute the third eigenvalue which has multiplicity two, we used $\mu=2$. For the first eigenfunction plot, we also added the extreme values.

\begin{figure}[!ht]
\subfigure[First eigenfunction of $A_1$]{
\includegraphics[width=0.31\textwidth]{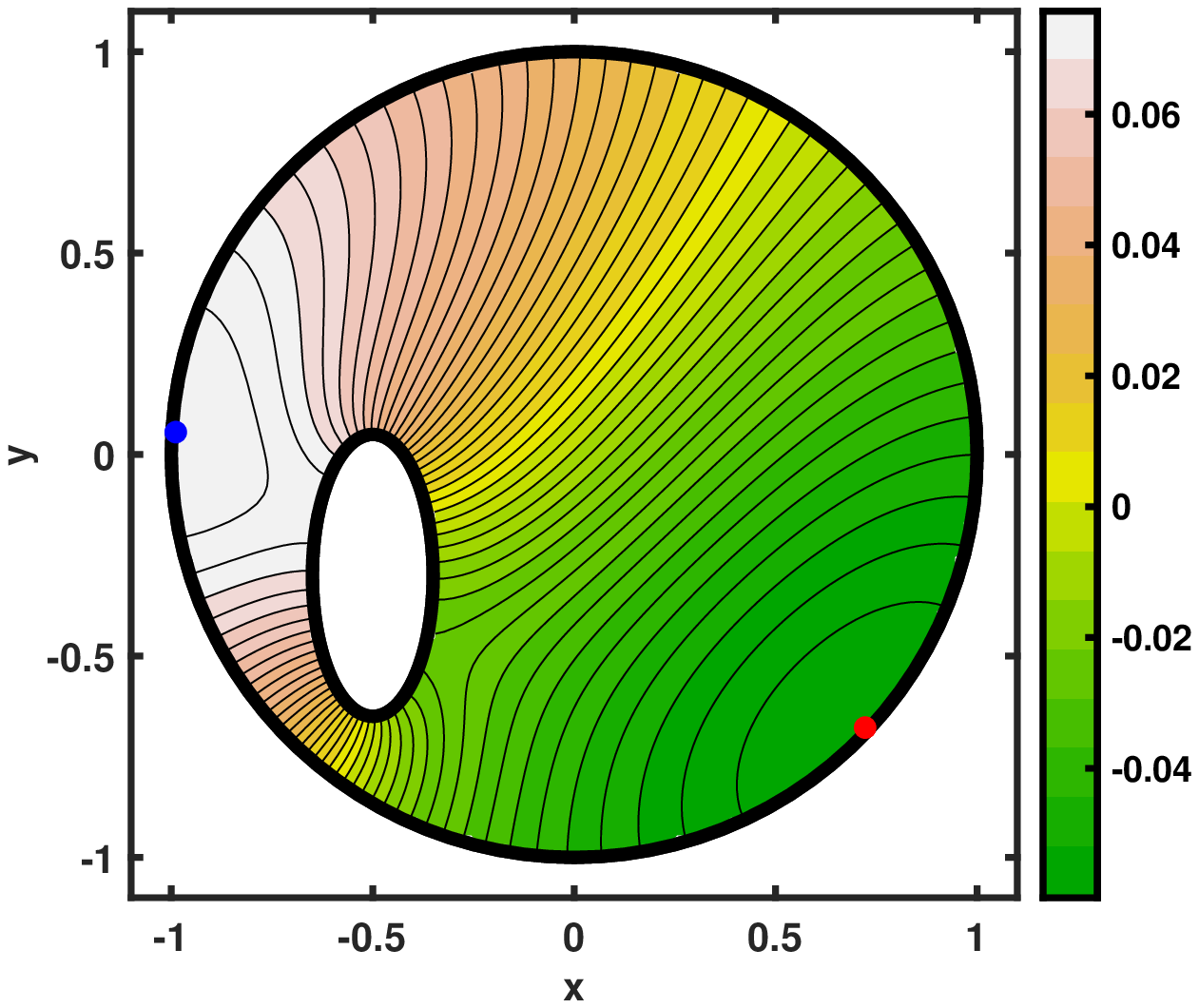}
}
\subfigure[First eigenfunction of $A_2$]{
\includegraphics[width=0.31\textwidth]{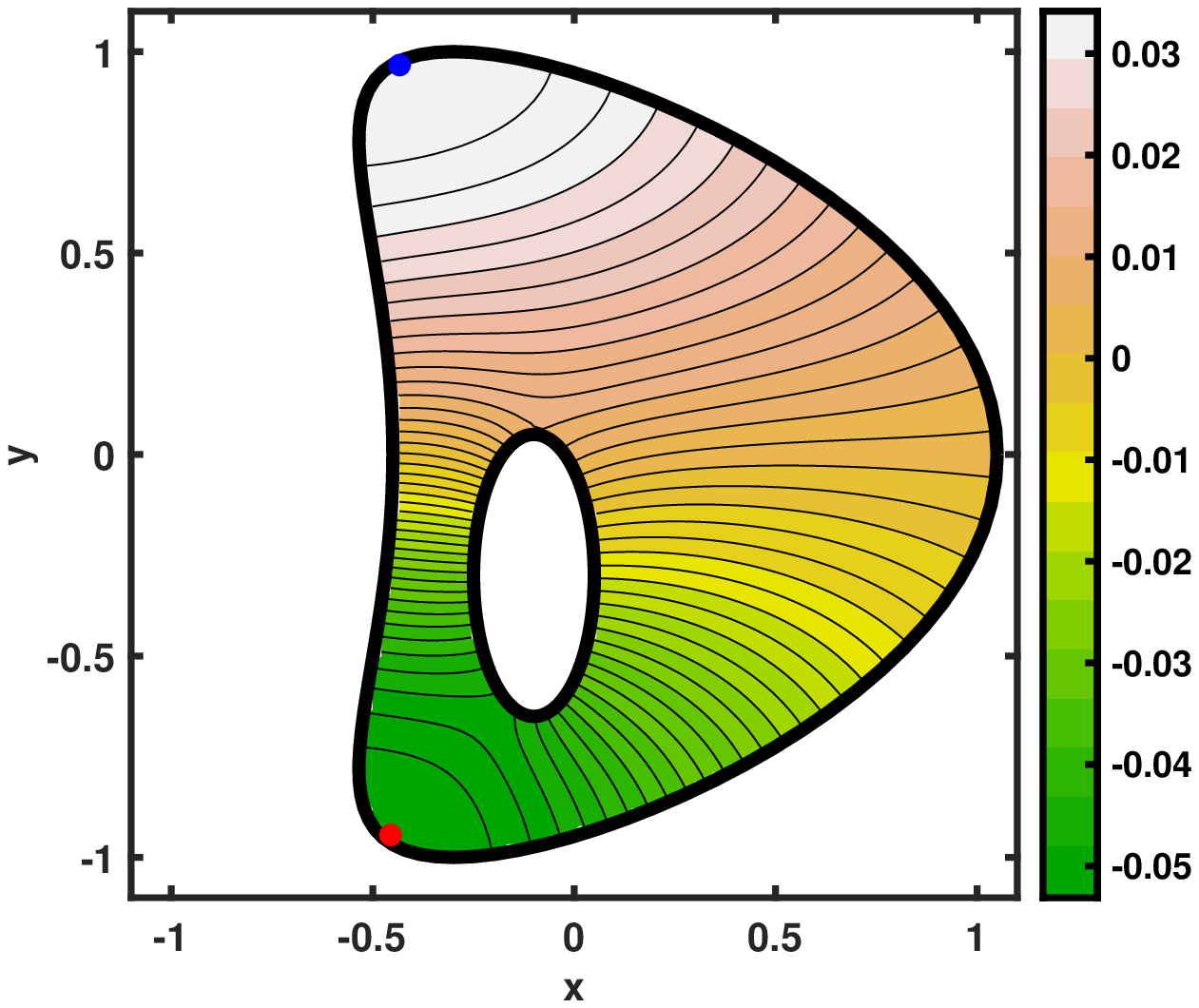}
}
\subfigure[First eigenfunction of $A_3$]{
\includegraphics[width=0.31\textwidth]{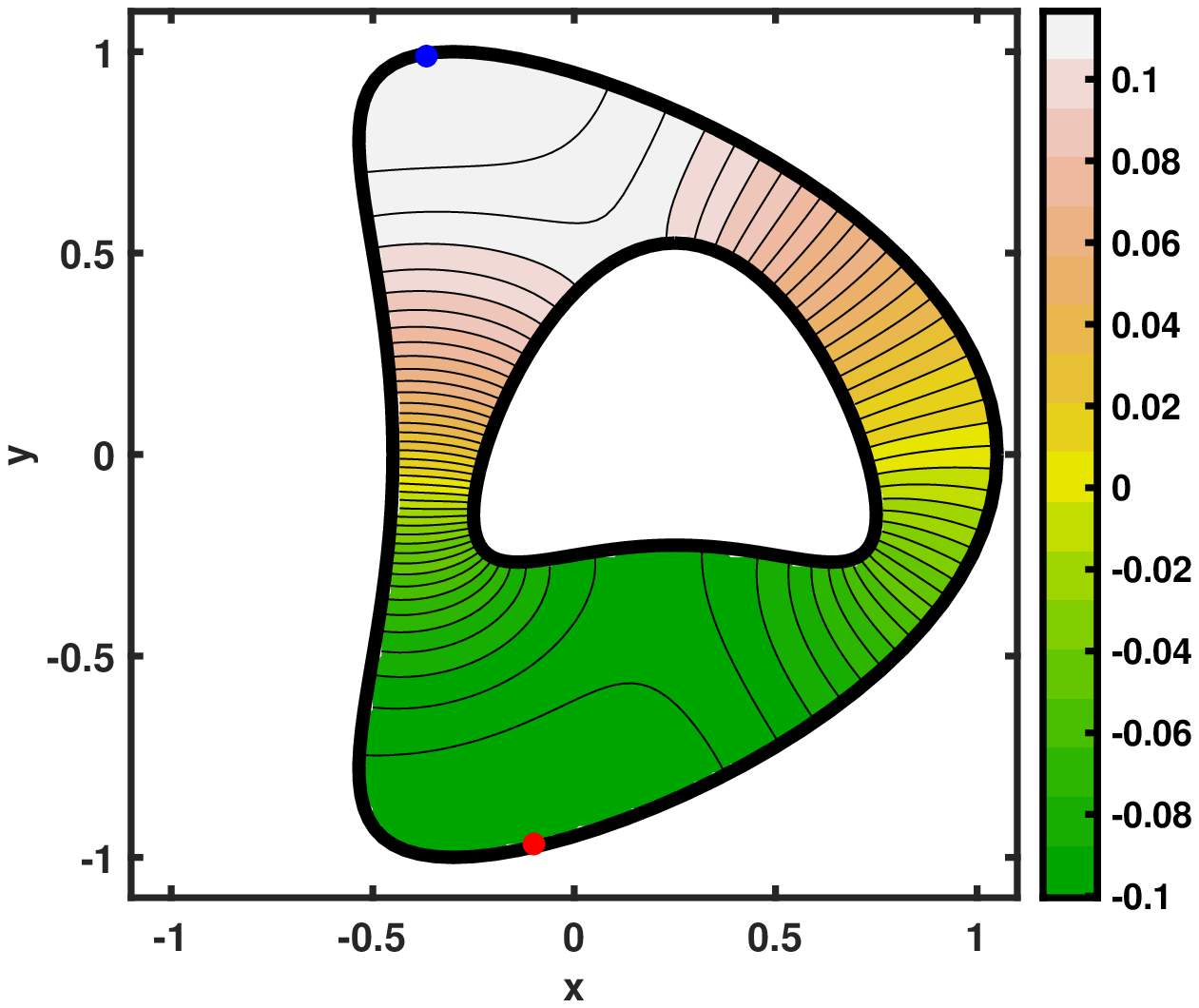}
}
\caption{\label{others} The eigenfunctions corresponding to the first non-trivial interior Neumann eigenvalues $1.662\,873$, $1.651\,571$, and $1.171\,590$ for the domains $A_1$, $A_2$, and $A_3$.
}
\end{figure}

We can see that the extreme values are again obtained on the boundary of the annulus. Now, we consider more complex non-simply connected domains. The first domain $A_1$ is given by a unit circle centered at the origin
removing an ellipse centered at $(-0.5,-0.3)^\top$ with semi-axis $0.15$ and $0.35$. We obtain the first non-trivial interior Neumann eigenvalue $1.662\,873$ and the eigenfunction within $\square_{1.1}$. The second domain $A_2$ is given by $D_3$ 
removing an ellipse centered at $(-0.1,-0.3)^\top$ with semi-axis $0.15$ and $0.35$. We get the first non-trivial interior Neumann eigenvalue $1.651\,571$ and the eigenfunction within $\square_{1.1}$. The third domain $A_3$ is given by $D_3$ 
removing a 90 degree counter-clockwise rotated version of $D_3$ scaled by $1/2$. We obtain the first non-trivial interior Neumann eigenvalue $1.171\,590$ and the eigenfunction within $\square_{1.1}$.
We used the parameter $\mu=1.5$ for all three domains. The eigenfunctions for $A_1$, $A_2$, and $A_3$ including the extreme values 
are illustrated in Figure \ref{others}.

As we can see again, the extreme values are obtained on the boundary. We also tried many other similar geometries with different kinds of holes and obtained similar results. The hot spots conjecture seems to hold.

Finally, we concentrate on a more complex geometry inspired by the work of Burdzy \cite[Figure 1]{burdzy}. He proved that there exists a bounded planar domain with one hole that fails the hot spots conjecture. 
However, the description of the proposed bounded domain with one hole is very technical and it is difficult to implement his domain to verify his theoretical result. His domain's boundary has many corners which would 
complicate the explicit construction of the boundary. Additionally, the domains are very thin.
Note that no numerical result supports his theoretical result.
In fact, up-to-date no numerical results have been shown yet.

We try to close this gap. We construct less complex domains that 
fail the hot spots conjecture and show it numerically. Further, we show the exact location of the extreme values.

The `teether' domain depicted in Figure \ref{complex} is constructed as follows: 
\begin{figure}[!ht]
\centering
\includegraphics[width=0.31\textwidth]{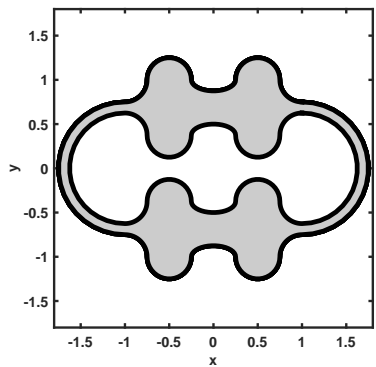}
\caption{\label{complex} A more advance bounded domain with a hole inspired by the work of Burdzy \cite[Figure 1]{burdzy}.}
\end{figure}
Let $E(x,y,a,b,t_1,t_2)$ be the ellipse centered at $(x,y)^\top$ with semi-axis $a$ and $b$ constructed for $\phi\in [t_1,t_2)$ using
the parametrization $(x+a\cos(\phi),y+b\sin(\phi))^\top$. Note that we also allow $t_1>t_2$ to guarantee the needed 
orientation of the curve.
The first half of the outer boundary is given by the pieces $E(4,4,1,1,-\pi/2,-\pi)$, $E(2,4,1,1,0,\pi)$, $E(0,4,1,1/2,0,-\pi)$, 
$E(-2,4,1,1,0,\pi)$, $E(-4,4,1,1,0,-\pi/2)$, and $E(-4,0,3,3,\pi/2,3\pi/2)$. Rotating this half by $\pi$ yields the second 
half of the outer
boundary. The first half of the inner boundary is given by the pieces $E(4,3/2,1,1,\pi/2,\pi)$, $E(2,3/2,1,1,0,-\pi)$, $E(0,3/2,1,1/2,0,\pi)$, 
$E(-2,3/2,1,1,0,-\pi)$, $E(-4,3/2,1,1,0,\pi/2)$, and $E(-4,0,5/2,5/2,\pi/2,3\pi/2)$. Rotating this half by $\pi$ yields the 
second half of the inner
boundary. Next, the orientation of the inner boundary is reversed. 
Finally, all coordinates of the boundary are multiplied with $1/4$. This yields our `teether' domain $C_1$.

\begin{figure}[!ht]
\subfigure[First eigenfunction of $C_1$]{
\includegraphics[width=0.31\textwidth]{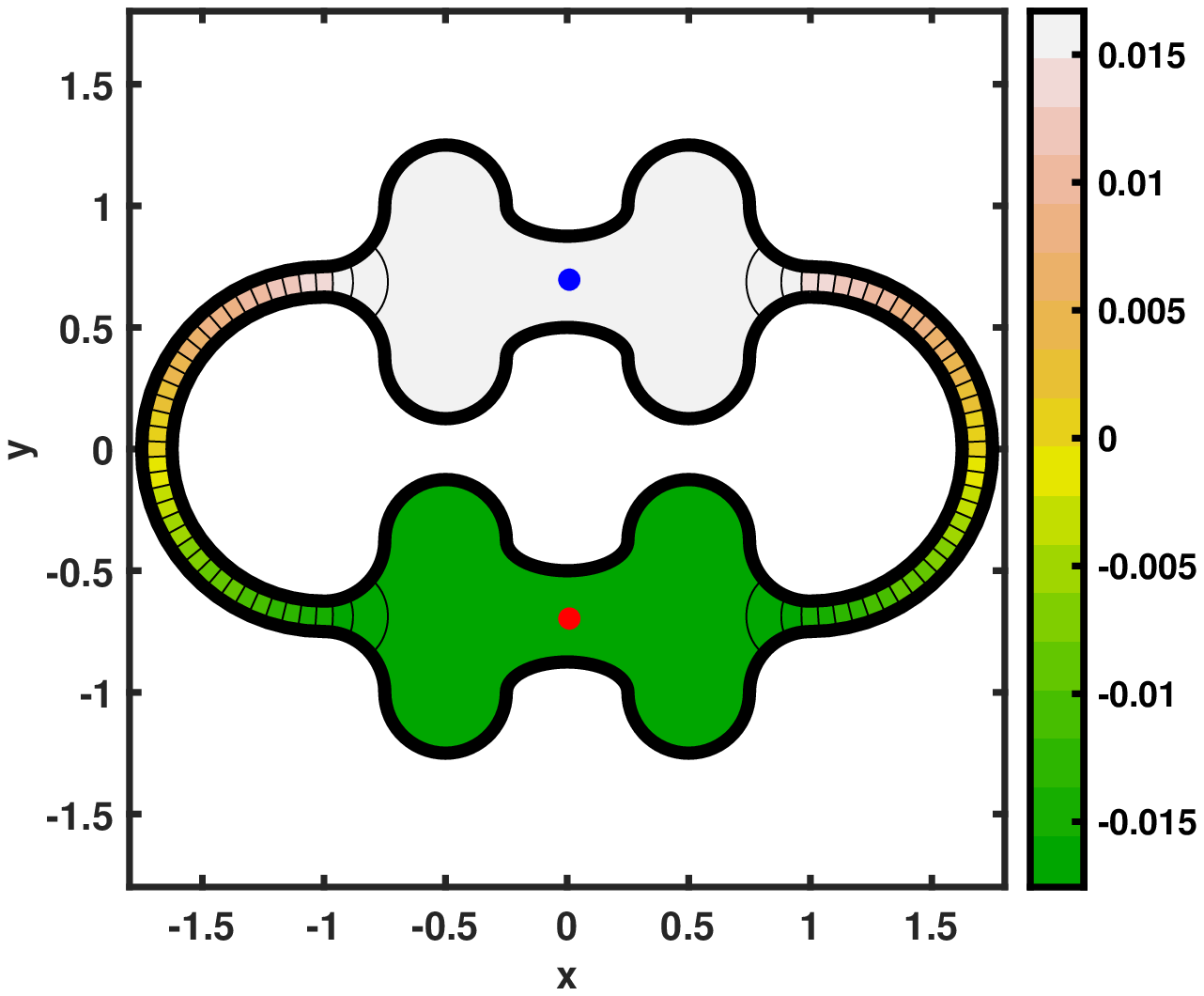}
}
\subfigure[Zoom around the maximum]{
\includegraphics[width=0.31\textwidth]{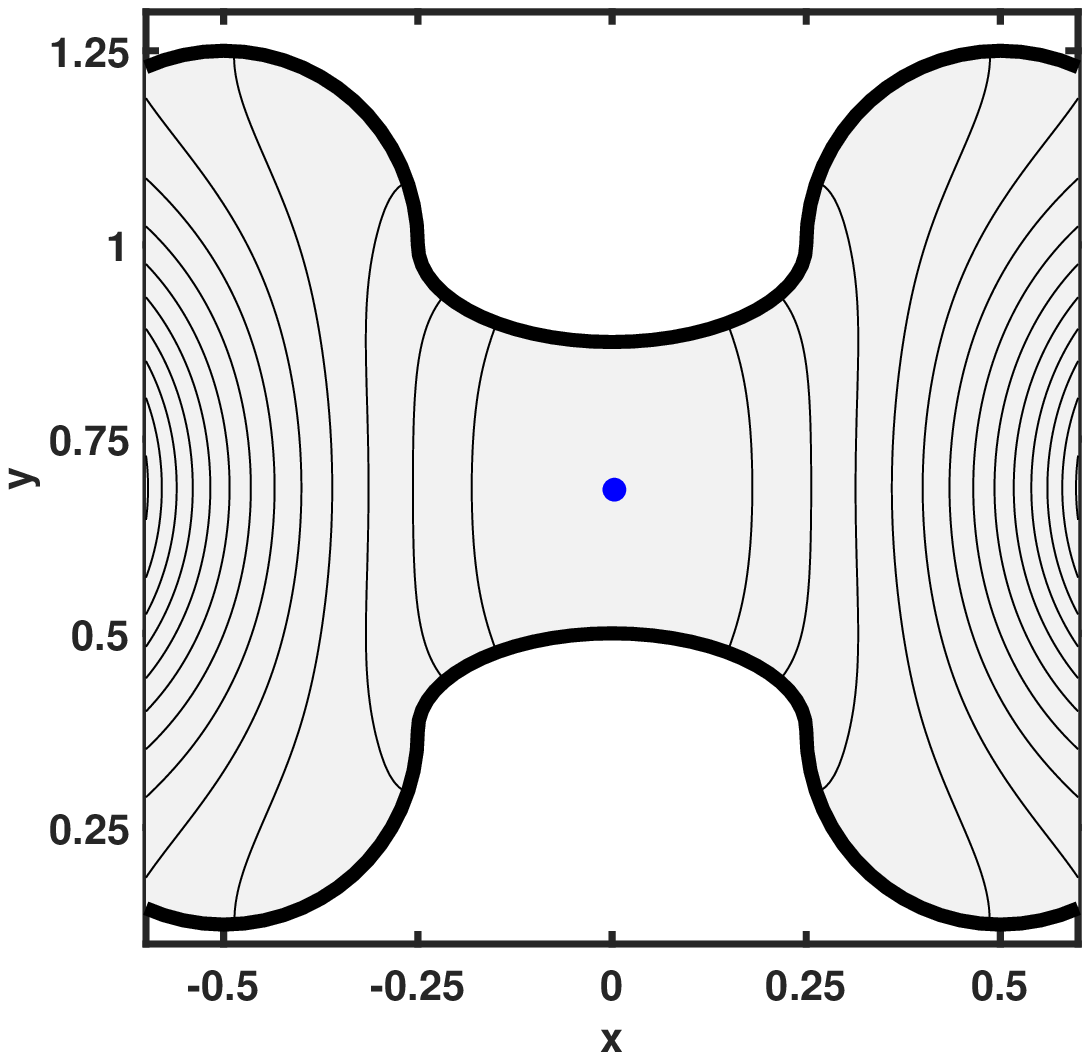}
}
\subfigure[Zoom around the minimum]{
\includegraphics[width=0.31\textwidth]{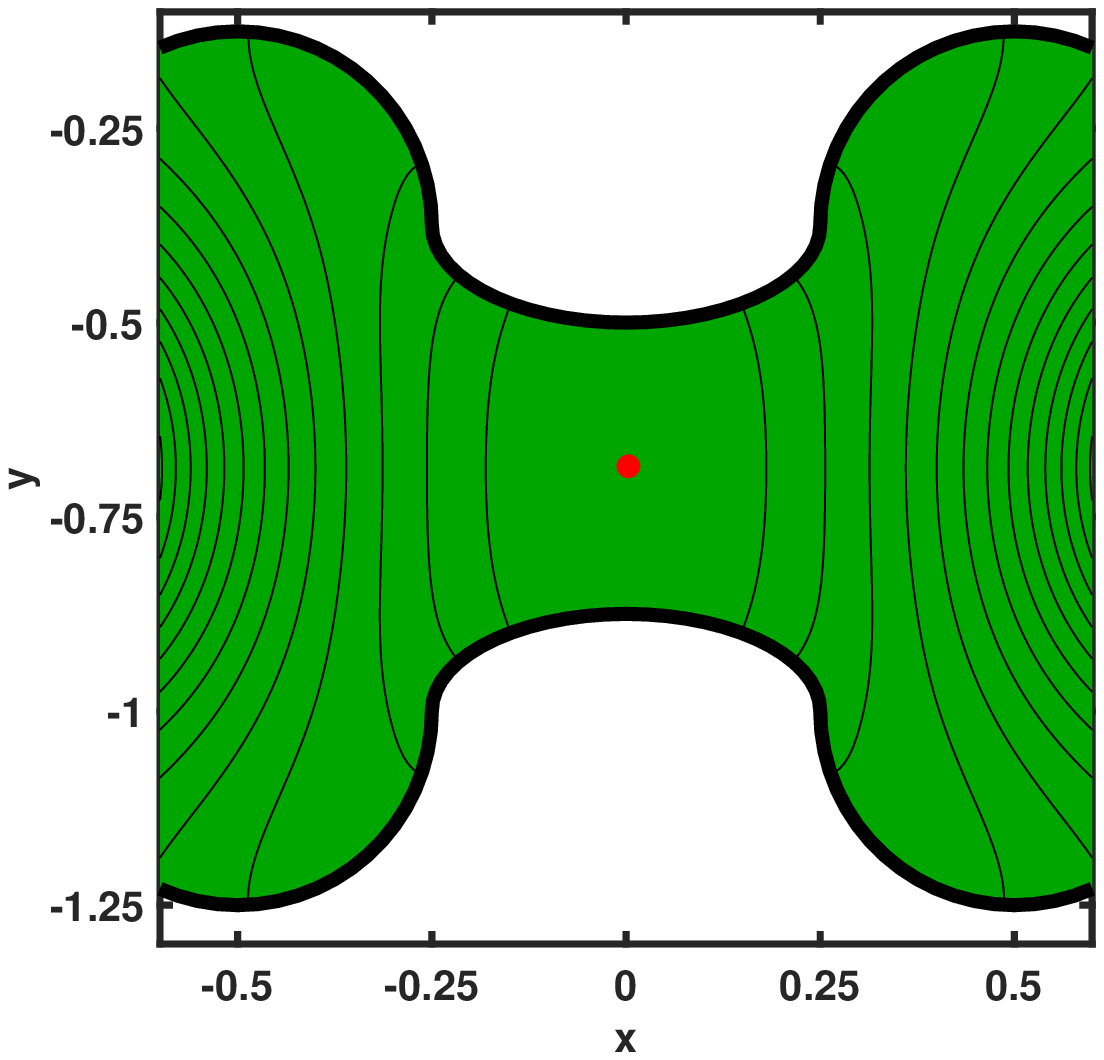}
}
\caption{\label{evcomplex} The eigenfunction and its zoomed version around the maximum and minimum corresponding to the first non-trivial interior Neumann eigenvalue $0.370\,708$ for the teether domain $C_1$.
}
\end{figure}

We use the parameter $\mu=0.8$, $\ell=20$ and $\square_{1.8}$. We obtain the first non-trivial interior Neumann eigenvalue $0.370\,708$. The corresponding eigenfunction including is extreme values is shown in Figure \ref{evcomplex}. 
Since the extreme values lie in very flat plateaus, we additionally show zoomed versions around the extreme values to better see that they are located inside the domain.

This shows that we are able to show numerically that there exists a bounded domain with one hole that fails the hot spots conjecture.

Next, we investigate some possible conditions needed to construct an example that fails the hot spots conjecture. With one bump it was not possible to obtain the extreme values inside the domain, say $C_2$.
We use the previous example, and remove one of the bump and its mirror version in the upper and lower part of the teether domain. We obtain the first non-trivial interior Neumann eigenvalue 
$0.534\,605$ and the corresponding eigenfunction within $\square_{1.6}$.
\begin{figure}[!ht] 
\subfigure[First eigenfunction of $C_2$]{
\includegraphics[width=0.31\textwidth]{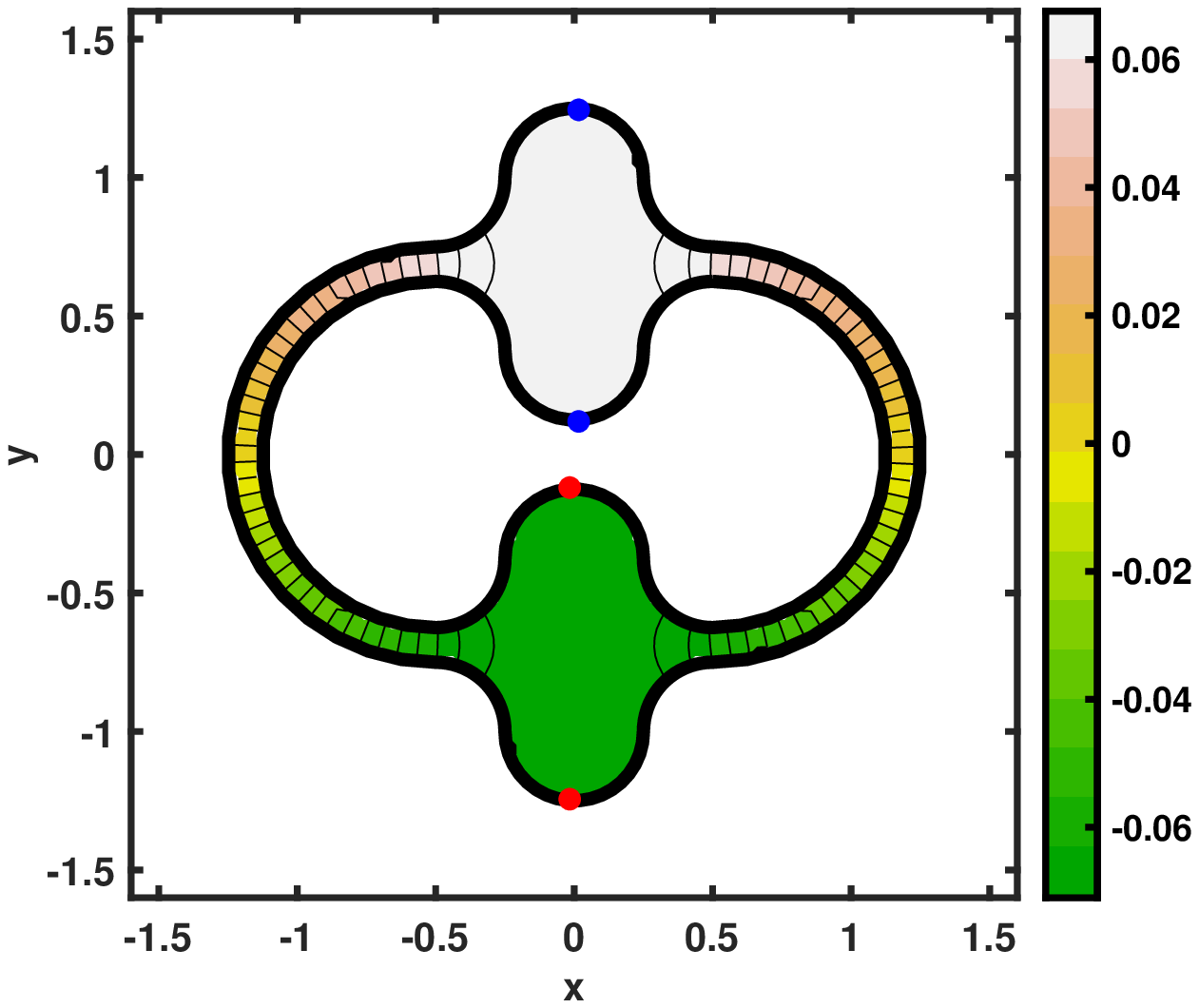}
}
\subfigure[Zoom around the maximum]{
\includegraphics[width=0.31\textwidth]{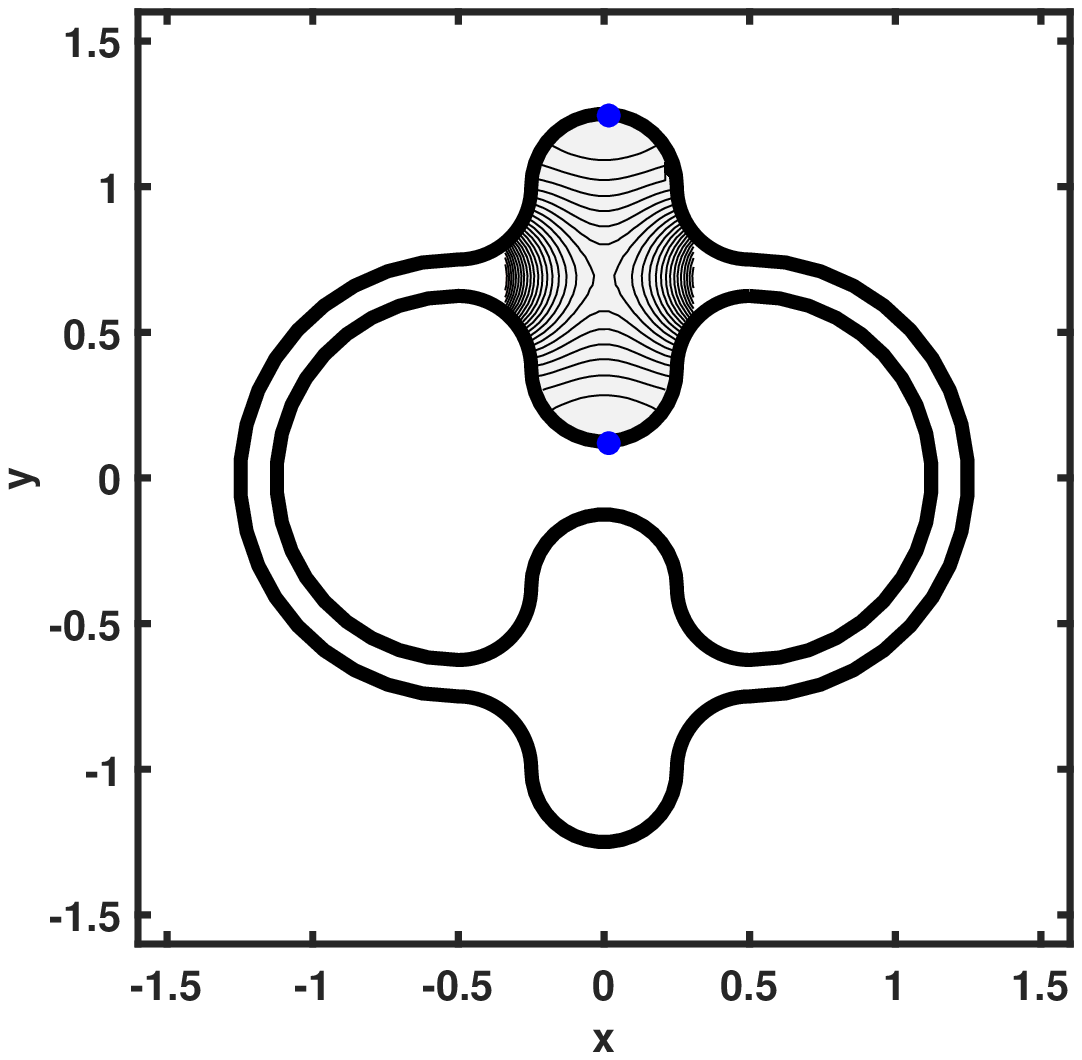}
}
\caption{\label{evcomplex2} The eigenfunction and its zoomed version around the maximum corresponding to the first non-trivial interior Neumann eigenvalue $0.534\,605$ for the domain $C_2$.
}
\end{figure}
As we can see again in Figure \ref{evcomplex2}, the values inside the bump area are very close to each other. The extreme values are attained on the boundary.

Right now, the proposed domain $C_1$ has two lines of symmetry. Now, we break the symmetry and show that we are still able to obtain a counter-example of the hot spots conjecture. 
We add a small amount of $0.0250$ to the semi axis of 
the ellipse that describes the upper left bump and obtain the domain $C_3$. The first non-trivial interior Neumann eigenvalue is $0.367\,496$.
\begin{figure}[!ht] 
\subfigure[First eigenfunction of $C_3$]{
\includegraphics[width=0.31\textwidth]{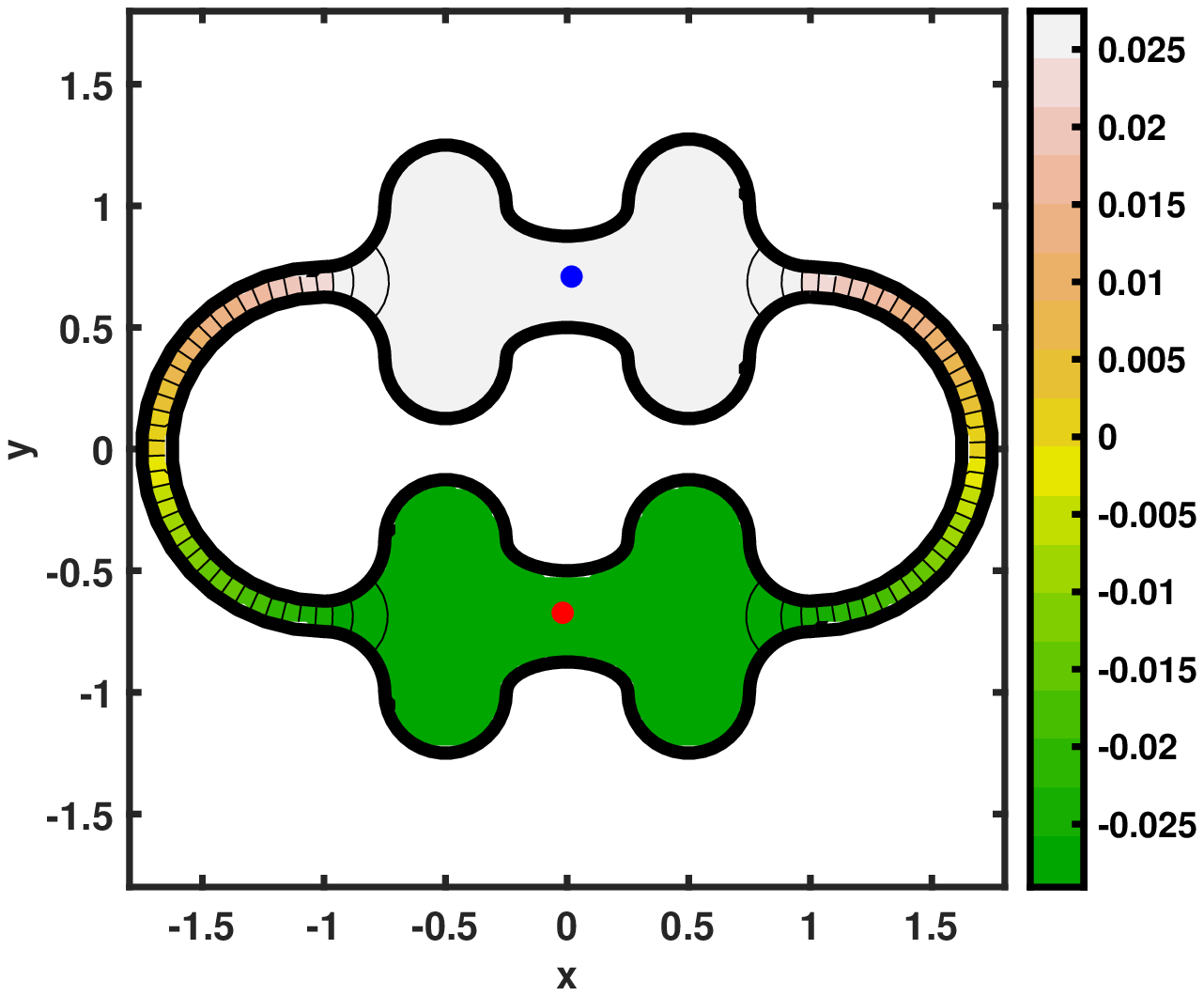}
}
\subfigure[First eigenfunction of $\widetilde{C}_3$]{ 
\includegraphics[width=0.31\textwidth]{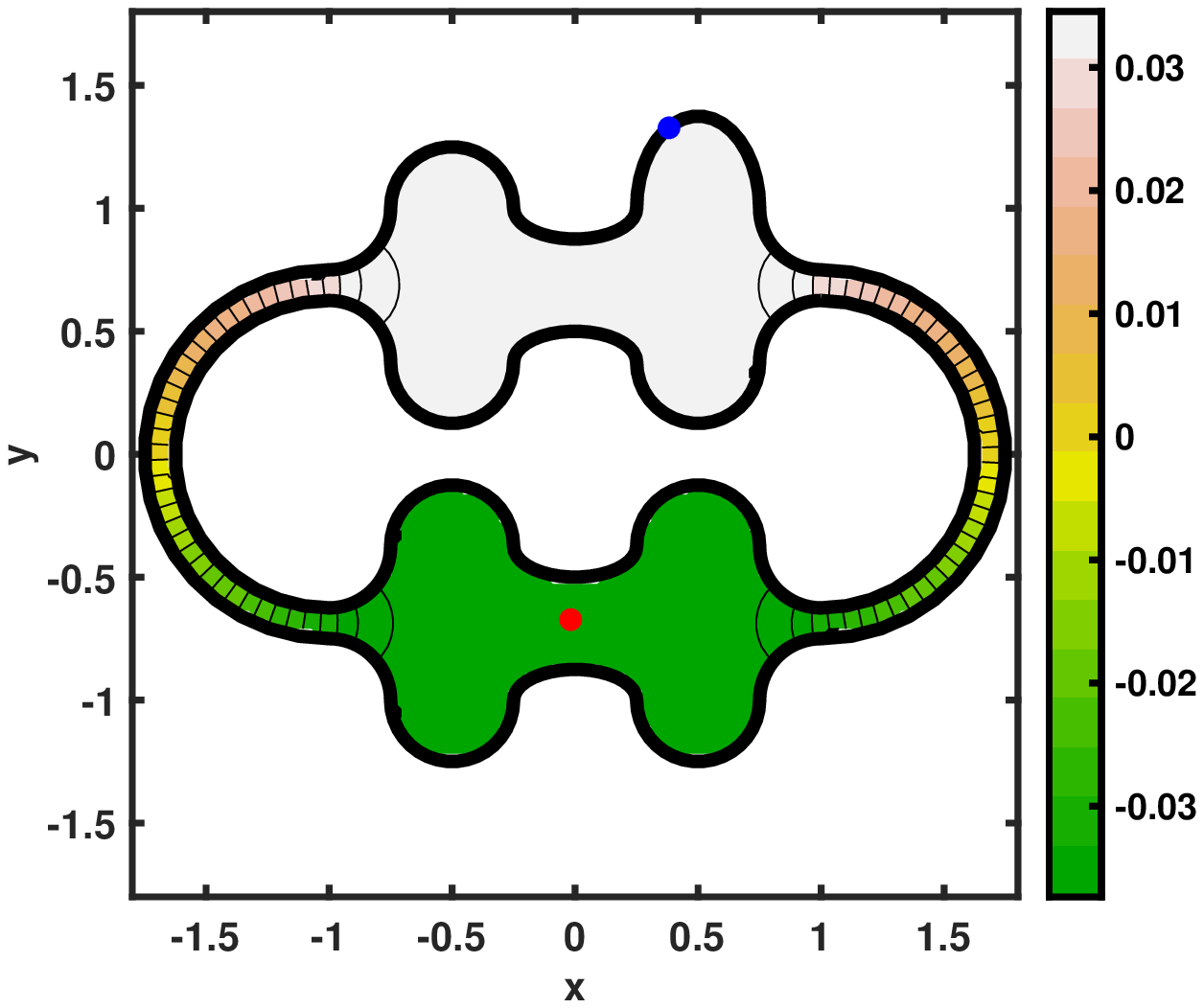}
}
\caption{\label{evcomplex3} The eigenfunction corresponding to the first non-trivial interior Neumann eigenvalue $0.367\,496$ and $0.370\,054$ for the domains $C_3$ and $\widetilde{C}_3$.
}
\end{figure}
As we can see in Figure \ref{evcomplex3}, the location of the hot spots of the eigenfunction within $\square_{1.8}$ are slightly changed, too, but they remain inside of $C_3$. Increasing the value from $0.0250$ to $0.1255$ 
will shift the maximal value from inside the domain $\widetilde{C}_3$ to the boundary while the minimum stays inside the domain. We obtain the eigenvalue $0.370\,054$.

\begin{figure}[!ht] 
\subfigure[First eigenfunction of $C_4 $]{
\includegraphics[width=0.31\textwidth]{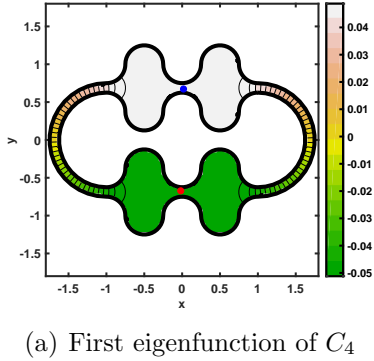}
}
\caption{\label{evcomplex4} The eigenfunction corresponding to the first non-trivial interior Neumann eigenvalue $0.384\,715$ for the domains $C_4$.
}
\end{figure}

Now, we show what happens if we make the gap between $E(0,4,1,1/2,0,-\pi)$ and $E(0,3/2,1,1/2,0,\pi)$ and its mirror counterpart smaller. We use $E(0,4,1,1,0,-\pi)$ and $E(0,3/2,1,1,0,\pi)$ instead to obtain $C_4$. We obtain the
first non-trivial interior Neumann eigenvalue $0.384\,715$ and its corresponding eigenfunction within $\square_{1.8}$ shown in Figure \ref{evcomplex4}.

Next, we show what happens if we make the bumps of the domain $C_4$ smaller. We construct the domain as before except that we introduce a new parameter $\delta>0$. The first half of the outer boundary is given by the pieces 
$E(4,3+\delta,1,\delta,-\pi/2,-\pi)$, $E(2,3+\delta,1,\delta,0,\pi)$, $E(0,3+\delta,1,\delta,0,-\pi)$, 
$E(-2,3+\delta,1,\delta,0,\pi)$, $E(-4,3+\delta,1,\delta,0,-\pi/2)$, and $E(-4,0,3,3,\pi/2,3\pi/2)$. Rotating this half by $\pi$ yields the second 
half of the outer
boundary. The first half of the inner boundary is given by the pieces $E(4,5/2-\delta,1,\delta,\pi/2,\pi)$, $E(2,5/2-\delta,1,\delta,0,-\pi)$, $E(0,5/2-\delta,1,\delta,0,\pi)$, 
$E(-2,5/2-\delta,1,\delta,0,-\pi)$, $E(-4,5/2-\delta,1,\delta,0,\pi/2)$, and $E(-4,0,5/2,5/2,\pi/2,3\pi/2)$. Rotating this half by $\pi$ yields the 
second half of the inner
boundary. Next, the orientation of the inner boundary is reversed. 
Finally, all coordinates of the boundary are multiplied with $1/4$. Note that $\delta=1$ yields the domain $C_4$.

\begin{figure}[!ht] 
\subfigure[First eigenfunction of $\tilde{C}_4$]{
\includegraphics[width=0.31\textwidth]{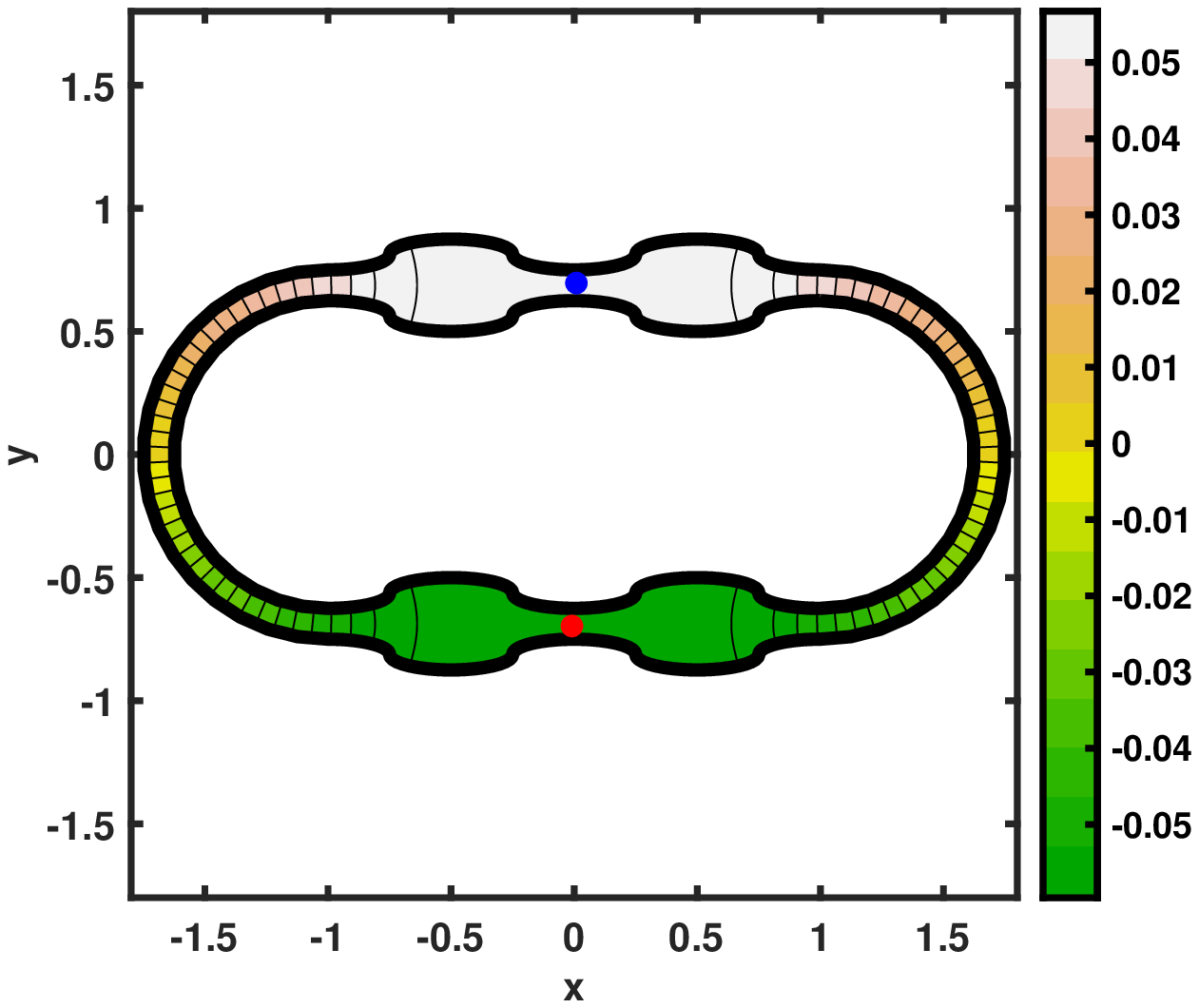}
}
\subfigure[First eigenfunction of $\hat{C}_4$]{
\includegraphics[width=0.31\textwidth]{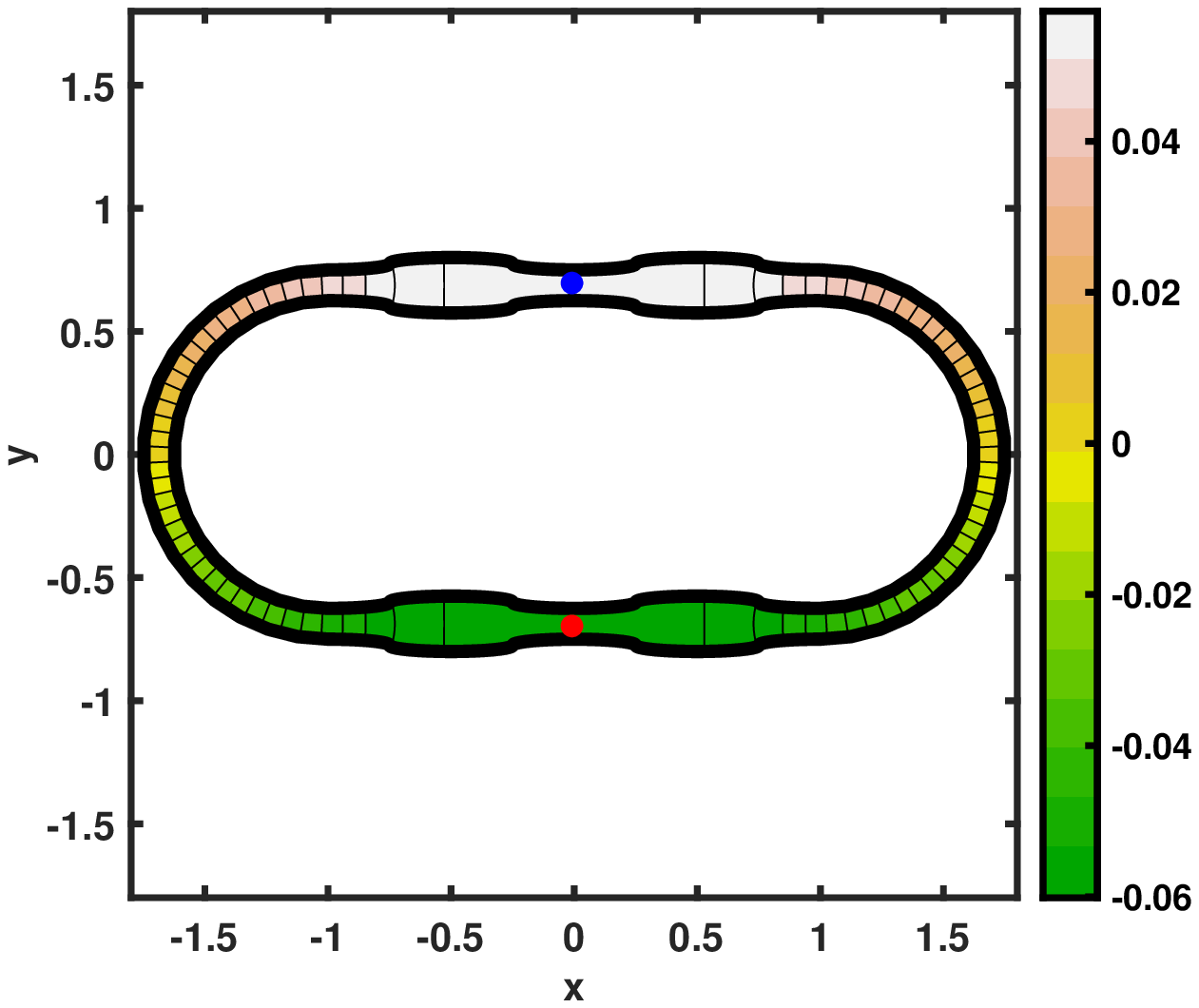}
}
\subfigure[First eigenfunction of $\bar{C}_4$]{
\includegraphics[width=0.31\textwidth]{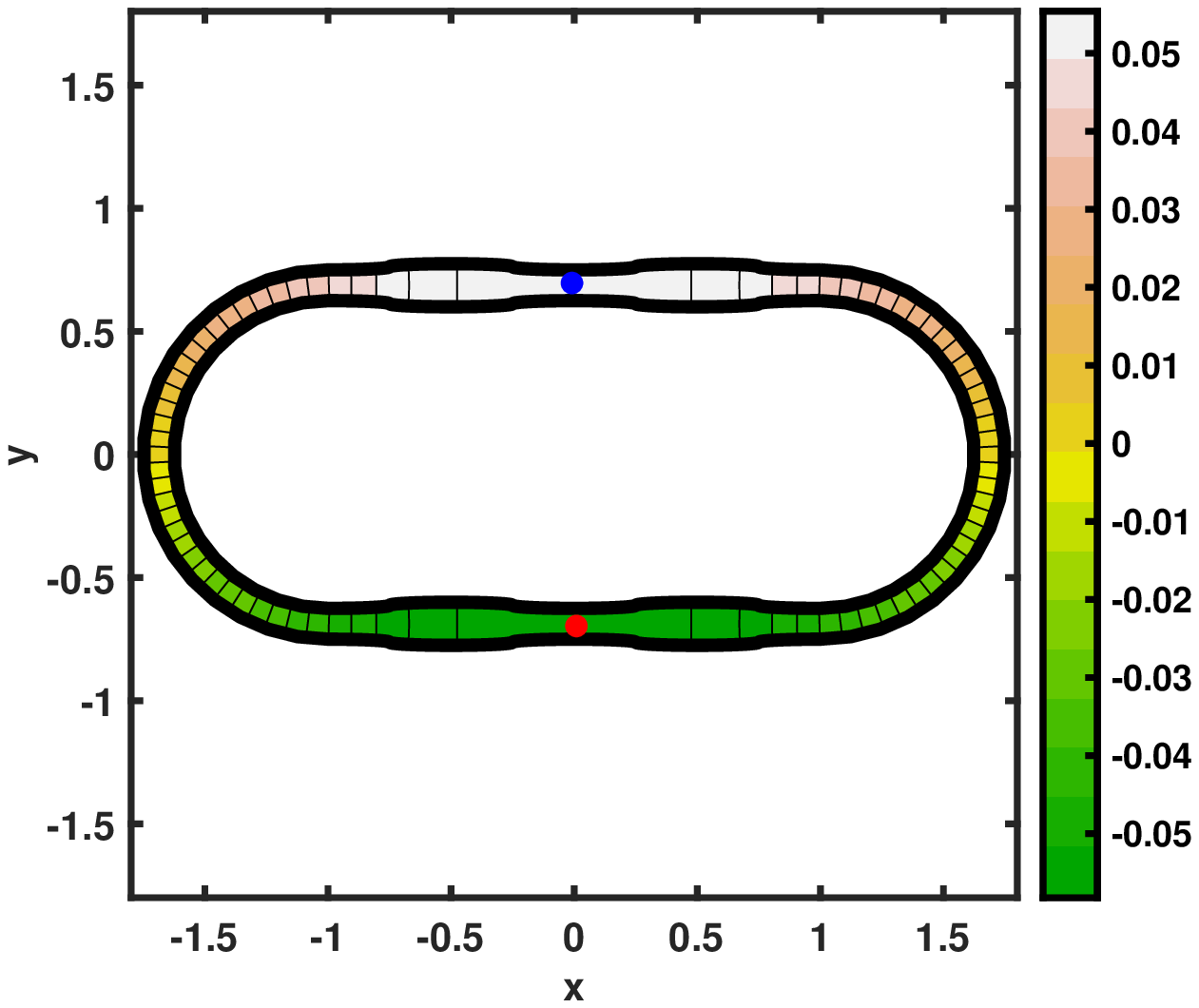}
}
\caption{\label{evcomplex6} The eigenfunction corresponding to the first non-trivial interior Neumann eigenvalue $0.578\,402$, $0.668\,373$, and $0.708\,633$ for the domains $\tilde{C}_4$, $\hat{C}_4$ and $\bar{C}_4$, respectively.
}
\end{figure}

We obtain the following results within $\square_{1.8}$ shown in Figure \ref{evcomplex6} for $\delta=1/4$, $\delta=1/10$, and $\delta=1/20$.
Surprisingly, the extreme values stay inside of the domain $\tilde{C}_4$, $\hat{C}_4$, and $\bar{C}_4$. Also note that the first non-trivial interior Neumann eigenvalue changes drastically.

\begin{figure}[!ht]  
\subfigure[First eigenfunction of $C_5$]{
\includegraphics[width=0.31\textwidth]{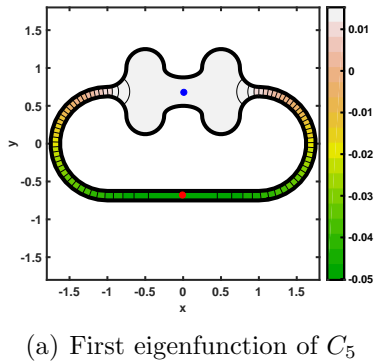}
}
\caption{\label{evcomplex5} The eigenfunction corresponding to the first non-trivial interior Neumann eigenvalue $0.563\,329$ for the domain $C_5$.
}
\end{figure}

Interestingly, we can also remove the bumps in the lower part of the teether domain $C_1$ and still obtain the extreme values inside the new domain $C_5$. The result is shown in Figure \ref{evcomplex5} within $\square_{1.8}$. 
The corresponding 
non-trivial interior Neumann eigenvalue is $0.563\,329$.

\begin{figure}[!ht]  
\subfigure[First eigenfunction of $S$]{
\includegraphics[width=0.31\textwidth]{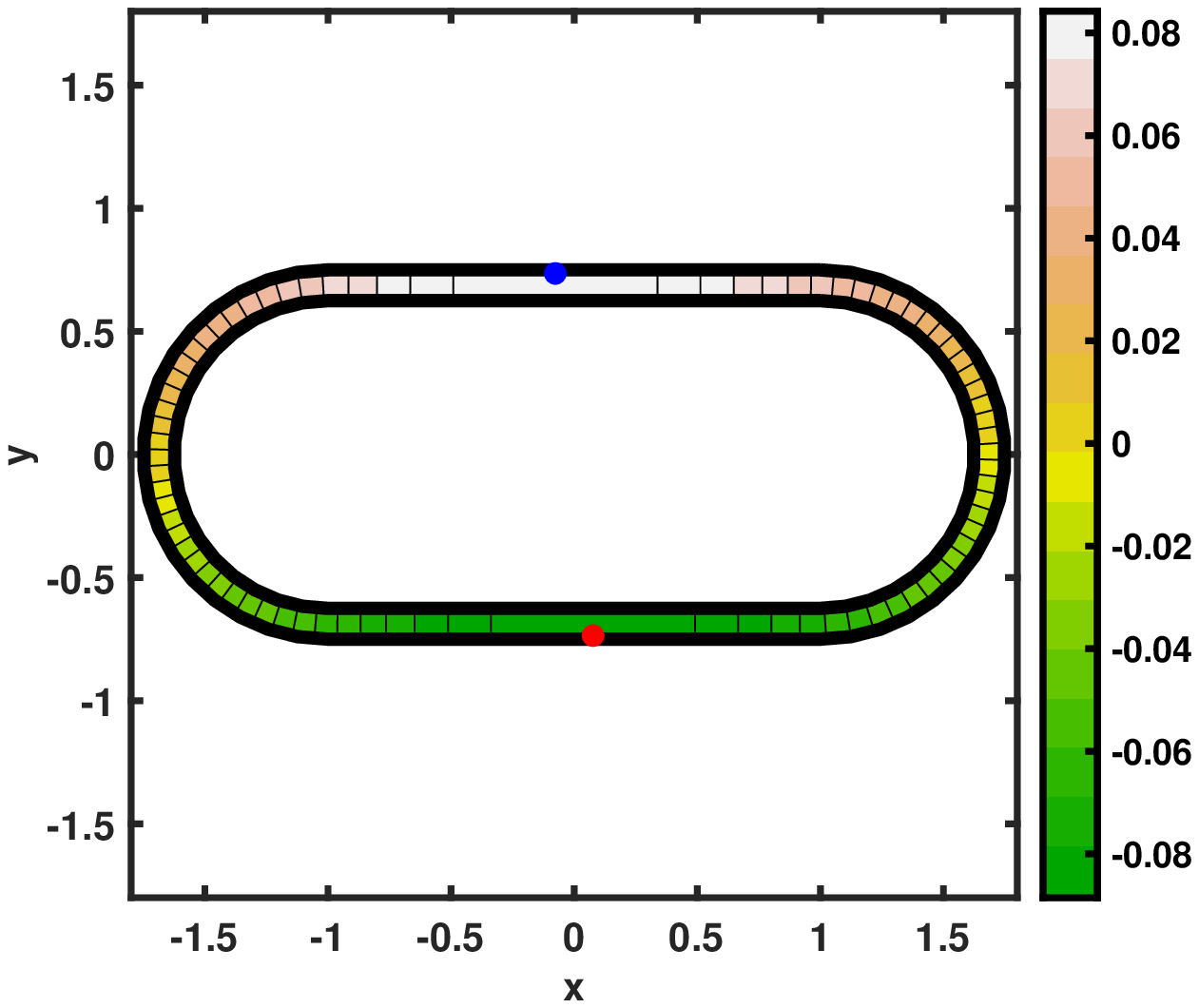}
}
\subfigure[Second eigenfunction of $S$]{
\includegraphics[width=0.31\textwidth]{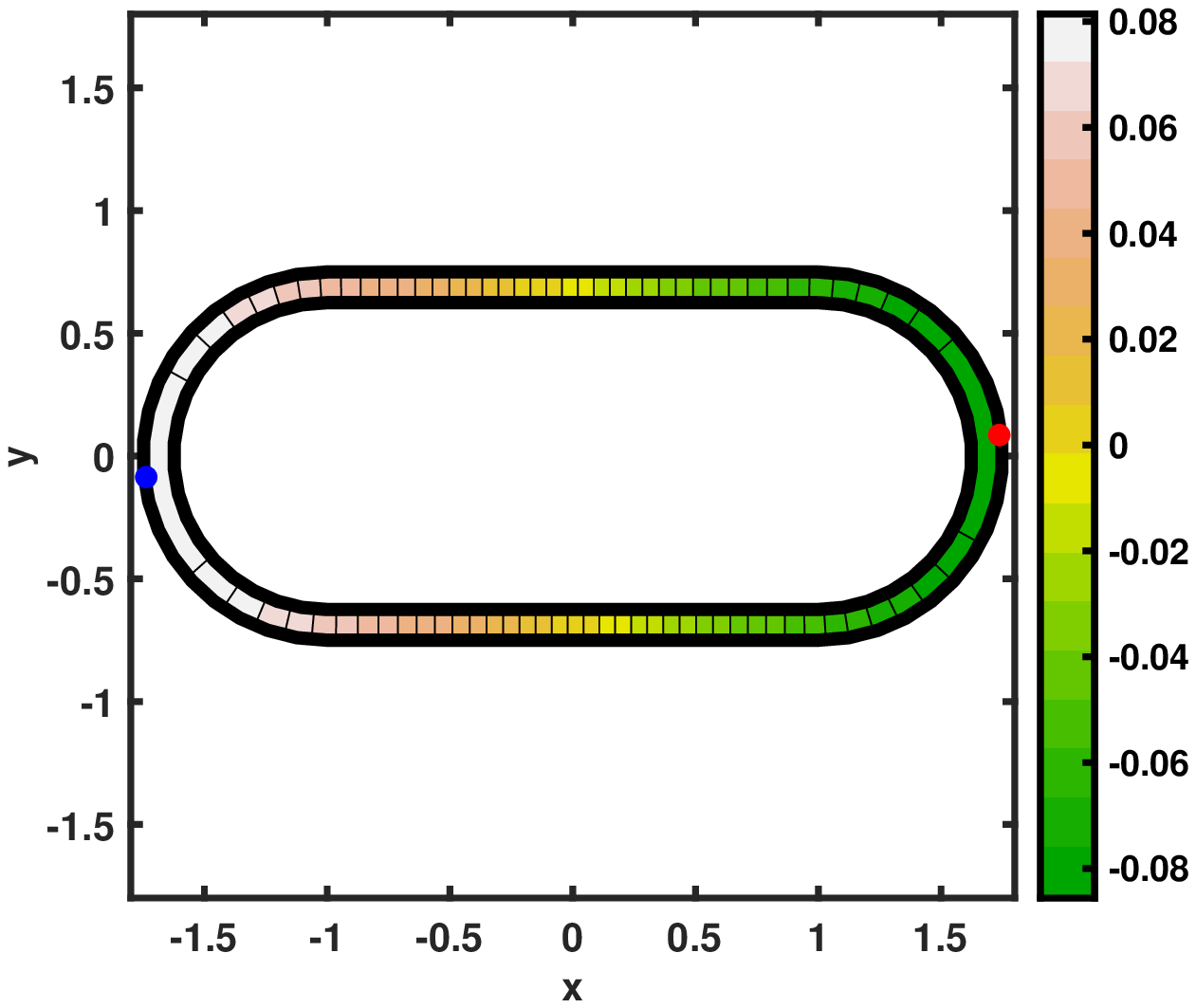}
}
\caption{\label{stadium} The eigenfunctions corresponding to the first and second non-trivial interior Neumann eigenvalue $0.755\,416$ and $0.756\,082$ for the domain $S$.
}
\end{figure}

The idea to also remove the bumps on the upper part yields the following results for the new `stadium' domain $S$ as shown in Figure \ref{stadium} within $\square_{1.8}$. 
Unfortunately, the extreme values are now on the boundary of $S$. The first non-trivial interior Neumann eigenvalues is $0.755\,416$. The second non-trivial eigenvalue is very close (but distinct) $0.756\,082$. 
We can see that the limiting process of making the bumps of $C_4$ smaller (see also $\tilde{C}_4$ and $\hat{C}_4$) yields the eigenfunction corresponding to the second non-trivial interior Neumann eigenfunction for the domain 
$S$. This is very unexpected.

\begin{figure}[!ht]  
\subfigure[First eigenfunction of $F_1$]{
\includegraphics[width=0.31\textwidth]{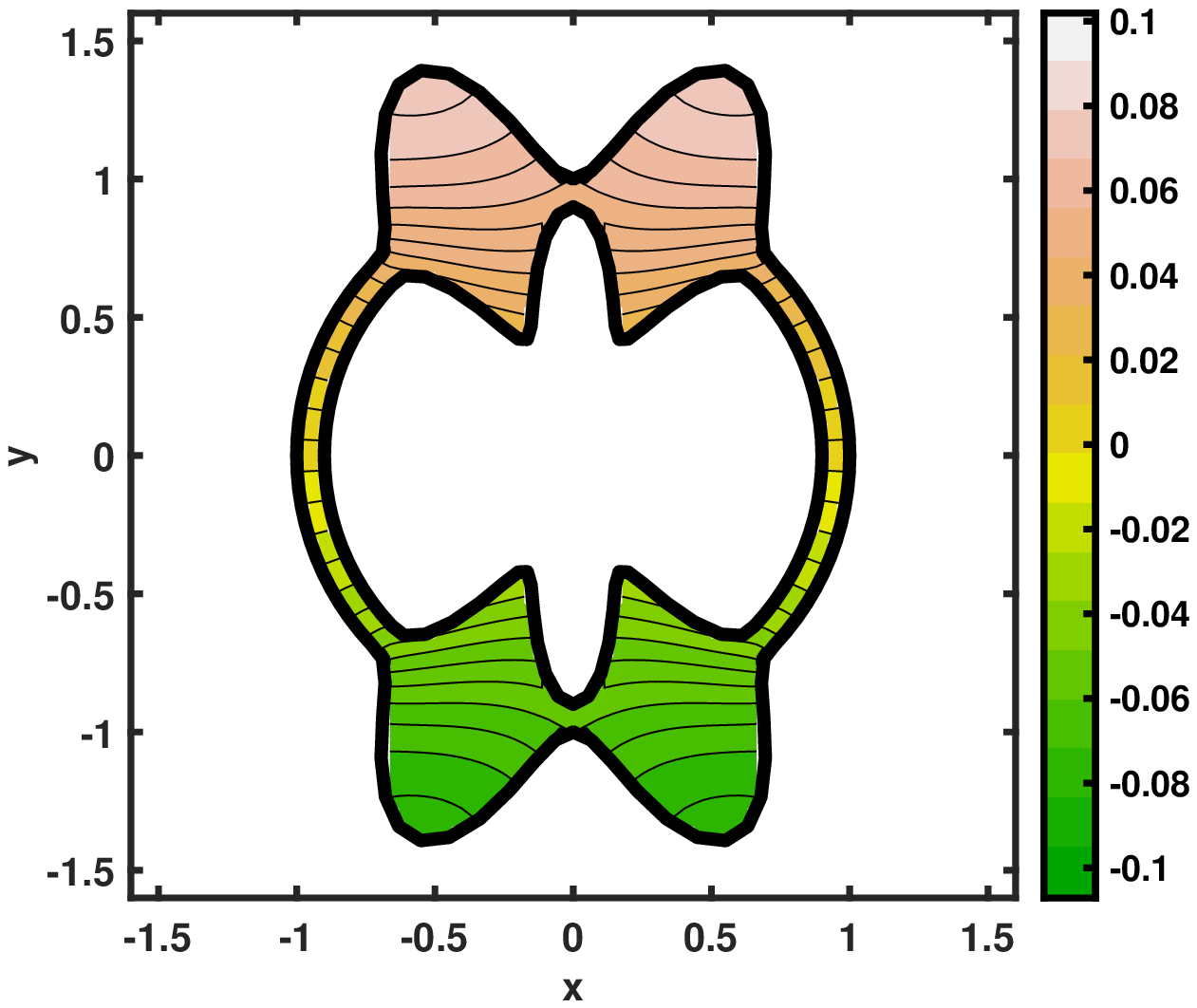}
}
\subfigure[First eigenfunction of $F_2$]{
\includegraphics[width=0.31\textwidth]{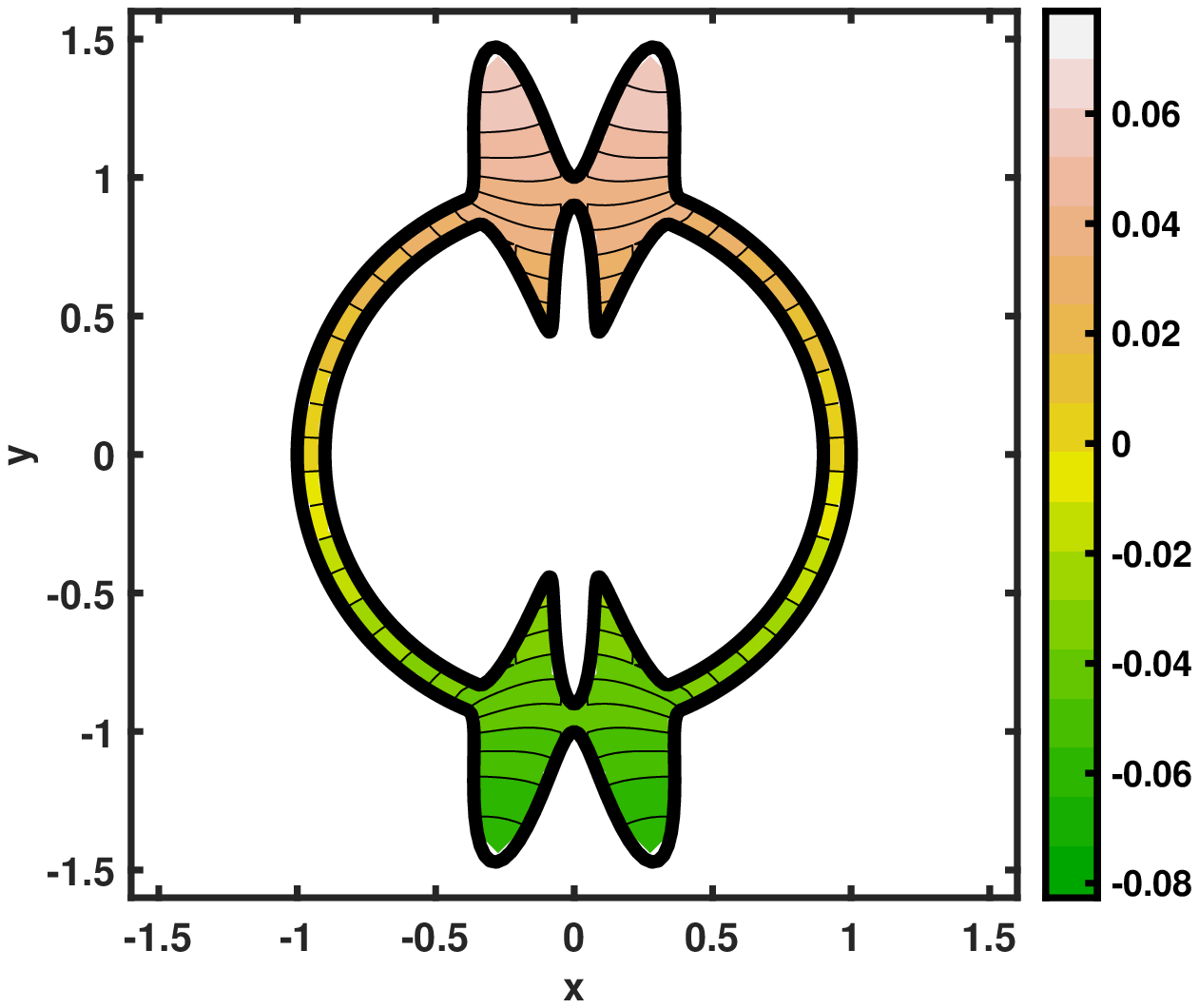} 
}                                                            
\caption{\label{flow} The eigenfunction corresponding to the first non-trivial interior Neumann eigenvalue $1.592\,787$ and $1.717\,098$ for the domains $F_1$ and $F_2$, respectively.
}
\end{figure}

Further, it seems that the heat flow out of the narrow part of the `pipes' has to be almost without inclination.
We use the parametrization $r_i(t)\cdotp (\sin(t),\cos(t))^\top$ with $r_i(t)=a_i(1+1/2\sin(\omega t\cdotp \mathbb{I}_{[(\omega/2-1)\pi/\omega,(\omega/2+1)\pi/\omega]\cup [(3\omega/2-1)\pi/\omega,(3\omega/2+1)\pi/\omega]}))$, $t\in [0,2\pi)$ with $a_1=1$ for the outer boundary and $a_2=0.9$ for the inner boundary. 
Here $\mathbb{I}$ denotes the indicator function. We use $\omega=4$ and $\omega=8$ to construct the domains $F_1$ and $F_2$. 
We obtain the first non-trivial interior Neumann eigenvalue $1.592\,787$ and $1.717\,098$ and the eigenfunction within $\square_{1.6}$. As we can see in Figure \ref{flow} the 
maximum and minimum value are obtained on the boundary. In fact, we have two maxima and two minima.

Next, we show in Table \ref{summ} for the different examples that show a failure of the hot spots conjecture the following information: 
The domain under consideration, the location of the global maximum and minimum inside the domain, and the ratio $\aleph_\textit{max}$ and $\aleph_\textit{min}$ which are defined as the ratio 
of the maximum inside the domain divided by the the maximum on the boundary and likewise for the minimum. All ratios will be larger than one, but we are interested in how close they are to one. Recall that the hot spots conjecture 
fails if we have $\aleph_\textit{max}>1$ and/or $\aleph_\textit{min}>1$. The maximum and minimum are calculated through the Matlab minimization routine fminsearch by passing the function given in (\ref{greensrep2}) to a 
tolerance of $10^{-10}$ for the step size and the difference of a function evaluation. For the first two domains, we use as starting value $(0,0.6)^\top$ and $(0,-0.6)^\top$, respectively. 
For the next three domains, we use $(0,0.69)^\top$ and $(0,-0.69)^\top$ and for the last domain, we use $(0,0.6875)^\top$ and $(0,-0.6875)^\top$.
\begin{table}[!ht]  
\caption{\label{summ} Location of the maximum and minimum inside the domain along with the ratios $\aleph_\textit{max}>1$ and/or $\aleph_\textit{min}>1$ for various domains $D$ that fail the hot spots conjecture.}
\begin{indented}
 \item[]\begin{tabular}{@{}rrrrr}
 \br
  $D$  & location max & location min & $\aleph_\textit{max}$ & $\aleph_\textit{min}$ \\
  \mr
    $C_1$         & $(-3.805_{-8},6.877_{-1})^\top$ & $( 3.299_{-8}, -6.877_{-1})^\top$ & $1+1.221_{-4}$ & $1+1.221_{-4}$\\
    $C_4$         & $( 8.126_{-8},6.875_{-1})^\top$ & $(-3.893_{-8}, -6.875_{-1})^\top$ & $1+1.196_{-5}$ & $1+1.196_{-5}$\\
    $\tilde{C}_4$ & $(-2.998_{-8},6.875_{-1})^\top$ & $(-2.028_{-8}, -6.875_{-1})^\top$ & $1+7.130_{-6}$ & $1+7.118_{-6}$\\
    $\hat{C}_4$   & $( 3.530_{-8},6.875_{-1})^\top$ & $( 2.824_{-8}, -6.875_{-1})^\top$ & $1+3.808_{-6}$ & $1+3.802_{-6}$\\
    $\bar{C}_4$   & $(-4.580_{-8},6.875_{-1})^\top$ & $(-4.787_{-8}, -6.875_{-1})^\top$ & $1+2.137_{-6}$ & $1+2.138_{-6}$\\ 
    $C_5$         & $(-1.968_{-7},6.877_{-1})^\top$ & ---------------                   & $1+1.438_{-3}$ & ---------------\\  
  \br
 \end{tabular}
 \end{indented}
\end{table}

As we can see the maximum and minimum are approximately located on the $y$-axis each centered between the bumps for the first four domains. Interestingly, we obtain for the ratios $1+\epsilon$ with very small $\epsilon>0$. The parameter 
$\epsilon$ decreases as the shape gets closer to the 'stadium' domain. Hence, we might conjecture that there exists a domain with one hole where we have 
$\aleph_\textit{max}=1+\epsilon$ and $\aleph_\textit{min}=1+\epsilon$ with $\epsilon\geq \epsilon_0>0$ given $\epsilon_0$ arbitrarily small.
For the last domain, we obtain the largest value for $\epsilon$. Precisely, we obtain $\epsilon_0>1+10^{-3}$. However, the location of the minimum might be on the boundary (or on the complete part of the $y$-axis).

\begin{figure}[!ht]  
\subfigure[First eigenfunction of $B$]{
\includegraphics[width=0.31\textwidth]{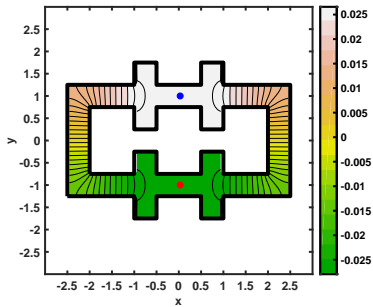}   
}
\caption{\label{brick} The eigenfunction corresponding to the first non-trivial interior Neumann eigenvalue $0.411\,448$ for 
the brick domains $B$.
}
\end{figure}

Finally, we also show a different easy to construct domain with one hole that fails the hot spots conjecture. The outer boundary is given by straight lines connecting the $21$ points 
$(0,1)^\top$, $(3,1)^\top$, $(3,0)^\top$, $(4,0)^\top$, $(4,1)^\top$,
$(6,1)^\top$, $(6,0)^\top$, $(7,0)^\top$, $(7,1)^\top$, $(10,1)^\top$,
$(10,6)^\top$, $(7,6)^\top$, $(7,7)^\top$, $(6,7)^\top$, $(6,6)^\top$,
$(4,6)^\top$, $(4,7)^\top$, $(3,7)^\top$, $(3,6)^\top$, $(0,6)^\top$, and $(0,1)^\top$. The inner boundary is given by the straight lines connecting the $21$ points
$(1,2)^\top$, $(3,2)^\top$, $(3,3)^\top$, $(4,3)^\top$, $(4,2)^\top$,
$(6,2)^\top$, $(6,3)^\top$, $(7,3)^\top$, $(7,2)^\top$, $(9,2)^\top$,
$(9,5)^\top$, $(7,5)^\top$, $(7,4)^\top$, $(6,4)^\top$, $(6,5)^\top$,
$(4,5)^\top$, $(4,4)^\top$, $(3,4)^\top$, $(3,5)^\top$, $(1,5)^\top$, and $(1,2)^\top$ and then stored in reversed order. The resulting coordinates are first shifted by $-(5,7/2)^\top$ and then scaled by $1/2$ 
to center the `brick' domain $B$ with respect to the origin. Using $600$ collocation 
points with the parameters $N=24$, $R=1/2$, $\mu=7/10$, and $\ell=20$ yields the eigenfunction shown in Figure \ref{brick} within $\square_3$ using a resolution of $100\times 100$.
As we can see, the maximal and minimal value are attained within the domain. Hence, we constructed another domain with one hole that fails the hot spots conjecture.

\subsection{Domains with more than one hole}
\begin{figure}[!ht]  
\subfigure[First eigenfunction of $C_{1,0.1}$]{
\includegraphics[width=0.31\textwidth]{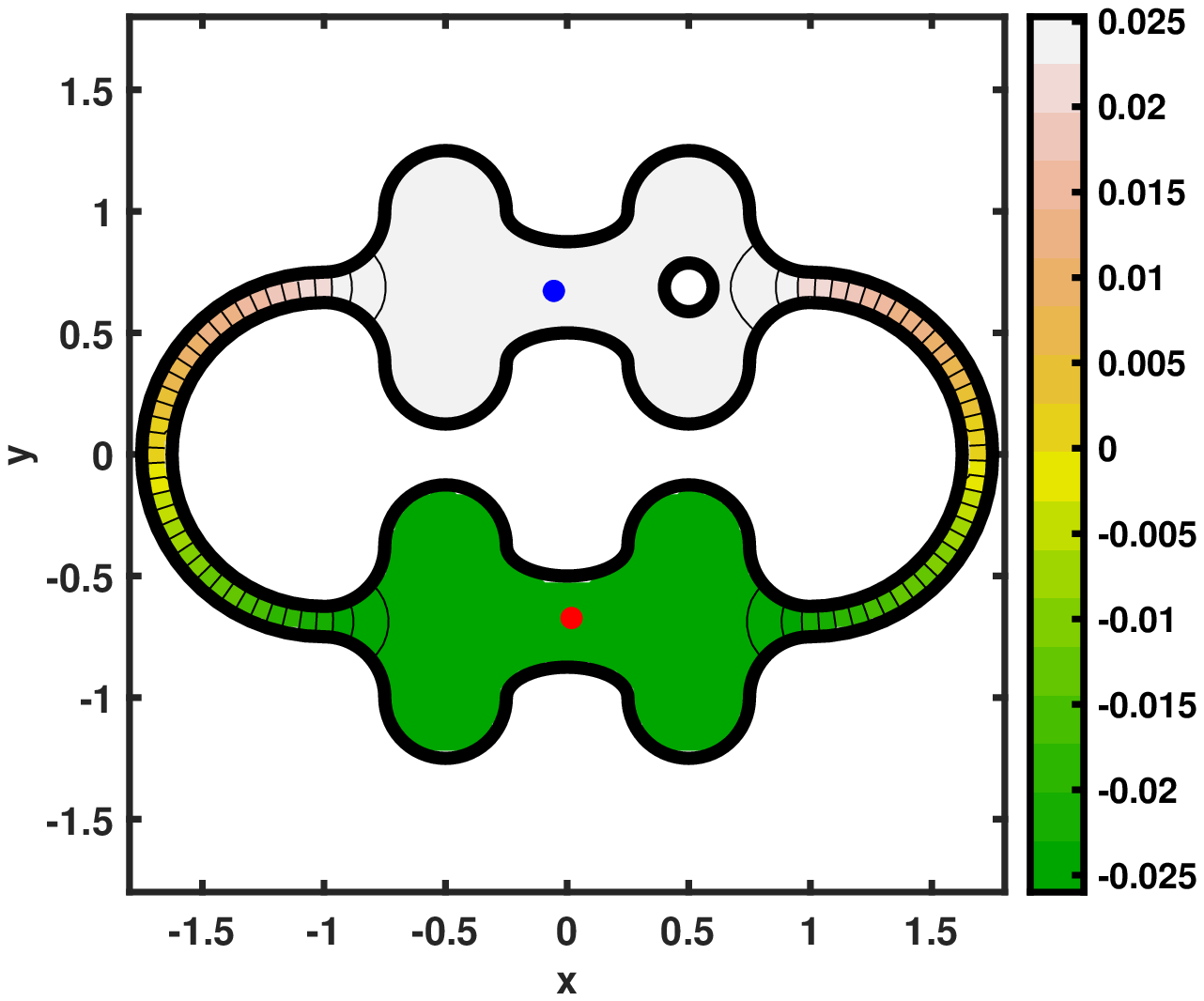}   
}
\subfigure[First eigenfunction of $C_{1,0.15}$]{
\includegraphics[width=0.31\textwidth]{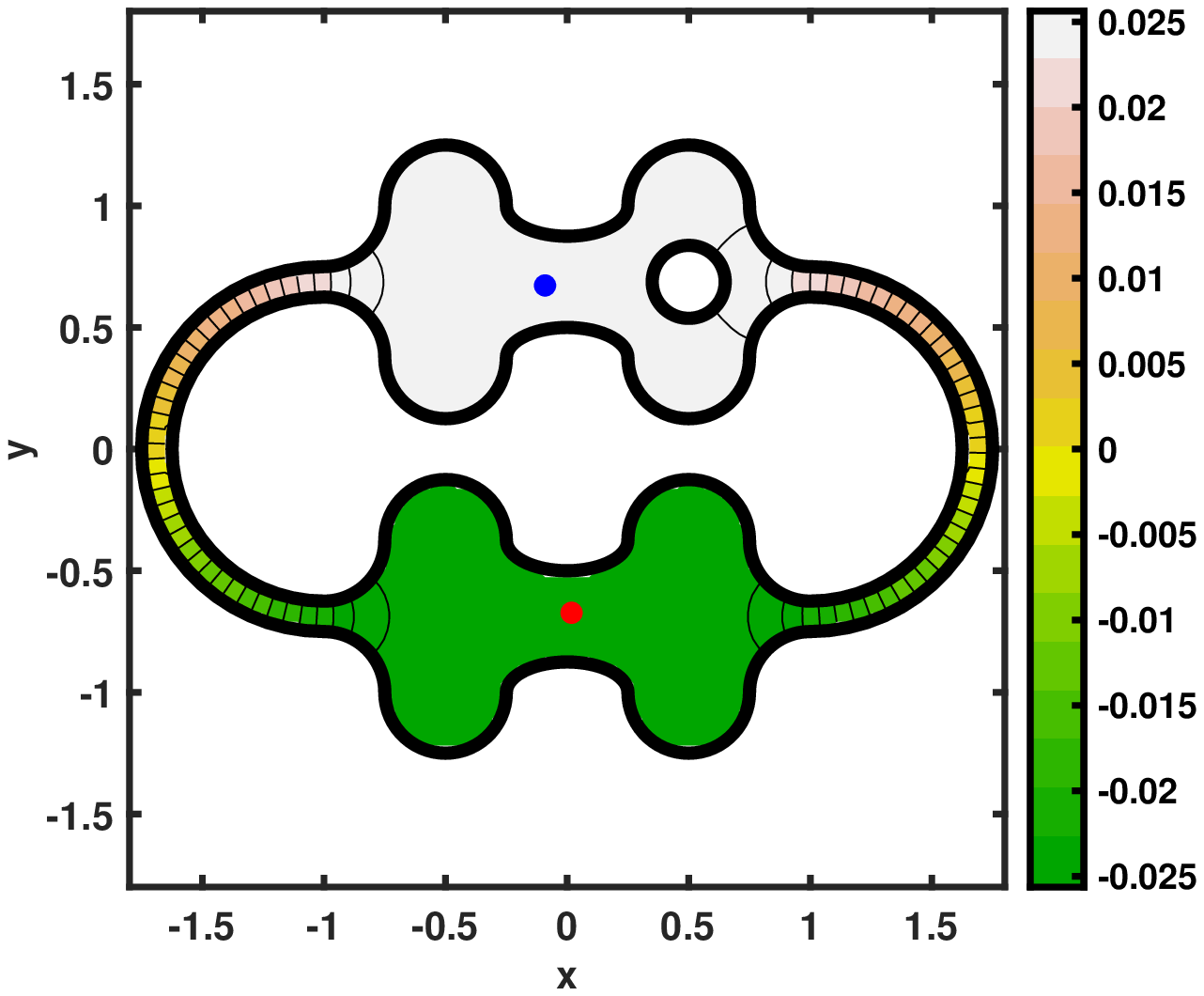}  
}
\subfigure[First eigenfunction of $C_{1,0.2}$]{
\includegraphics[width=0.31\textwidth]{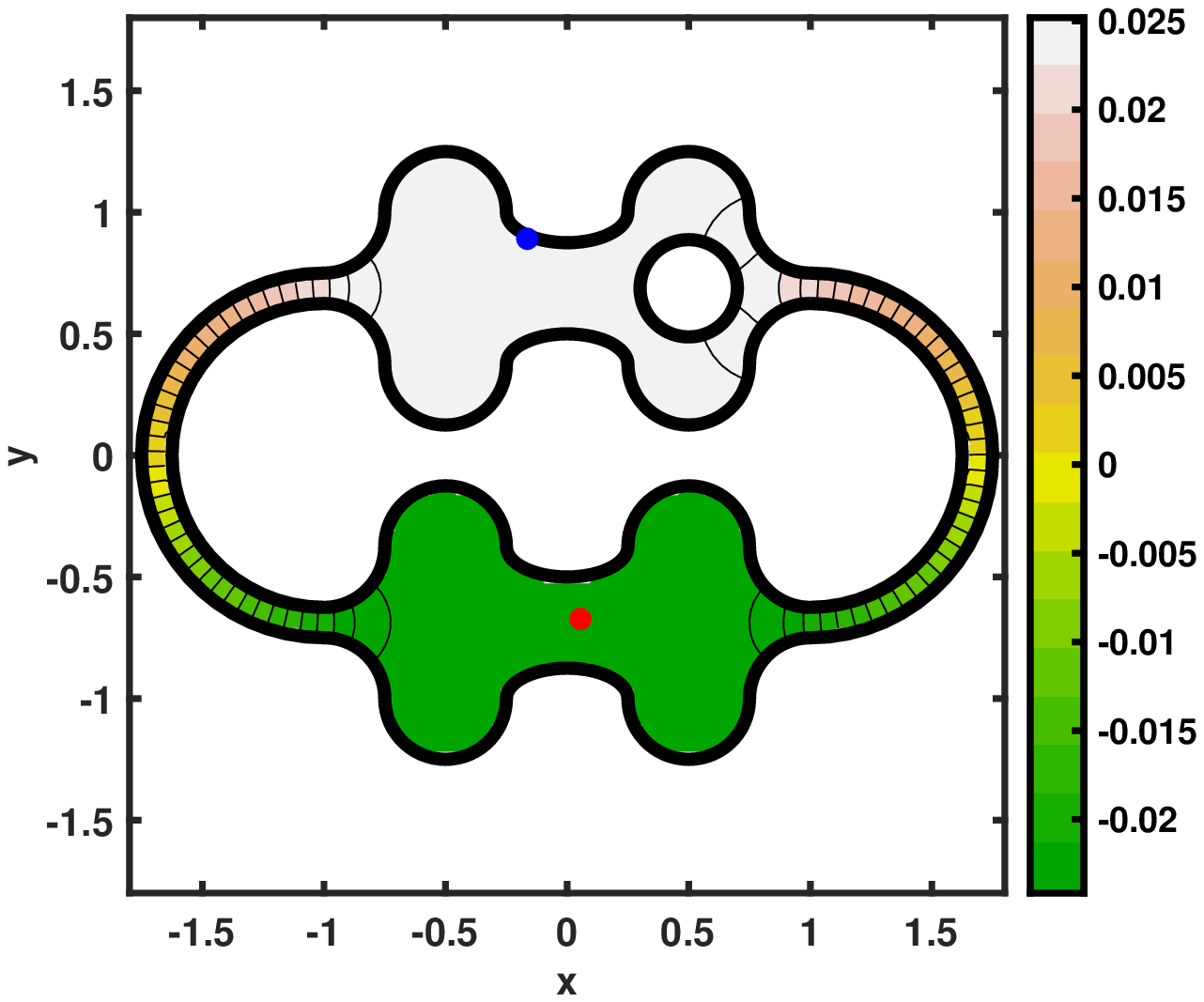}  
}
\caption{\label{holes2} The eigenfunction corresponding to the first non-trivial interior Neumann eigenvalue $0.372\,580$, $0.374\,917$, and $0.378\,131$ for 
the domains $C_{1,0.1}$, $C_{1,0.15}$ and $C_{1,0.2}$, respectively.
}
\end{figure}

\begin{figure}[!ht]  
\subfigure[First eigenfunction of $\widetilde{C}_{1,0.1}$]{
\includegraphics[width=0.31\textwidth]{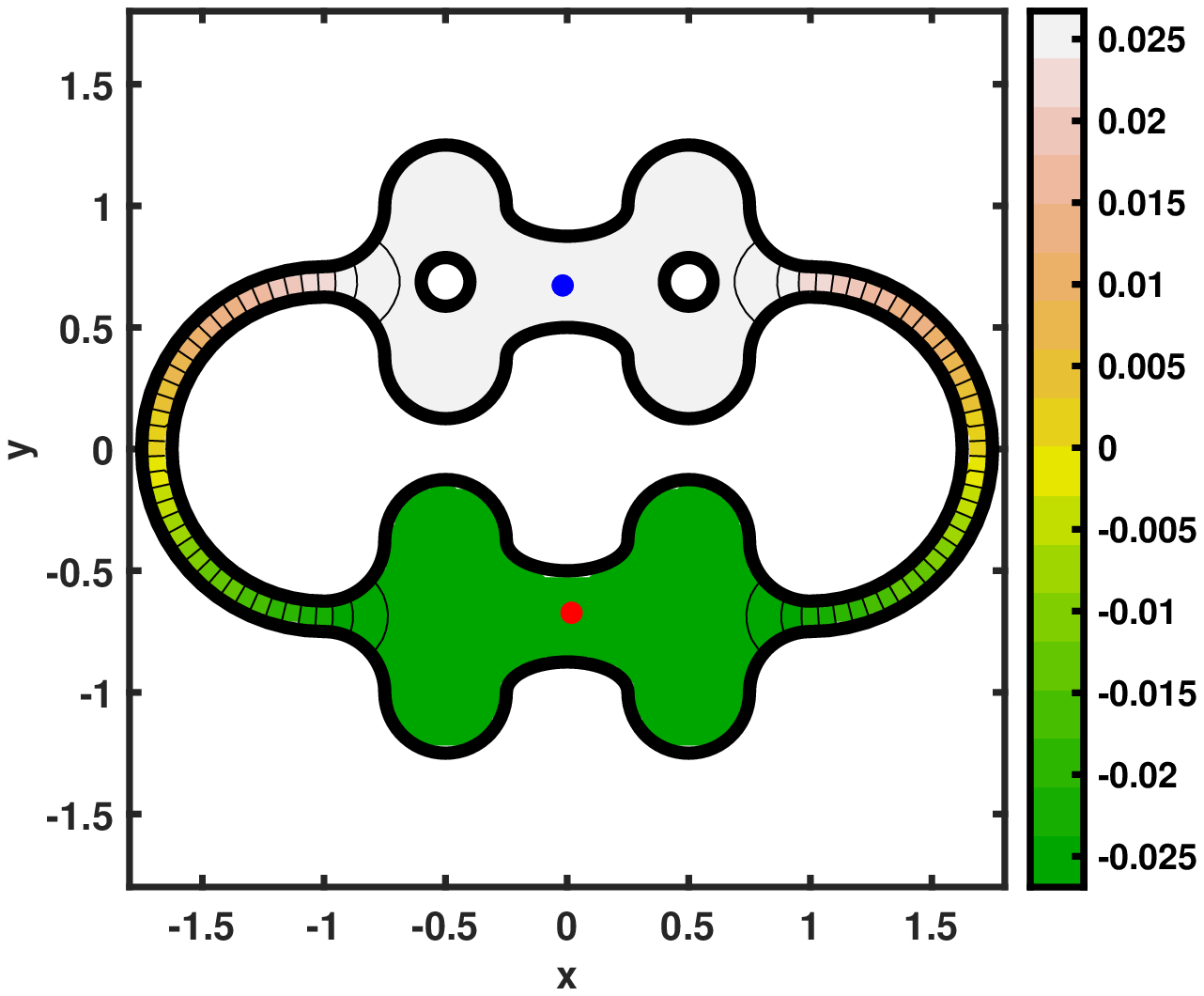}   
}
\subfigure[First eigenfunction of $\widetilde{C}_{1,0.15}$]{
\includegraphics[width=0.31\textwidth]{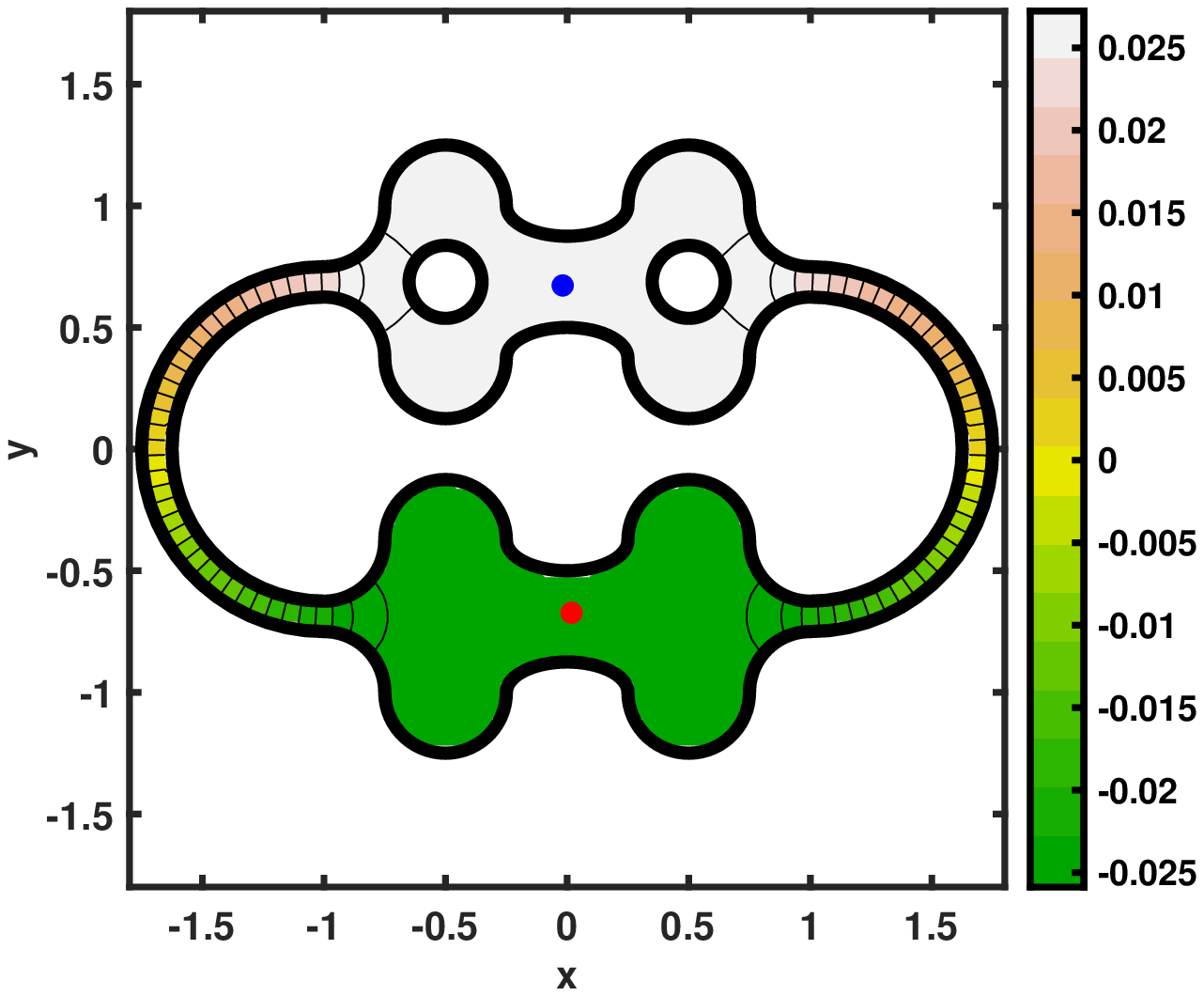}  
}
\subfigure[First eigenfunction of $\widetilde{C}_{1,0.2}$]{
\includegraphics[width=0.31\textwidth]{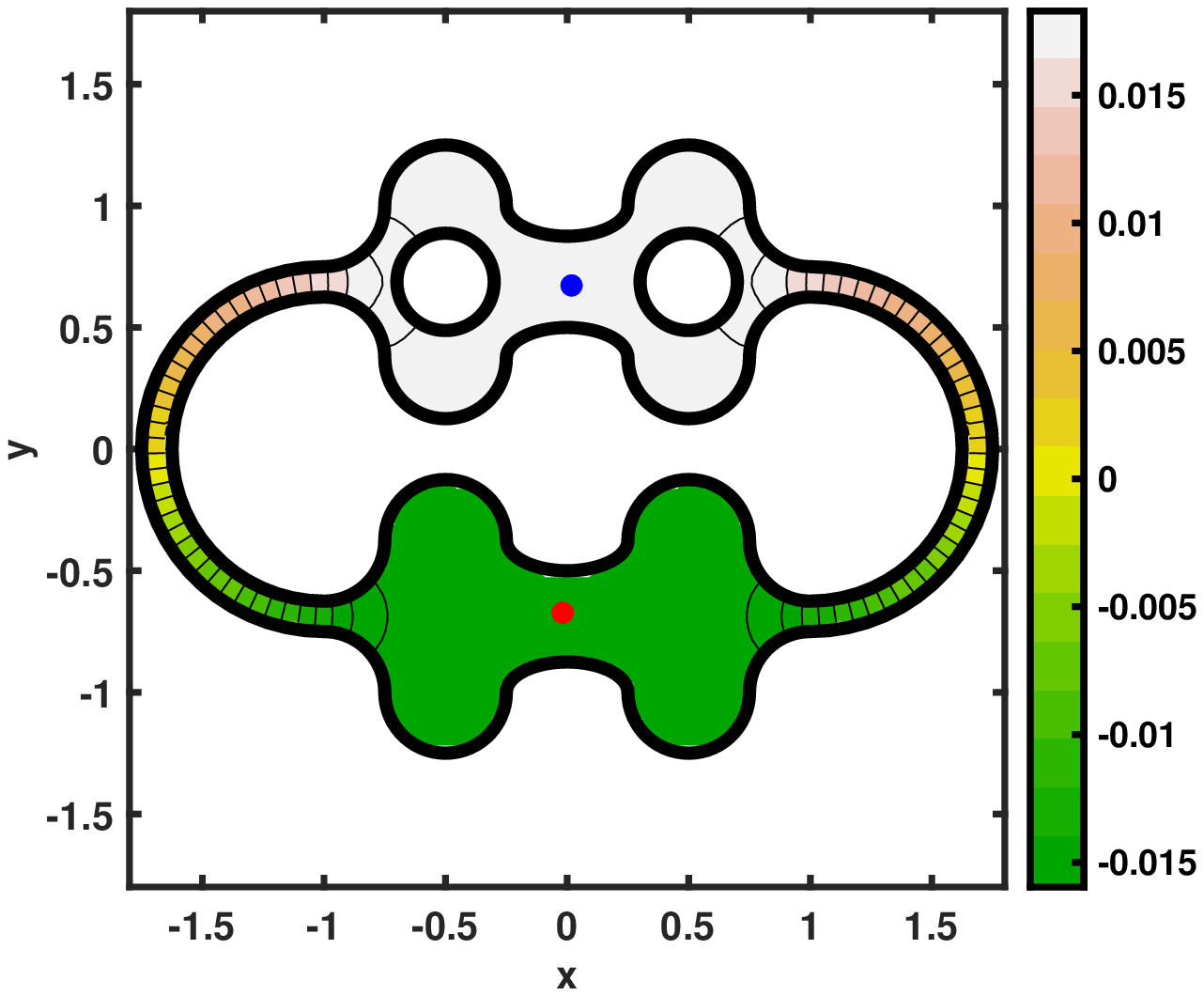}  
}
\caption{\label{holes3} The eigenfunction corresponding to the first non-trivial interior Neumann eigenvalue $0.374\,552$, $0.379\,670$, and $0.387\,461$ for 
the domains $\widetilde{C}_{1,0.1}$, $\widetilde{C}_{1,0.15}$ and $\widetilde{C}_{1,0.2}$, respectively.
}
\end{figure}

\begin{figure}[!ht]  
\subfigure[First eigenfunction of $C_{1,0.1}^{\mathrm{mir}}$]{
\includegraphics[width=0.31\textwidth]{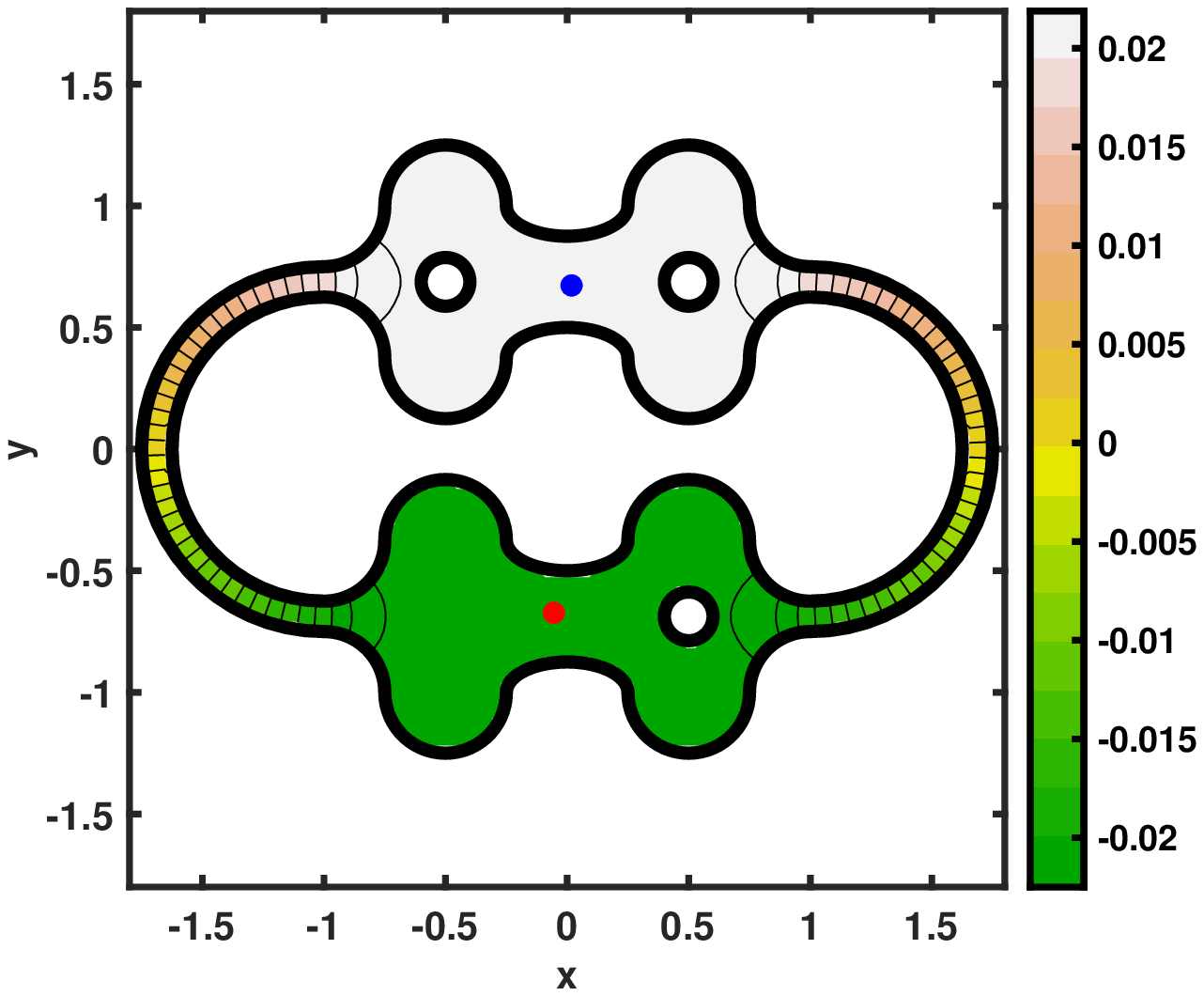}  
}
\subfigure[First eigenfunction of $C_{1,0.15}^{\mathrm{mir}}$]{
\includegraphics[width=0.31\textwidth]{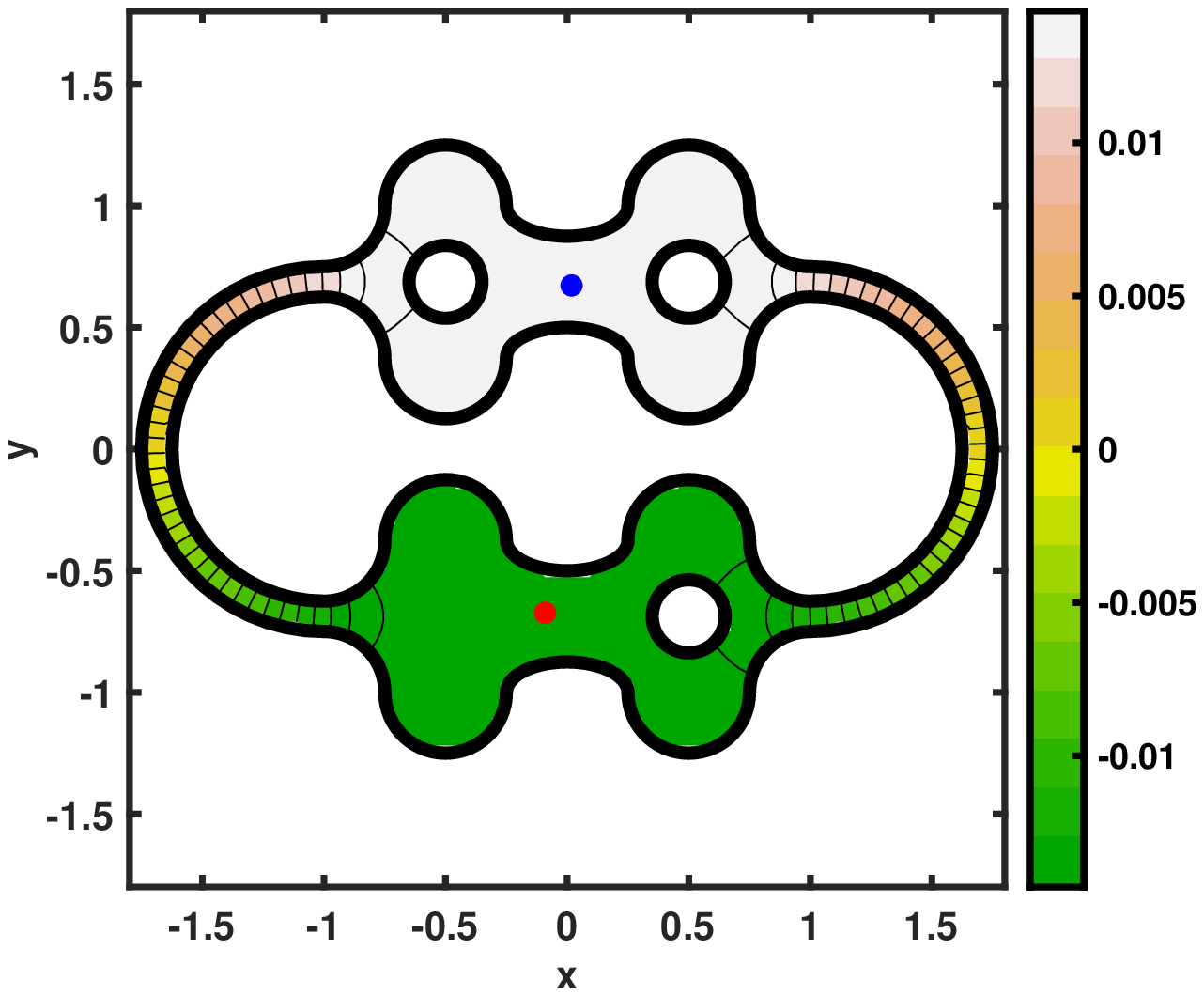}  
}
\subfigure[First eigenfunction of $C_{1,0.2}^{\mathrm{mir}}$]{
\includegraphics[width=0.31\textwidth]{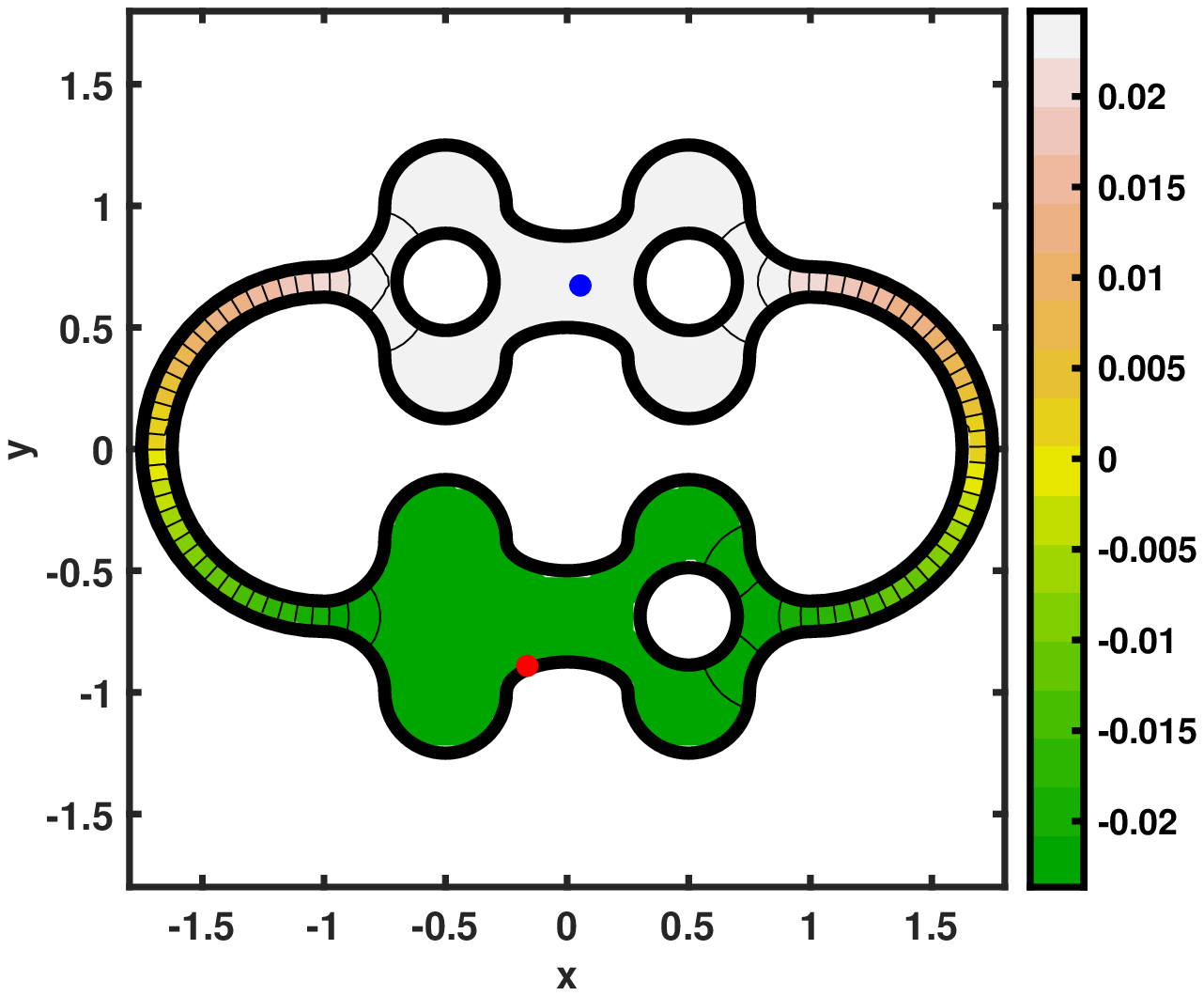}  
}
\caption{\label{holes4} The eigenfunction corresponding to the first non-trivial interior Neumann eigenvalue $0.376\,406$, $0.383\,781$, and $0.394\,522$ for 
the domains $C_{1,0.1}^{\mathrm{mir}}$, $C_{1,0.15}^{\mathrm{mir}}$ and $C_{1,0.2}^{\mathrm{mir}}$, respectively.
}
\end{figure}

\begin{figure}[!ht]  
\subfigure[First eigenfunction of $\widetilde{C}_{1,0.1}^{\mathrm{mir}}$]{
\includegraphics[width=0.31\textwidth]{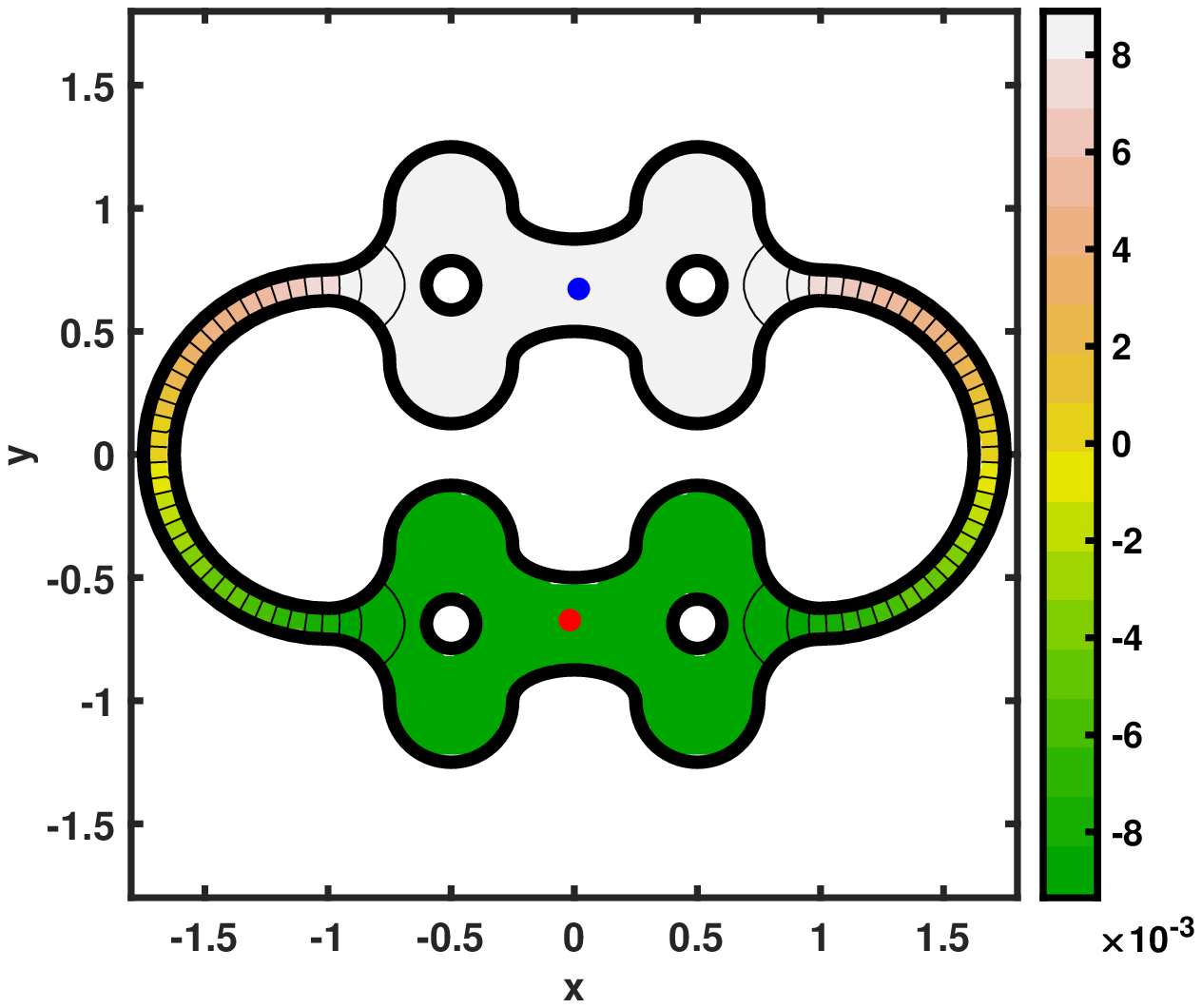}  
}
\subfigure[First eigenfunction of $\widetilde{C}_{1,0.15}^{\mathrm{mir}}$]{
\includegraphics[width=0.31\textwidth]{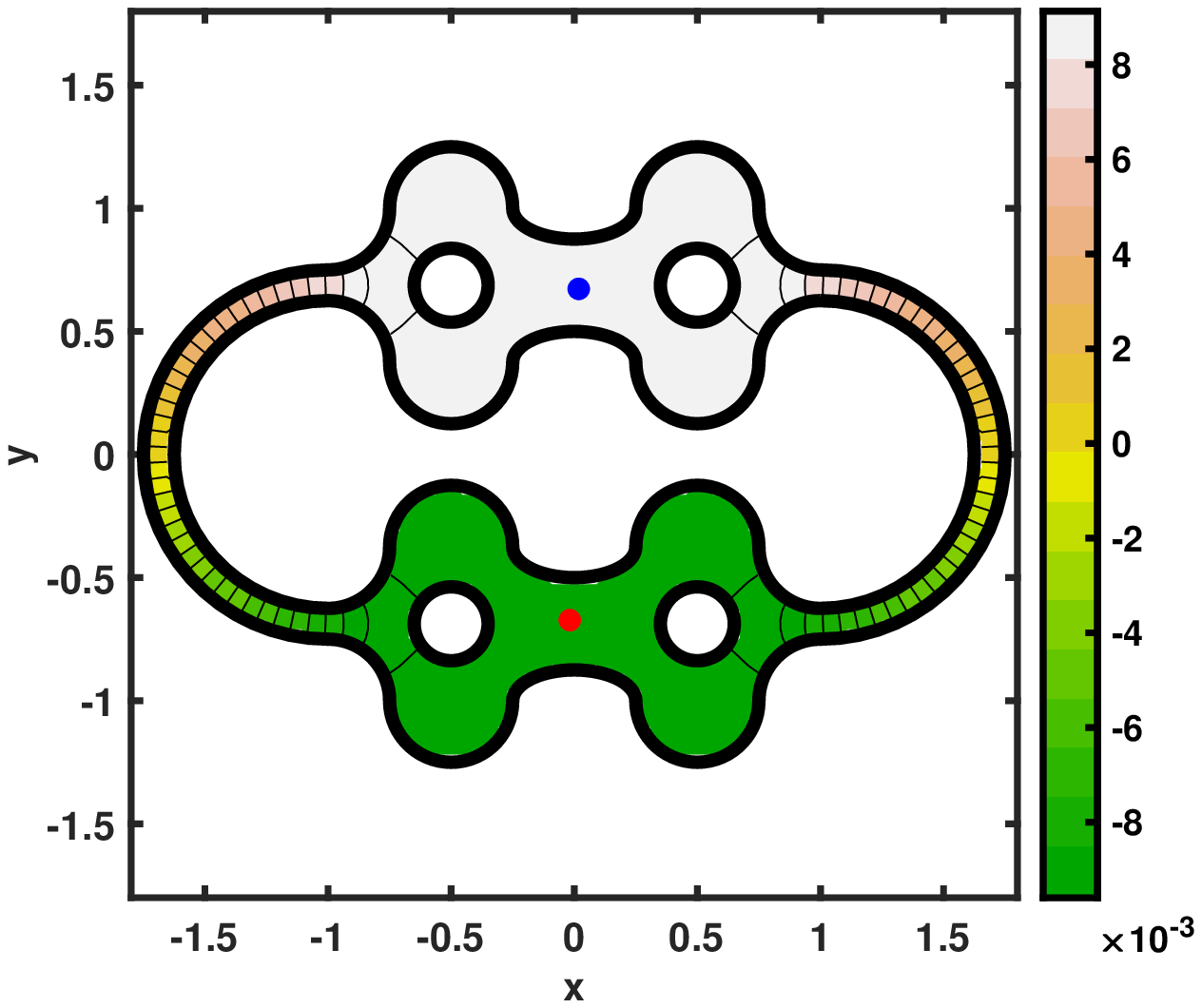}  
}
\subfigure[First eigenfunction of $\widetilde{C}_{1,0.2}^{\mathrm{mir}}$]{
\includegraphics[width=0.31\textwidth]{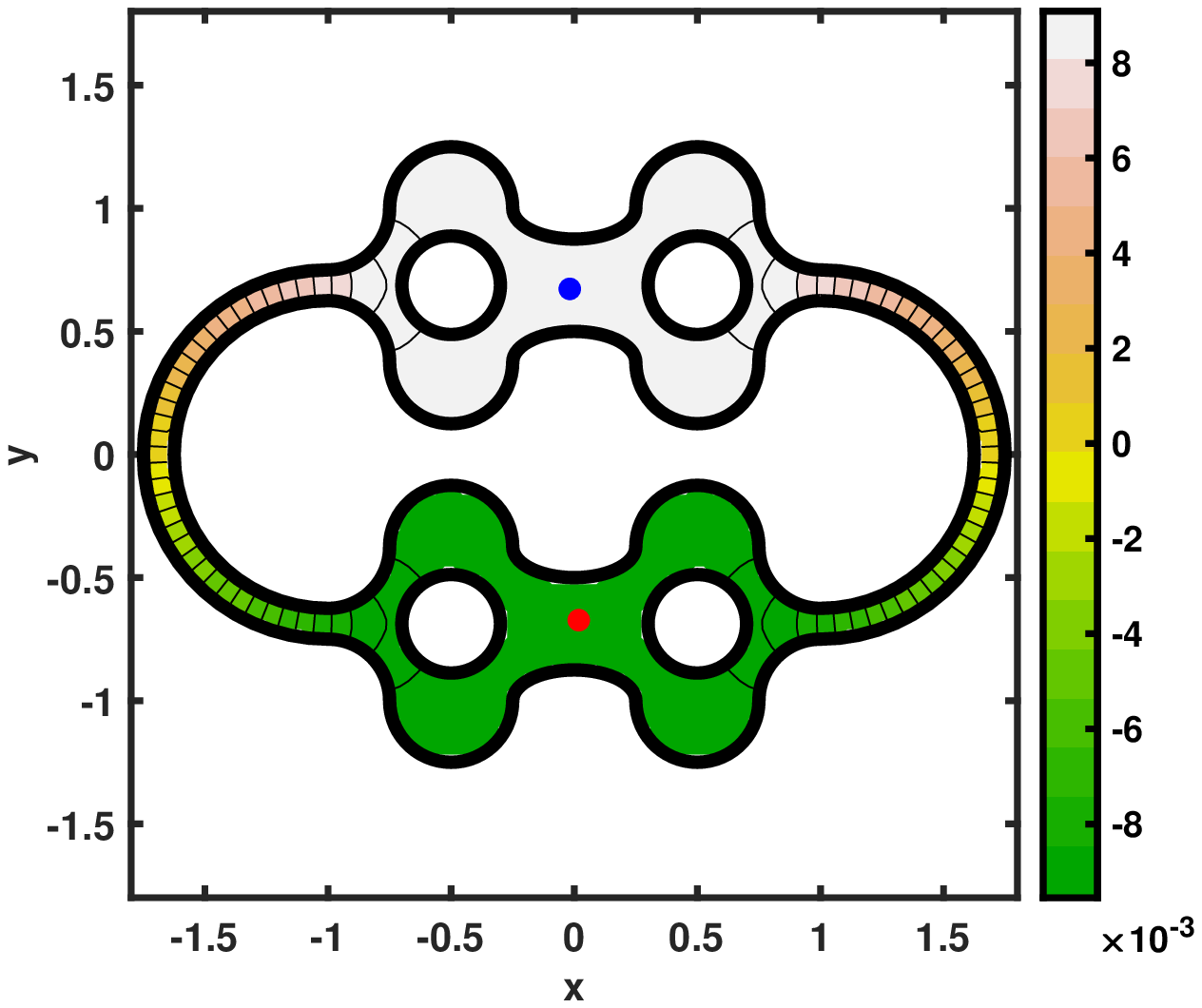}  
}
\caption{\label{holes5} The eigenfunction corresponding to the first non-trivial interior Neumann eigenvalue $0.378\,360$, $0.388\,438$, and $0.403\,504$ for 
the domains $\widetilde{C}_{1,0.1}^{\mathrm{mir}}$, $\widetilde{C}_{1,0.15}^{\mathrm{mir}}$ and $\widetilde{C}_{1,0.2}^{\mathrm{mir}}$, respectively.
}
\end{figure}

Finally, we show without further discussion that we are also able to construct examples with more than one hole where the hot spots conjecture fails to hold. Using the teether domain $C_1$ and removing a circle 
centered at $(1/2,11/16)^\top $ with radius $\mathfrak{R}$, yields a domain with two holes, say $C_{1,\mathfrak{R}}$. Using $\mathfrak{R}=0.1$, $\mathfrak{R}=0.15$, and $\mathfrak{R}=0.2$, gives the 
results shown in Figure \ref{holes2} where we used the same 
set of parameters as for the results for $C_1$. 
As we can see, the maximum and minimum still remain inside the domains $C_{1,0.1}$ and $C_{1,0.15}$. 
The location of the maximum is slightly shifted to the right whereas the minimum is slightly shifted to the left for the first two cases. If the radius of the removed circle is large, 
then the maximum goes to the boundary whereas the 
minimum stays inside the domain $C_{1,0.2}$ as shown in the last contour plot of Figure \ref{holes2}.
Next, we construct domains with three holes. Therefore, we use the previous domain $C_{1,\mathfrak{R}}$ and mirror the circular hole at the $y$-axis. This yields the domain $\widetilde{C}_{1,\mathfrak{R}}$. Using the same parameters as 
before yields the results shown in Figure \ref{holes3}.
As we can see now, the maximum and minimum remain inside all the three domains $\widetilde{C}_{1,0.1}$, $\widetilde{C}_{1,0.15}$, and $\widetilde{C}_{1,0.2}$. Next, we consider domains with four holes. Therefore, 
we use the domains 
$C_{1,0.1}$, $C_{1,0.15}$ and $C_{1,0.2}$ and mirror the upper right circular hole at the $x$-axis. This yields the domains $C_{1,0.1}^{\mathrm{mir}}$, $C_{1,0.15}^{\mathrm{mir}}$ and $C_{1,0.2}^{\mathrm{mir}}$, respectively. The results are presented in
Figure \ref{holes4}.
As expected, we obtain the minimal and maximal value inside the domain for the two domains $\widetilde{C}_{1,0.1}^{\mathrm{mir}}$ and $\widetilde{C}_{1,0.15}^{\mathrm{mir}}$ whereas the maximum is on the 
boundary for the domain $\widetilde{C}_{1,0.2}^{\mathrm{mir}}$. The domains with five holes $\widetilde{C}_{1,0.1}^{\mathrm{mir}}$, $\widetilde{C}_{1,0.15}^{\mathrm{mir}}$ and $\widetilde{C}_{1,0.2}^{\mathrm{mir}}$ 
are constructed by mirroring the domains $\widetilde{C}_{1,0.1}$, $\widetilde{C}_{1,0.15}$ and $\widetilde{C}_{1,0.2}$ at the $x$-axis. The results are shown in Figure \ref{holes5}.
As we can see, the maximal and minimal values are attained inside all the considered domains with five holes.

The extension for the construction of domains having more than five holes which do not satisfy the hot spots conjecture is now straightforward.

\section{Summary and outlook}\label{suma}
In this paper, a detailed description is given on how to 
compute the first non-trivial eigenvalue and its corresponding eigenfunction for the Laplace equation with Neumann boundary condition for a given domain with one hole. The problem is reformulated as a non-linear eigenvalue problem involving boundary integral equations thus reducing a two-dimensional problem to a one-dimensional problem. Due to superconvergence we are able to achieve highly accurate approximations both for the eigenvalue and the eigenfunction. With this method at hand, we can compute the eigenvalue and eigenfunction for several different constructed domains. This gives the possibility to find domains with one hole failing the hot spots conjecture and investigate the influence of varying the domain. 
Some interesting observation can be made such as that the ratio between the maximal/minimal value inside the domain and the maximal/minimal value on the boundary can be $1+10^{-3}$.
The Matlab codes including the produced data are available at github $$\texttt{https://github.com/kleefeld80/hotspots}$$
and researchers can run it on their own constructed domains and reproduce the numerical results within this article. 
This might give new ideas whether one can find assumptions in order to prove or disprove the hot spots conjecture.
The extension for domains with more than one hole is straightforward. For the sake of completeness they are given at the end of the numerical results section for domains with up to five holes, but without detailed discussion. 

It would be interesting to check whether it is possible to construct three-dimensional domains with one hole that fail the hot spots conjecture, too. 
The software for a domain without a hole would already be available and only needs to be extended (see \cite{kleefeldlin2,kleefeldhabil}). The consideration of other partial differential equations in two or three dimensions whose fundamental solution is known together with Neumann boundary condition could be numerically investigated as well.

\ack
The author thanks Prof. Stefan Steinerberger from the University of Washington, Seattle (USA) for the fruitful discussions during the preparation of the manuscript.

\section*{References}	
\bibliographystyle{abbrv}
\bibliography{ip-biblio}

\end{document}